\documentclass[oneside,reqno, twoside, 11pt]{amsart}  
\usepackage[pdftex]{graphicx}
\usepackage{color}
\usepackage[T1]{fontenc}  
\usepackage{lmodern}
\usepackage[english]{babel}
\usepackage{mathrsfs}
\usepackage{appendix}  

\usepackage{upgreek}   

\usepackage{amsfonts}
\usepackage{extpfeil}
\usepackage{verbatim}

\usepackage{BOONDOX-cal}

\usepackage{microtype}
\usepackage[new]{old-arrows}
\usepackage{amsmath}
\usepackage[all, cmtip]{xy}
\usepackage{tikz-cd}

\usepackage{amsthm}
\usepackage{mathtools}

\usepackage{bbm}
\usepackage{paralist}
\usepackage[T1]{fontenc}
\usepackage[colorlinks,citecolor=magenta,pagebackref=true,urlcolor=magenta, bookmarksdepth=3]{hyperref}
\usepackage[nodisplayskipstretch]{setspace}
\usepackage{breqn}

\usepackage{nicefrac}

\usepackage[font=small,labelfont=bf]{caption}
\usepackage[labelformat=simple]{subcaption}

\usepackage{epsfig}

\DeclareFontFamily{OT1}{pzc}{}
\DeclareFontShape{OT1}{pzc}{m}{it}{<-> s * [1.10] pzcmi7t}{}
\DeclareMathAlphabet{\mathpzc}{OT1}{pzc}{m}{it}

\makeatletter
\newtheorem*{rep@theorem}{\rep@title}
\newcommand{\newreptheorem}[2]{%
	\newenvironment{rep#1}[1]{%
		\def\rep@title{#2~\ref{##1}}%
		\begin{rep@theorem}}%
		{\end{rep@theorem}}}

\makeatother

\theoremstyle{plain}
\newtheorem*{thm*}{Theorem}
\newtheorem*{cor*}{Corollary}
\newtheorem{thm}{Theorem}[section]

\newtheorem{cor}[thm]{Corollary}
\newtheorem{lem}[thm]{Lemma}
\newtheorem*{lem*}{Lemma}
\newtheorem{prp}[thm]{Proposition}
\newtheorem{conj}[thm]{Conjecture}

\theoremstyle{definition}
\newtheorem{dfn}[thm]{Definition}
\newtheorem{ex}[thm]{Example}

\newtheorem{rem}[thm]{Remark}
\newtheorem{rmk}[thm]{Remark}
\newtheorem{con}[thm]{Convention}

 \newcommand{\Z}{\mathbb{Z}}
 \newcommand{\R}{\mathbb{R}}

\newcommand{\sbseq}{\subseteq}
\newcommand{\spseq}{\supseteq}

\newcommand{\vanish}[1]{}

\def\sbs\subset
\def\sbseq{\subseteq}

\def\langle{\left<}
\def\rangle{\right>}

\def\({\left(}
\def\){\right)}
\def\no={\,{\,|\!\!\!\!\!=\,\,}}

\def\no={\,{\,|\!\!\!\!\!=\,\,}}
\def\sbseq{\subseteq}

\def\sbseq{\subseteq}
\def\sbs\subset
\def\spseq{\supseteq}

\newcommand{\xqedhere}[2]{%
\rlap{\hbox to#1{\hfil\llap{\ensuremath{#2}}}}}

\newcommand\Defn[1]{\textbf{#1}}
\newcommand{\cm}[1]{}

\newcommand\mbf[1]{\mathbf{#1}}

\newcommand\mr[1]{\mathrm{#1}}

\newcommand{\fld}{\mathbbm{k}}

\newcommand\vol{\mr{vol}}
\renewcommand\deg{\mr{deg}}
\newcommand{\bigslant}[2]{{\raisebox{.3em}{$#1$} \Big/ \raisebox{-.3em}{$#2$}}}
\renewcommand\emptyset{\varnothing}
\newcommand\x{\mathbf{x}}

\DeclareMathOperator{\cone}{cone}

\title[Lattice polytopes and Lefschetz properties]{Lattice polytopes and semigroup algebras: Generic Lefschetz properties and Parseval-Rayleigh identities}
\author[Karim Adiprasito]{Karim Alexander Adiprasito}
\address{{Karim Adiprasito \ \emph{and}\ Vasiliki Petrotou}, Sorbonne Université and Université Paris Cité, CNRS, IMJ-PRG, F-75005 Paris, France}
\email{adiprasito@imj-prg.fr \emph{and} petrotou@imj-prg.fr}

\author[Stavros Papadakis]{Stavros Argyrios Papadakis}
\address{{Stavros Papadakis}, Department of Mathematics, University of Ioannina, Ioannina, 45110, Greece}
\email{spapadak@uoi.gr}
\author[Vasiliki Petrotou]{Vasiliki Petrotou}

\date{17.09.2025}

\subjclass[2010]{}

\begin{document}

\begin{abstract}
We study semigroup algebras associated to lattice polytopes. 

\vspace*{2mm}

We begin by generalizing and refining work of Hochster, and describe the volume maps of these algebras, that is, their fundamental classes, in terms of Parseval-Rayleigh identities and differential equations, which we prove to be equivalent. 

\vspace*{2mm}

We use these descriptions to establish strong Lefschetz properties. 

\vspace*{2mm}

A consequence is the resolution of several conjectures concerning unimodality properties of the $h^\ast$-polynomial of lattice polytopes.
\end{abstract}

\maketitle

\newcommand{\AR}{\mathcal{A}}
\newcommand{\IR}{\mathcal{I}}
\newcommand{\JR}{\mathcal{J}}
\newcommand{\BR}{\mathbf{I}}
\newcommand{\CR}{\mathcal{C}}
\newcommand{\Mu}{M}
\newcommand{\Soc}{\mathcal{S}\hspace{-1mm}\mathcal{o}\hspace{-1mm}\mathcal{c}}
\newcommand{\Socl}{{\Soc^\circ}}

\section{Lattice polytopes and semigroup algebras}

The protagonist of this paper is a convex polytope $P$ all whose vertices lie in the lattice~$\mathbb{Z}^d$. These polytopes are called lattice polytopes, and are of tremendous importance throughout mathematics, see \cite{AK, Barvinok, BB, BeckRobins, CoxLittleSchenck, KKMS, PK}. Not getting distracted with an endless list, let us simply note this: One of the main points of focus in this context is often the function \[E_P(i)\ :=\ \#\{iP\cap \mathbb{Z}^d\},\]
firstly at nonnegative integers $i$. It is one of the fundamental facts of Ehrhart theory that this function is a polynomial. 

It is often convenient to encode this polynomial into a generating function
\[\mathrm{Ehr}_P(t)\ :=\ \sum_{i=0}^\infty E_P(i) t^i.\]
It is in particular fruitful to write it this way because another polynomial appears: We can write 
\[\mathrm{Ehr}_P(t)\ =\ \frac{h^\ast_P(t)}{(1-t)^{d+1}},\]
where $d$ is the dimension of $P$ and 
$\ h^\ast_P(t)=h_0^\ast+h_1^\ast t+\ldots+h_d^\ast t^d\ $\ is a polynomial of degree at most $d$. A purely combinatorial perspective would be to only try to understand the coefficients of these polynomials, but that view is incomplete.

This is because $\mathrm{Ehr}_P(t)$ and  $h^\ast_P(t)$ have an algebraic interpretation, and it is really this algebraic structure that is a greater mystery, and that has attracted almost as much attention as the numerical question about understanding the aforementioned functions. The story goes like this:

Embed the polytope $P$ in $\R^{d}\times \R$ at height $1$, that is, in the affine hyperplane $\R^{d}\times \{1\}$. Consider the cone over $P \times \{1\}$. 
\[\cone(P)\coloneqq \R_{\geq 0} (P \times \{1\}).\] 
We obtain a (discrete) semigroup associated to $P$:
\[\cone(P)\cap(\Z^{d}\times \Z).\]
This, in turn, generates a semigroup algebra 
\[\fld^\ast[P]\coloneqq \fld^\ast[\cone(P)\cap(\Z^{d}\times \Z)],\] 
graded by the last coordinate which we will refer to as the \Defn{height}.

We refer to \cite{BG2} for an introduction into the subject. 

Here $\fld$ is any field, though we generally assume the field to be infinite to ensure the existence of an Artinian reduction. In this case $\fld^\ast[P]$ is Cohen-Macaulay by Hochster's theorem \cite{Hochster} with Hilbert series $\mathrm{Ehr}_P(t)$. For a choice of linear system of parameters $\theta_1,\ldots,\theta_{d+1}\in\AR^1(P)$, the Artinian reduction
\[\AR^\ast(P)\coloneqq \bigslant{\fld^\ast[P]}{\langle \theta_1,\ldots,\theta_{d+1}\rangle\fld^\ast[P]}\] 
has $\dim\AR^k(P)=h^\ast_k$.  It follows that the $h^{\ast}_k$ are nonnegative. 

It has been a central question in the theory of lattice polytopes to determine additional properties for the coefficients of the $h^\ast$-polynomial. In particular, Ohsugi and Hibi conjectured that under two special conditions, the coefficients form a unimodal sequence \cite{OH2}, see also \cite{Braun, SL}.

The first of these conditions is that $P$ has the \Defn{integer decomposition property}, short \Defn{IDP}: every lattice point of $\cone(P)$ is a nonnegative integral combination of lattice points in $P \times \{1\}$, or equivalently that $\fld^\ast[P]$ is generated in degree one. We shall often identify $P \times \{1\}$ and $P$, and in particular identify the lattice points of $P$ with those of $P \times \{1\}$; in particular, the monomial generators of $\fld^1[P]$ are in correspondence to lattice points of $P$.

The second property is the \Defn{reflexive} property: there is a lattice point $p$ in $\Z^{d}\times\{1\}$ such that
\[\cone^\circ(P)\cap(\mathbb{Z}^d\times \mathbb{Z}) =\ p+\cone(P)\cap(\mathbb{Z}^d\times \mathbb{Z}) ,\]
where $\cone^\circ(P)$ is the interior of $\cone(P)$.

This is equivalent to $\fld^\ast[P]$ being algebraically Gorenstein with socle degree $d$. In other words, the Artinian reduction is a Poincar\'e duality algebra with socle degree $d$. 

Hence, it also implies, and is in fact known to be equivalent, to a palindromic symmetry analogous to the Dehn-Sommerville relations for polytopes:
\[h_k^{\ast}=h_{d-k}^{\ast}\  \text{for all}\ k\leq \nicefrac{d}{2}.\]

Numerically on the level of the $h^\ast$-vector, the restriction to reflexive polytopes rather than all Gorenstein polytopes is without loss of generality. Bruns and Römer \cite{BR07} showed that for every Gorenstein polytope there is a reflexive polytope with the same $h^\ast$-vector.

We resolve the following conjecture of Ohsugi and Hibi.

\begin{conj}[Ohsugi-Hibi \cite{OH2}]\label{conjecture}
	For any IDP reflexive lattice polytope $P\subset\R^d$, the coefficients of the $h^{\ast}$-polynomial are unimodal:
	\[ h_0^{\ast}\leq h_1^{\ast}\leq\ldots\leq h_{\lfloor \nicefrac{d}{2}\rfloor}^{\ast} = h_{\lceil \nicefrac{d}{2}\rceil}^{\ast} \geq \ldots \geq h_d^{\ast}.\]
\end{conj}

This is the updated form of a conjecture of Hibi \cite{Hibi92}, after Musta{\c t}a and Payne gave an example showing the necessity of the IDP assumption \cite{MP}. These conjectures in turn go back to a more general one of Stanley \cite{Stanley89}, who proposed that the unimodality may hold for a general Gorenstein standard graded integral domain, no doubt motivated by the $g$-conjecture. There are too many partial results in this direction, resolving the Ohsugi-Hibi conjecture for special polytopes, see for instance \cite{Athanasiadis, BR07, BDS, MP, OH, OH2}.

In fact, we shall prove statements that are more powerful than this, and apply to more general cases. For instance, in the case of polytopes that have only the integer decomposition property, we still obtain monotone decreasing coefficients in the second half, i.e.
\[ h^{\ast}_{\lfloor \nicefrac{d+1}{2}\rfloor}\geq \ldots \geq h_d^{\ast}. \]
More importantly, we actually prove algebraic theorems on $\fld^\ast[P]$ that illuminate why and how the integer decomposition property enters. And some that, in particular, and to some degree, apply to all lattice polytopes.

The first (to be stated, yet not to be proven) algebraic result is this, and it immediately implies Conjecture~\ref{conjecture}:

\begin{thm}\label{thm:lef}
	If $P$ is an IDP reflexive polytope, and the characteristic of $\fld$ is $2$ or $0$, then there is an Artinian reduction $\AR^\ast(P)$ of  ${\fld}^\ast[P]$ that has the Lefschetz property, i.e., there is a linear element $\ell\in\AR^1(P)$ such that for any $k\leq \nicefrac{d}{2}$, the map
	\[ \AR^k(P) \xrightarrow{\ \cdot \ell^{d-2k}\ }\ \AR^{d-k}(P)\]
	is an isomorphism.
\end{thm}

The integer decomposition property enters subtly, and we will make clear where and when it happens.

Back to the topic at hand: \emph{generic} shall mean that the Artinian reduction is taken by linear forms 
\[\theta_1,\ldots,\theta_{d+1},\quad \text{where} \quad \theta_i = \sum_{p\in P\cap\mathbb{Z}^d} \theta_{i,p} \x_p,\] with algebraically independent coefficients $\theta_{i,p}$, which necessitates passing to a transcendental field extension $\widetilde{\fld}=\fld(\theta_{i,p})$ of $\fld$. In particular, when we speak of the \emph{generic} Artinian reduction of $\fld^\ast[P]$, then we first pass to the larger field $\widetilde{\fld}$ and associated semigroup algebra $\widetilde{\fld}^\ast[P]$, and perform the Artinian reduction of $\widetilde{\fld}^\ast[P]$ with respect to the linear system $(\theta_i)$.

\begin{con}
Since we consider polytopes of different size, $\widetilde{\fld}$ is only defined within the context, and shall simply denote a sufficiently large transcendental field extension of $\fld$ whenever we need it for the purposes of genericity. For most results, where no specific linear system of parameters is needed, and in particular for those purely concerning classical commutative algebra, an Artinian reduction is just that (without a priori restrictions to the l.s.o.p.), and takes place over the field $\fld$. In contrast, anything \emph{generic} will be highly specific, and understood to consider statements over transcendental extension $\widetilde{\fld}$ of $\fld$ of appropriate size. For instance, when we speak of the generic element $\ell$ in $\fld^1[P]$, then it is not actually an element in $\fld^\ast[P]$ but instead lives in $\widetilde{\fld}^\ast[P]$ where new algebraically independent variables parametrize each coefficient of $\ell$, and is therefore uniquely defined in that field extension.
\end{con}
 
This is a rather extreme choice of linear system of parameters (l.s.o.p.), necessitated by the proof via anisotropy. We want to emphasize the importance of the choice of l.s.o.p., as it makes a crucial difference for the ring we end up working with. 

For a specific, and somewhat canonical, choice of l.s.o.p., $\AR^\ast(P)$ is isomorphic to the orbifold Chow ring of the associated toric Deligne-Mumford stack \cite{BCS}. In contrast to this choice, we make use of the generic Artinian reduction here. The special choice of linear system, as well as consequences particular to that specific linear system, will be discussed in \cite{APP4}.

As graded vector spaces, the results are isomorphic and the inequalities on the dimensions of the graded pieces remain unaffected. 

Yet, as observed in \cite{BD16}, whether an Artinian reduction admits a Lefschetz element depends on the Artinian reduction, and not only on $\fld^\ast[P]$. Braun and Davis gave an example of an IDP reflexive simplex and an Artinian reduction of the associated semigroup algebra which does not even admit a weak Lefschetz element. 

The linear system of parameters chosen there is however not the canonical system for the orbifold Chow ring, nor is it generic.

As the composition of individual multiplications with $\ell$, the Lefschetz isomorphism gives us an injection in the first half and a surjection in the second half, and thus the desired inequalities on the dimensions of the graded pieces:

\begin{thm}[announced in \cite{APPS}]
	The $h^\ast$-polynomial of a reflexive IDP lattice polytope of dimension $d$ has a unimodal sequence of coefficients. Moreover, we have that  \[(h_{i}^{\ast}-h_{i-1}^{\ast})_{1\leq  i\leq \nicefrac{d}{2}}\]  is an $M$-vector in the sense of Macaulay \cite{Macaulay}: It is the Hilbert polynomial of a commutative graded algebra generated in degree one.
\end{thm}	

Here the $M$-vector property too immediately follows from the Lefschetz property in the usual fashion \cite{Stanley87}:
The vector of differences $\max\{h^\ast_i- h^\ast_{i-1},0\}$ is the Hilbert vector of a standard graded algebra, namely 
\[\bigslant{\mathcal{A}^\ast(P)}{\ell \mathcal{A}^\ast(P)}.\]

\subsection*{Idea and setup} The overall idea is based on the recent works of Adiprasito \cite{AHL}, Papadakis and Petrotou \cite{APP, PP}, in that we reduce the Lefschetz property to a property of pairings, introduced as biased pairings in \cite{AHL} and anisotropy in \cite{PP}. 

However, our work requires a critical new ingredient: The aforementioned works are much simplified because we are gifted detailed knowledge of the rings involved, including their fundamental class (also called the \Defn{volume map} in this setting) based on a wealth of previous works that describe the Chow rings of toric varieties, first from the perspective of algebraic geometry, then using combinatorics \cite{Brion, KFT}. 

In the case of lattice polytopes, we have little to work with.  The case of simplicial spheres (and then (pseudo)manifolds and cycles) made use of explicitly combinatorial techniques to reach the desired goal, and the algebra we investigate here is not immediately as governed by a combinatorial structure as the previous one, being cut out by binomials rather than monomials.

And so while Hochster \cite{Hochster} studied the canonical module for the semigroup algebras associated to lattice polytopes, he did not provide enough for us to proceed. In particular, we need a sufficiently explicit description of the fundamental class, that is, the isomorphism between the top nontrivial cohomology and the ground field. We know that such a description exists in principle, of course, see \cite{Bruns-Herzog} and \cite[Theorem 13.4.7]{CoxLittleSchenck}. That said, there is no direct description of this polynomial. Worse, it is not canonically defined; two definitions may differ up to a unit of the base field. We build on Hochster's work and provide a sufficiently explicit description in two ways.

First, we recall the Kustin-Miller normalization of \cite{APP2}. We proposed there a way to define fundamental class uniquely (and not just up to a scalar) at least within a certain field extension that parametrizes Artinian reductions. 

Second, we give two descriptions, which we prove to be equivalent: The first is an identity of Parseval-Rayleigh type. The second is a system of differential equations.

The paper is organized as follows. We start in Section~\ref{sec:lattice} by generalizing the setup in order to state our main theorem and deduce the individual numerical corollaries within Ehrhart theory. For the sake of completeness, we then recall in Section~\ref{sec:ani} the necessary parts of the machinery of \cite{AHL,APP, PP} in order to prove the Lefschetz statements by way of anisotropy. 

Following this setup, we give our new contributions to the theory. Section~\ref{sec:deg} contains the normalization of the fundamental class and an auxiliary identity, while Section~\ref{sec:PC} contains the key identity of Parseval-Rayleigh type, and a differential equation for the volume map we prove to be equivalent. We finish with a discussion of open questions in Section~\ref{sec:discussion}. The appendix contains additional material, such as alternative proofs and illuminating facts that are not necessary on the way to the proof of the main results, but can be helpful independently. In particular we reveal that while semigroup algebras of different polytopes are not immediately related by algebraic maps, they do satisfy some interesting maps that almost behave like pullbacks.

\tableofcontents

\section{The algebraic results, reviewed}\label{sec:lattice}

For $P$ reflexive, the semigroup algebra $\fld^\ast[P]$ is Gorenstein of Krull dimension equal to the dimension of the polytope plus one \cite{Bruns-Herzog}: After an Artinian reduction using a linear system of parameters of length equal to the Krull dimension, we arrive at a Poincar\'e duality algebra of socle degree $d$, that is, a graded ring whose top nontrivial degree is $d$, is one-dimensional as a vector space and so that the ring becomes a Poincar\'e duality algebra with respect to this copy of $\fld$. 

For general IDP polytopes, the situation is a little more delicate. One can force Poincar\'e duality however, using the usual trick: we allow for relative objects. 
The $\fld^\ast[P]$-module defined by
\[ \fld^\ast[P,\partial P]\coloneqq\fld^\ast[\cone^\circ(P)\cap (\mathbb{Z}^d\times \mathbb{Z})],\ \quad \cone^\circ(P)\ =\ \cone(P)\setminus\partial \cone(P),\]
is the canonical module of the Cohen-Macaulay ring $\fld^\ast[P]$, see \cite{Hochster}, or more explicitly \cite[Chapter 6]{Bruns-Herzog}.

After Artinian reduction, we are left with a perfect bilinear pairing
\[\AR^k(P)\ \times\ \AR^{d+1-k}(P,\partial P)\ \longrightarrow\ \AR^{d+1}(P,\partial P).\]
Let us state the first key result:

\begin{thm}\label{thm:rellef}
	If $P$ is an IDP polytope of dimension $d$, and the characteristic of $\fld$ is $2$ or $0$, then some Artinian reduction $\AR^\ast(P)$ of  ${\fld}^\ast[P]$ has the relative Lefschetz property, i.e., there exists a linear element $\ell\in\AR^1(P)$ such that for all $k\leq \nicefrac{d+1}{2}$,
	\[\AR^k(P,\partial P)\ \xrightarrow{\ \cdot \ell^{d+1-2k}\ }\ \AR^{d+1-k}(P)\]
	is an isomorphism. 
\end{thm}

Theorem~\ref{thm:rellef} includes surjections 
\[\AR^{j}(P)\xrightarrow{\cdot\ell}\AR^{j+1}(P)\quad \text{and}\quad \AR^k(P)\xrightarrow{\cdot\ell^{d+1-2k}}\AR^{d+1-k}(P)\] for $j\geq \nicefrac{d+1}{2}$ and $k\leq \nicefrac{d+1}{2}$ and thus we obtain the following result.

\begin{cor}\label{cor:rellly}
	The $h^{\ast}$-polynomial of an IDP lattice polytope $P$ of dimension $d$ has monotone decreasing coefficients in the second half, i.e.,
	\[ h^{\ast}_{\lfloor \nicefrac{d+1}{2}\rfloor}\geq \ldots \geq h_d^{\ast}  \geq h_{d+1}^{\ast}=0. \]
	Moreover, for all $k \le \nicefrac{d+1}{2}$, we have
	\[h^{\ast}_{k}\ \ge\ h^{\ast}_{d+1-k}.\]
\end{cor}

A corollary then concerns Stapledon's $a$-polynomial \cite{Stapledon}.

\begin{cor}
	For any IDP lattice polytope, the $a$-polynomial has unimodal coefficients.
\end{cor}

The last part also follows from Theorem~\ref{thm:lefsphere}, by observing that the $a$-polynomial corresponds exactly to the $h^\ast$-polynomial of $\partial P$ as a {lattice complex}. Here, a lattice complex is a polyhedral complex built out of lattice polytopes, such that the lattice structure agrees on intersections of faces \cite{BG}. We will return to this object in a bit.

In addition, we obtain some interesting consequences if we know at what height the canonical module is generated: Suppose all minimal lattice points of $\cone^\circ(P)$, with respect to the order induced by the semigroup $\cone(P) \cap (\mathbb{Z}^d\times \mathbb{Z})$, are of height at most~$j$, or in other words, $\fld^\ast[P,\partial P]$, the semigroup algebra of interior lattice points of $\cone(P)$, is generated in degree $\le j$ as an ideal over $\fld^\ast[P]$.

\begin{thm}\label{thm:level}
	The $h^\ast$-polynomial of an IDP lattice polytope $P$ of dimension $d$ with $\cone^\circ(P) \cap (\mathbb{Z}^d\times \mathbb{Z})$ generated at height $\leq j$ has monotone increasing coefficients in the initial part, i.e.,
	\[h^{\ast}_0\leq\ldots\leq h^{\ast}_{\lceil\frac{d+1-j}{2}\rceil}.\]
	Moreover, for $k\leq \frac{d+1-j}{2}$ we have 
	\[h_k^\ast\leq h^\ast_{d+1-j-k}.\]
\end{thm}

This is not an immediate consequence of Theorem~\ref{thm:rellef}, but instead follows from an analogous Lefschetz Theorem~\ref{thm:levellef}. 
This in particular includes the result on Gorenstein IDP lattice polytopes Theorem~\ref{thm:lef}, which is the case when  $\AR^\ast(P,\partial P)$ is generated by a single element.

Combining Theorem~\ref{thm:level} and Corollary~\ref{cor:rellly}, we have:

\begin{cor}
	The $h^\ast$-polynomial of an IDP lattice polytope $P$ with $\cone^\circ(P) \cap (\mathbb{Z}^d\times \mathbb{Z})$ generated at height $\leq 3$ has a unimodal sequence of coefficients.
\end{cor}

These results extend to sheafified versions, and in fact, these form an important step: 

Consider an abstract polytopal complex $X$ (that is, a strongly regular CW complex \cite{AB}) whose cells (also called \Defn{faces}) are lattice polytopes, with the property that the lattices agree in common intersections. This is what we call a \Defn{lattice complex}. A \Defn{subcomplex} is a down-closed subset of elements of the lattice complex, that is, $Y$ is a subcomplex of $X$ if whenever $A\in Y \subset X$, and $B$ is a face of $A$, then $B\in Y$. The pair $(X,Y)$ is also called a relative lattice complex, and its elements are the elements of $X$ not in $Y$. Unless we specify further, a statement made about lattice complexes applies to both absolute and relative lattice complex. We obtain naturally also a cone (or fan) over lattice complexes, which is the cone over its elements with the natural attaching map.

We obtain an analogous ring $\fld^\ast[X]$, that naturally generalizes the face ring (or Stanley-Reisner ring):

It is defined quite simply as the direct sum of the individual semigroup rings of the cells, identified at the common intersections of faces, so that two monomials $\x_a$ and $\x_b$, for $a$ and $b$ that do not lie in a common face of $X$, multiply to $0$ (see also Section~\ref{subs!dfnofvariousideals}). We will still use the notation $\x_a\x_b=\x_{a+b}$ in this case, with the convention that $\x_{a+b}=0$ if $a$ and $b$ do not lie in a common face of $X$.

This is a convenient notion: For instance, we have the short exact sequence
\[0\ \longrightarrow\ \fld^\ast[P,\partial P]\ \longrightarrow\ \fld^\ast[P]\ \longrightarrow\ \fld^\ast[\partial P]\ \longrightarrow\ 0.\]

We will simply call this the \Defn{semigroup algebra of $X$}, though it is also known as the \Defn{toric face ring} \cite{IR, Stanley87}. 
We summarize essential results of \cite{BBR, IR}: 
\begin{thm}\label{thm:BBR}
	This semigroup algebra over a lattice complex $X$ is Cohen-Macaulay if $X$ is topologically Cohen-Macaulay. 
	
	In particular, it is Cohen-Macaulay if $X$ is a sphere $\Sigma$ (in which case the ring is also algebraically Gorenstein) or a ball $\Delta$ or a relative ball $(\Delta, \partial \Delta)$, in which case the Poincar\'e pairing applies to the pair of spaces $\AR^\ast(\Delta, \partial \Delta)$ and $\AR^\ast(\Delta)$.
\end{thm}

Motivated by this, we put special focus on lattice balls and spheres, that is, lattice complexes homeomorphic, or $\mathbb{Z}$ homology-equivalent, to a ball resp.\ sphere\footnote{In fact, for all results over characteristic $2$ that concern lattice balls or spheres, one only needs $\mathbb{Z}/2\mathbb{Z}$ homology-equivalence to balls/spheres, and these results extend to characteristic $0$ in absence of torsion, but to simplify statements, we do not concern ourselves too much with this minor generalization.}. We call such complexes IDP if every single one of its cells is IDP. 

In analogy with the proofs of the $g$-theorem of \cite{AHK,APP, PP}, we have:
\begin{thm}\label{thm:lefsphere}
	If $X$ is an IDP lattice sphere or ball of dimension $d$, and the characteristic of $\fld$ is $2$ or $0$, then some Artinian reduction $\AR^\ast(X)$ of  ${\fld}^\ast[X]$ has the Lefschetz property, that is, we have an isomorphism
	\[\AR^k(X,\partial X)\ \xrightarrow{\ \cdot \ell^{d+1-2k}\ }\ \AR^{d+1-k}(X)\]
	with respect to some $\ell\in \AR^{1}(X)$ and all $k\leq \nicefrac{d+1}{2}$.
\end{thm}

Similar results hold for manifolds and cycles (again in direct analogy to \cite{APP}) but seem less immediately relevant here. They can be worked out easily using the methods we provide here, however.

Before we finish this introductory overview, let us briefly discuss the beautiful main insight. And perhaps surprisingly, this holds for lattice polytopes whether they are IDP or not:

\begin{thm}[The Parseval-Rayleigh identity]\label{thm:parsevalcomplex}
	For a lattice $d$-ball or sphere $\Delta$, and $\alpha$ a lattice point in $\cone^\circ(\Delta)\cap (\mathbb{Z}^d\times \{d+1\})$, we have in $\AR^\ast(\Delta,\partial \Delta)$ over characteristic $2$ that
	\begin{equation}\label{eq:id}
		\vol (\x_\alpha) \   =\    \sum_{ \beta \in (\Delta\cap \mathbb{Z}^d)^{d+1}}   \vol (\x_{\frac{\alpha+\beta}{2}})^2   \theta^\beta.
	\end{equation}
\end{thm}

Here, $\vol$ is a certain canonical choice of isomorphism between $\AR^{d+1}(\Delta,\partial \Delta)$ and the ground field we shall specify later, see Section~\ref{sec:deg}, and quite naturally, $\cone^\circ(\Delta)$ is the cone over the interior of $\Delta$. We shall specify the other specifics
of this result in Section~\ref{sec:PC}, but for now simply appreciate how it relates different elements in the semigroup algebra in a surprisingly non-homogeneous way.

Of course, there is a mystery: why then do we need the integer decomposition property? We will clarify this, and also discuss what can be done without it in Section~\ref{sec:beyond}. 

\section{From anisotropy to the Lefschetz property}\label{sec:ani}

To deduce the Lefschetz property from anisotropy, we employ a reduction found originating in \cite{AHL}: It is enough to demonstrate a nondegeneracy property of the Poincar\'e pairing at certain ideals. It is useful to introduce an intermediate property, which we call the \Defn{Hall-Laman relations}, which describe anisotropy of the Hodge-Riemann bilinear form.

The overall strategy is to show that (suitable) anisotropy implies the (suitable) Hall-Laman relations which imply the Lefschetz property for suitable/generic linear elements. The final implication is rather easy, but the first takes a thought: we use the lifting trick of \cite{AHL} to describe Hodge-Riemann pairings of arbitrary semigroup algebras of lattice balls in terms of the middle Poincar\'e pairing in the semigroup algebra of a higher dimensional lattice ball, see Lemma~\ref{lem:midred}.

We note that we give here the derivation for Theorem~\ref{thm:rellef}. For Theorem~\ref{thm:level} we need a different form of the Lefschetz property, and consequently a different form of anisotropy. Nevertheless, since the basic building blocks are the same, we focus on the derivation of the former here.

\subsection{Everything in its right place: a reminder}    \label{subs!dfnofvariousideals}
\phantom {==}
\vspace{5pt}

Let us just remind ourselves, in situ, of the notions going forward.
\subsubsection*{The rings}  Assume $P \subset \mathbb{R}^d$ is a lattice polytope. We do not assume that 
$\dim P = d$.  We denote the lattice point set of $P$ by $V(P)$.  
We consider the convex  cone $\cone(P)$ over $P$ defined by 
\[
\cone(P)\ =\ \{  \sum_{1 \leq t \leq r } \lambda_t  (u_t,1) :   r \geq 1, 
u_t \in P,  \lambda_t  \in     \mathbb{R}_{\geq 0}   \}   \subset \mathbb{R}^{d+1}.
\] 
Then, under addition $ \cone(P) \cap \mathbb{Z}^{d+1}$ is a submonoid of $\mathbb{Z}^{d+1}$ and we set
\[
\fld^\ast[P]\ =\ \fld^\ast[ \cone(P) \cap \mathbb{Z}^{d+1}]. 
\]
In other words,  $\fld^\ast[P ]$ is the graded $\fld$-algebra associated to the monoid  $\cone(P) \cap \mathbb{Z}^{d+1}$.

Consider now a lattice polytopal complex $X$, that is, a collection of lattice polytopes attached along attaching relations $\sim$ and such that the lattice structures coincide.

We denote by $\fld^\ast[X]$ the ring which is obtained as
\[
\fld^\ast[X]\ =\ \bigslant{\bigoplus_{P\ \in\ X} \fld^\ast[P] } {\langle x_\alpha-x_\beta: \alpha \sim \beta\rangle}.
\]
The multiplication, which it suffices to examine on the level of monomials, is naturally described by
$\x_\alpha\cdot \x_\beta = \x_{\alpha+\beta}\in \fld^\ast[P]\subset \fld^\ast[X]$ if $\alpha$ and $\beta$ lie in some $\cone(P)$ for $P\in X$, and $\x_\alpha\cdot \x_\beta=\x_{\alpha+\beta}=0$ otherwise. 

We call this the \Defn{semigroup algebra of $X$}, and it coincides with the \Defn{toric face ring} of Ichim and Römer \cite{IR}. It is naturally associated to $\cone(X)$, which is defined as the cone over the elements of $X$ with the natural attaching maps (and similar for relative complexes). Such collections of cones are also called (abstract) fans. 

For a pair of lattice complexes $\Psi:=(X,Y),\ Y\subset X$, we define naturally the module 
$\fld^\ast[\Psi]$ as the kernel of the natural surjection
$\fld^\ast[X]\ \twoheadrightarrow\ \fld^\ast[Y]$. 

We denote by $\AR^\ast(P), \AR^\ast(X), \AR^\ast(\Psi)$ etc.\ the Artinian reductions of the respective rings, and unless further specified, we have a linear system of parameters indexed by coefficients $\theta_{i,j}$.

\subsubsection*{Some natural maps}  Assume $\Delta$ is a lattice polytopal ball or sphere of dimension $d$.   
The inclusion  map  
\begin{equation}   \label{eqn!map000}
	{\fld}^\ast [\Delta,  \partial \Delta]  \ \hookrightarrow\  {\fld}^\ast [\Delta]  
\end{equation}            
induces, for all $k \geq 0$, natural maps
\begin{equation}   \label{eqn!map001}
	\AR^{k}(\Delta ,\partial \Delta) \ \longrightarrow \AR^{k}(\Delta)
\end{equation} 
which, in general, are not injective. For example, for $k=d+1$, we have $ \AR^{k}(\Delta)=0$ while
$\AR^{k}(\Delta ,\partial \Delta)$ is isomorphic to the ground field (that is, $\fld$ or $\widetilde{\fld}$ depending on whether we are just considering an arbitrary Artinian reduction over the original field, or the generic Artinian reduction over the field extension).

The multiplication map 
\begin{equation}   \label{eqn!map002}
	{\fld}^\ast [\Delta]  \times    {\fld}^\ast [\Delta,  \partial \Delta]   \longrightarrow {\fld}^\ast [\Delta,  \partial \Delta]   
\end{equation} 
induces, for all $k,l$, multiplication maps
\begin{equation}   \label{eqn!map003}
{\fld}^k  [\Delta] \times    {\fld}^l  [\Delta,  \partial \Delta]  \longrightarrow {\fld}^{k+l}[\Delta,  \partial \Delta].
\end{equation} 
There is an induced multiplication map
\begin{equation}   \label{eqn!map004}
	\AR^\ast  (\Delta) \times    \AR^\ast  (\Delta,  \partial \Delta)  \longrightarrow  \AR^\ast(\Delta,  \partial \Delta)
\end{equation} 
giving   $\AR^\ast  (\Delta,  \partial \Delta) $ the structure of $\AR^\ast  (\Delta)$-module. This map
induces, for all $k,l$,  multiplication  maps 
\begin{equation}   \label{eqn!map005}
	\AR^k  (\Delta) \times    \AR^l  (\Delta,  \partial \Delta)  \longrightarrow  \AR^{k+l}(\Delta,  \partial \Delta).
\end{equation} 

We remark that if  $k+l = d+1$, the pairing 
\begin{equation}   \label{eqn!map0033}
	\AR^k  (\Delta) \times    \AR^{d+1-k} (\Delta,  \partial \Delta)  \longrightarrow  \AR^{d+1}(\Delta,  \partial \Delta) 
\end{equation}   
is perfect, see \cite{Hochster} and \cite[Chapter 6]{Bruns-Herzog}.

Assume $p$ is a homogeneous polynomial of degree $s$ with $0 \leq s \leq d+1$.
For $k \geq 0$ we set
\begin{equation}   \label{dfn!lsoaldaga}
	p \AR^k  (\Delta, \partial \Delta )  = \{ p u :  u \in  \AR^k  (\Delta, \partial \Delta ) \} 
	\subset   \AR^{k+s}  (\Delta, \partial \Delta ) 
\end{equation}
and
\begin{equation}   \label{dfn!msikdoas}
	p \AR^k  (\Delta)  = \{ p u :  u \in  \AR^k  (\Delta ) \} 
	\subset   \AR^{k+s}  (\Delta). 
\end{equation}
Note that for all $k$, there is a well-defined perfect pairing 
\begin{equation}   \label{eqn!mapois90a}
	\phi_{p}  :  p \AR^k  (\Delta) \times   p \AR^{d+1-k-s} (\Delta,  \partial \Delta)  \longrightarrow  
	\AR^{d+1}(\Delta,  \partial \Delta)
\end{equation}   
defined by  $ \phi_{p}(pa, pb) =  pab$ for all $a \in \AR^k  (\Delta) $
and $b \in \AR^{d+1-k-s} (\Delta,  \partial \Delta) $. It is clear that this pairing is well-defined: Assume first  that  $c \in  \AR^{k}  (\Delta)$ 
has the property that
$pc =0$ in $\AR^{k+s}  (\Delta)$.  Then for all $b \in \AR^{k} (\Delta,  \partial \Delta)$
we get that $pcb = 0$.   Assume now  that  $c \in  \AR^{d+1-k-s}  (\Delta, \partial \Delta)$ 
has the property that
$pc =0$ in $\AR^{d+1-k}  (\Delta,\partial \Delta)$.  Then for all 
$a \in \AR^{k} (\Delta)$
we get that $pac = apc = 0$.

We also note it is a perfect pairing. Assume first  $a \in  \AR^{k}  (\Delta)$ 
has the property that $pa \not= 0 $ in $\AR^{k+s}  (\Delta)$.  Then, since the pairing
in Equation~(\ref{eqn!map0033})  is perfect, there exists $b \in \AR^{d+1-k-s}  (\Delta, \partial \Delta)$
such that $pab \not= 0$. Hence  $ \phi_{p} (pa, pb) \not= 0$. 
Assume now $b \in \AR^{d+1-k-s}  (\Delta, \partial \Delta) $ 
has the property that $pb \not= 0 $ in $\AR^{d+1-k}  (\Delta, \partial \Delta)$.  
Then, since the pairing
in Equation~(\ref{eqn!map0033})  is perfect, there exists $a \in \AR^{k}  (\Delta)$
such that $apb \not= 0$. Hence  $ \phi_{p} (pa, pb) \not= 0$.

\begin{rem}\label{rem:annGorenstein}
	Another way to think about the pairing of Equation~\ref{eqn!mapois90a} is to consider it as a perfect pairing 
	\[\bigslant{\AR^k  (\Delta)}{\mr{ann}\, p} \times  \bigslant{\AR^{d+1-k-s} (\Delta,  \partial \Delta)}{\mr{ann}\, p} \ \longrightarrow\  
	\bigslant{\AR^{d+1-s} (\Delta,  \partial \Delta)}{\mr{ann}\, p} .\]
	Of course, this is just a shift of degree: 
	${\AR^k  (\Delta)}/{\mr{ann}\, p}\ \cong\ p \AR^k  (\Delta)$ and \[\bigslant{\AR^{d+1-k-s} (\Delta,  \partial \Delta)}{\mr{ann}\, p}\ \cong\ p \AR^{d+1-k-s} (\Delta,  \partial \Delta).\]
\end{rem}

\begin{rem}   Assume $\mathcal{I}^\ast$ is a nonzero graded  $ \AR^\ast  (\Delta)$-submodule of  $ \AR^\ast  (\Delta, \partial \Delta ) $.
	Then, there exists $k \geq 0$ and nonzero  $u \in \mathcal{I}^{k}$.  Using the perfect pairing
	\[
	\AR^k  (\Delta) \times    \AR^{d+1-k} (\Delta,  \partial \Delta)  \longrightarrow  \AR^{d+1}(\Delta,  \partial \Delta)  
	\]
	of Equation~(\ref{eqn!map0033}), 
	there exists $w \in   \AR^{d+1-k} (\Delta)$ such that $wu$ is a nonzero element of $\AR^{d+1}(\Delta,  \partial \Delta)$. 
	Hence,   $wu$ is a nonzero element of $ \mathcal{I}^{d+1}$. Since $\AR^{d+1}(\Delta,  \partial \Delta)$ is a $1$-dimensional
	vector space over $ \fld$, it follows that $ \mathcal{I}^{d+1} = \AR^{d+1}(\Delta,  \partial \Delta). $
\end{rem}  

\subsection{Anisotropy} 
The prototype of anisotropy is the following. 

\begin{thm}\label{thm:ani}
	If $X$ is an IDP lattice ball or sphere of dimension $d$, and the characteristic of $\fld$ is $2$ or $0$, then the 
	generic Artinian reduction $\AR^\ast(X, \partial X)$ of ${\fld}^\ast[X,\partial X]$ has the anisotropy property. 
	This means that  for every nonzero $u\in\AR^{k}(X,\partial X)$ of degree $k\le \nicefrac{(d+1)}{2}$, we 
	have \[u^2\ \neq\ 0,\]
	in other words  $u^2$ is a nonzero element of $\AR^{2k}(X,\partial X)$.
	Moreover, if $m$ is a monomial of degree $\le d+1-2k$ such that $m u$ is nonzero in
	$\AR^{k+ \deg(m)}(X,\partial X)$, then 
	\[m u^2\ \neq\ 0, \]
	in other words   $m u^2$ is a nonzero element of $\AR^{2k+ \deg(m)}(X,\partial X)$.
\end{thm}

Most of the remainder of this paper is devoted to proving this result 
(see Section~\ref{sec:proofani}), and related results, see Section \ref{sec:anilevel}. Let us note that it is enough
to prove this theorem in characteristic $2$, see \cite[Section~3.2]{KLS}. 
Before we do that, however, we follow the derivation of the Lefschetz property from it, based on the ideas of \cite{AHL, APP, PP}. 

\subsection{The Hall-Laman relations}
Consider a lattice ball $\Delta$ of dimension $d$.
Assume  that $k\le \frac{d+1}{2}$,  $\ell \in \AR^1(\Delta)$ and $\mathcal{I}^\ast\subset \AR^\ast(\Delta,\partial \Delta)$
is a nonzero graded submodule. 
We say that $\AR^\ast(\Delta,\partial \Delta)$ satisfies the \Defn{Hall-Laman relations} 
for the triple $(k, \ell,   \mathcal{I}^\ast)$ if the pairing
\begin{equation}\label{eq:sl}
	\begin{array}{ccccc}
		\mathcal{I}^k& \times &\mathcal{I}^k & \longrightarrow &\ \AR^{d+1}(\Delta,\partial \Delta)\\
		a	&		& b& {{\xmapsto{\ \ \ \ }}} &\ ab\ell^{d+1-2k}
	\end{array}
\end{equation}
is nondegenerate. We say that $\AR^\ast(\Delta,\partial \Delta)$ satisfies the \Defn{absolute Hall-Laman relations}
with respect to the pair    $(k, \ell)$
if for {\it all} nonzero graded submodules $\mathcal{I}^\ast\subset \AR^\ast(\Delta,\partial \Delta)$
the Hall-Laman relations 
are true for the triple $(k, \ell,   \mathcal{I}^\ast)$.  

\begin{prp}  \label{prop!wkidlsfg}
	The absolute Hall-Laman relations are true for the pair $(k, \ell)$
	if and only if the Hodge-Riemann bilinear form 
	\begin{equation*}\label{eq:sl}
		\begin{array}{rccccc}
			Q_{\ell,k}:& \AR^k(\Delta,\partial \Delta)& \times &\AR^k(\Delta,\partial \Delta)
			& \longrightarrow &\ \AR^{d+1}(\Delta,\partial \Delta) \\
			&	a	&		& b& {{\xmapsto{\ \ \ \ }}} &\ ab\ell^{d+1-2k}
		\end{array}
	\end{equation*}
	is anisotropic in the following sense: if 
	$u\in \AR^k(\Delta,\partial \Delta)$ is not zero, then $Q_{\ell,k}(u,u)\neq 0$.
\end{prp} 

\begin{proof}
	The absolute Hall-Laman relations state that the Hodge-Riemann bilinear form 
	does not degenerate at any  nonzero graded submodule  $\mathcal{I}^\ast$ of
	$\AR^\ast(\Delta,\partial \Delta)$. 
	But it clearly suffices to verify this fact at principal ideals generated by 
	single elements $u$ in $\AR^k(\Delta,\partial \Delta)$, which in turn is the property 
	of anisotropy. The other direction is clear.
\end{proof}

Hence, since every anisotropic symmetric bilinear form is nondegenerate, 
we get from  Proposition~\ref{prop!wkidlsfg} that the absolute Hall-Laman 
relations for the pair $(k,\ell)$  imply a Lefschetz type property
for $\AR^\ast(\Delta,\partial \Delta)$    at degree $k$.   

We introduce refinements of  these properties:
Assume $p$ is a homogeneous polynomial of degree $s$ in $\fld^\ast[\Delta]$,   $k$ is an integer with $k\le \frac{d+1-s}{2}$, 
$\ell \in  \AR^1(\Delta)$ and $\mathcal{I}^\ast \subset \AR^\ast(\Delta,\partial \Delta)$ is a nonzero graded submodule.
We say that $\AR^\ast(\Delta,\partial \Delta)$ satisfies the 
\Defn{Hall-Laman relations with respect to the quadruple  $(k,p,\ell,  \mathcal{I}^\ast)$} 
if the pairing
\begin{equation*}\label{eq:dsdsl}
	\begin{array}{rccccc}
		& p \AR^k(\Delta,\partial \Delta)& \times & p \AR^k(\Delta,\partial \Delta)
		& \longrightarrow &\ \AR^{d+1}(\Delta,\partial \Delta) \\
		& p	a	&		& p b& {{\xmapsto{\ \ \ \ }}} &\ pab\ell^{d+1-s-2k}
	\end{array}
\end{equation*}
(which is well-defined by the obvious argument)  is perfect when 
restricted to  $p\mathcal{I}^k  \times p\mathcal{I}^k$.   If this is true for every nonzero graded
submodule   
$\mathcal{I}^\ast \subset \AR^\ast(\Delta,\partial \Delta)$   we say $\AR^\ast(\Delta,\partial \Delta)$ satisfies the 
\Defn{absolute Hall-Laman relations with respect to the triple  $(k,p,\ell)$}.  
Arguing similarly as in the proof of Proposition~\ref{prop!wkidlsfg} 
this is equivalent
to saying that for all  $u \in \AR^{k}(\Delta,\partial \Delta)$ such that  $pu\in \AR^{k+s}(\Delta,\partial \Delta)$ is nonzero
it holds that  $pu^2\ell^{d+1-s-2k}$ is a nonzero element of $\AR^{d+1}(\Delta,\partial \Delta)$.

Recall 
the definition of $p \AR^k  (\Delta, \partial \Delta ) $ in 
Equation~(\ref{dfn!lsoaldaga}),  the  definition of $p \AR^k  (\Delta) $ in 
Equation~(\ref{dfn!msikdoas}) and the perfect pairing $\phi_p$ 
in  Equation~(\ref{eqn!mapois90a}).
For a vector subspace $U$ of $\AR^{k}  (\Delta, \partial \Delta )$ we set
\begin{equation}   \label{dfn!toslaosfsa}
	p U   = \{ p u :  u \in U \} 
	\subset   p \AR^{k}  (\Delta, \partial \Delta )
\end{equation}
and
\begin{equation}   \label{dfn!paoeosaf}
	{\mr{ann}(p U)} = \{ p b  :   b \in  \AR^{d+1-k-s}(\Delta)    \;  
	\; \text{ and }  \; pab = 0      \; \text{ for all }  \;a \in U \}.   
\end{equation}
We have that $pU$ is a vector subspace of $p \AR^{k} (\Delta, \partial \Delta )$  and 
${\mr{ann}(p U)}$ is a vector subspace of  $p\AR^{d+1-k-s} (\Delta)$. Moreover,
\begin{equation}   \label{dfn!aitidsd}
	{\mr{ann}(p U)} = \{ x \in p \AR^{d+1-k-s}(\Delta)\;  : \;  \phi_p ( x, z) = 0 
	\; \text{ for all }  \;z \in pU \}. 
\end{equation}

\begin{lem}\label{lem:triv}
	The following are equivalent:
	\begin{compactenum}[(1)]
		\item $\AR^\ast(\Delta,\partial \Delta)$ satisfies the Hall-Laman relations with respect to the quadruple   
		$(k,p,\ell, \mathcal{I}^\ast )$.
		\item The map 
		\begin{equation} \label{eqn!frerwedsasf}
			p\mathcal{I}^k\ \xrightarrow{\ \cdot \ell^{d+1-s-2k}\ } 
			\bigslant{p \AR^{d+1-k-s}(\Delta)}{\mr{ann} (p\mathcal{I}^{k})}
		\end{equation}
		is an isomorphism. 
	\end{compactenum}
\end{lem}

Here, we set  $t =   d+1-k-s$. 
The map in Equation~(\ref{eqn!frerwedsasf}) sends $pa$ to the class of $pa \ell^{t-k}$ in
${p \AR^{t}(\Delta)}/{\mr{ann} (p\mathcal{I}^{k})}$ for all
$a \in  \mathcal{I}^k$.  Moreover,  it  is the composition of the 
inclusion  map $p\mathcal{I}^k \rightarrow p \AR^{k}(\Delta,\partial \Delta)$ map
with the map induced by multiplication with $ \ell^{t-k}$
\[p \AR^{k}(\Delta,\partial \Delta) \xrightarrow{\ \cdot \ell^{t-k}\ } p \AR^{t}(\Delta,\partial \Delta)\]
(compare Equation~(\ref{eqn!map004})), with the natural map 
$p \AR^{t}(\Delta,\partial \Delta) \rightarrow p \AR^{t}(\Delta)$
(compare Equation~(\ref{eqn!map001})),  with the natural quotient map 
\[ p\AR^{t}(\Delta)\ \longrightarrow \ \bigslant{p\AR^{t}(\Delta)}{\mr{ann}(p\mathcal{I}^{k})}.\]

\begin{proof}   We define
	\[
	\rho_1 :  p \mathcal{I}^{k}  \times   p \mathcal{I}^{k}
	\longrightarrow \ \AR^{d+1}(\Delta,\partial \Delta)\\
	\]
	with $\; \rho_1  ( pa, pb) = pab\ell^{d+1-s-2k} \; $ for all  $a,b \in \mathcal{I}^{k}$.
	As observed before, $\rho_1$ is well-defined.

	Recall that if  $\rho : V \times W \rightarrow \fld$ is a perfect pairing of
	finite dimensional $\fld$-vector spaces, $X$ is a vector subspace of $V$,
	and we set  
	\[
	\mr{ann}(X) = \{ w \in W:  \rho (x,w) = 0 \;  \text{ for all } \;  x \in X  \},
	\]
	then it follows that  the pairing $\rho$ induces a perfect pairing $ \rho' : X \times (W/\mr{ann}(X))  \rightarrow \fld$
	such that  $\rho' (x, [w]) = \rho (x,w)$ for all $x \in X, w \in W$. Applying that to the perfect pairing  $ \phi_{p}$  defined in 
	Equation~(\ref{eqn!mapois90a}) we get an induced perfect pairing
	\[
	\rho_2 :  p \mathcal{I}^{k}  \times    \bigslant{p \AR^{d+1-k-s}(\Delta)}{\mr{ann} (p\mathcal{I}^{k})}
	\longrightarrow \ \AR^{d+1}(\Delta,\partial \Delta).
	\]
	
	We define 
	\[
	\upsilon :  p \mathcal{I}^{k}  \times   p \mathcal{I}^{k} \longrightarrow 
	p \mathcal{I}^{k}  \times   \bigslant{p \AR^{d+1-k-s}(\Delta)}{\mr{ann} (p\mathcal{I}^{k})}
	\] 
	by  $\upsilon  (p a, p b)  $ =  $(pa, pb  \ell^{d+1-s-2k})$ for all $a,b \in \mathcal{I}^{k}$. 
	It is easy to see that $\upsilon$ is well-defined and that
	\[
	\rho_1 = \rho_2 \circ \upsilon.
	\]
	Since $\rho_2$ is a perfect pairing, we get that $\rho_1$ is a perfect pairing if and only if
	$\upsilon$ is  an isomorphism.  The result follows.
\end{proof}

\subsection{Pyramids}  

To prove the Hall-Laman relations, let us introduce an auxiliary construction:

\begin{dfn}[Pyramids]  \label{}
	Given a lattice polytope $P$ with ambient lattice $\Z^d$, the \Defn{pyramid} over \Defn{base} $P$ is 
	constructed as the convex hull of $P\times \{0\}$ and $\mbf{a}=(0,\cdots,0,1)$, the \Defn{apex}, 
	in $\Z^d \times \Z$. We also denote this as $\mr{pyr}_{\mbf{a}}P$. The pyramid over a lattice complex is the collection of pyramids over its elements.
	
	If $X$ is a lattice complex, then $\mr{pyr}_{\mbf{a}} X$ is the lattice complex consisting of cones $\mr{pyr}_{\mbf{a}}P,$ $P\in X$ together with $X$ itself. The pyramid over a lattice ball is also a lattice ball. The following
	lemma is immediate from the definitions. Naturally, a relative complex $(X,Y)$ has, associated to it, a pyramid $(\mr{pyr}_{\mbf{a}} X,\mr{pyr}_{\mbf{a}} Y)$.	
\end{dfn}

\begin{lem}
	The pyramid over a lattice polytope is IDP if and only if the base is IDP.  The pyramid over a lattice polytopal complex is IDP
	if and only if the base is IDP.
\end{lem}

This is an immediate consequence of the following observation:

\begin{lem}\label{lem:pyropyro}
Every lattice point of $\cone (\mr{pyr}_{\mbf{a}} P) \cap (\Z^d \times \Z \times \Z)$ is of the form $\lambda \mbf{a} + z$, where $\lambda$ is a nonnegative integer and 
\[z\in \cone (P \times \{0\}) \cap (\Z^d \times \Z \times \Z).\]
\end{lem}

\subsection{The pyramid lemma}
On the level of semigroup algebras, a pyramid corresponds to the introduction of a new indeterminate, 
corresponding to the apex. Consider a lattice complex $X$ of dimension $d$.
If $X$ is a lattice ball  then $\mr{pyr}_{\mbf{a}} X$ is a lattice ball of dimension $d+1$. Consider the semigroup algebra over $X$.

We use the linear system of parameters $\Theta\x=(\theta_{i,j}) \x$ parametrized as usual by a matrix $(\theta_{i,j})$ with $d+2$ rows, where $i$ ranges from $1$ to $d+2$, 
$j$ ranges over the lattice points of $\mr{pyr}_{\mbf{a}} X$ and $\x$ is the vector of indeterminates, that is,  the column matrix with $(1,j)$-entry equal to $\x_j$. 

Without loss of generality we assume that the last column of the matrix $(\theta_{i,j})$ corresponds to the apex ${\mbf{a}}$, and is equal to transpose of $(0,0,\cdots,0,1)$.

We consider the Artinian reduction $\AR^\ast(\mr{pyr}_{\mbf{a}} X)$ of  ${\fld}^\ast[\mr{pyr}_{\mbf{a}} X]$ with respect to this linear system
of parameters, 
and the Artinian reduction $\AR^\ast(X)$ of $\fld^\ast[X]$ with respect to the linear system of parameters 
specified by the first $d+1$  rows of $\Theta\x$. 

Set $\overline{h}=\theta_{d+2}$, the last entry of the linear system of parameters; it plays a special role. With respect to the above linear 
systems of parameters,  the following lemma is a straightforward computation.

\begin{lem}[Pyramid lemma]\label{lem:pyramid}
	Consider $X$ as above, and $Y$ a possibly empty subcomplex of $X$. Consider the pair $\Psi=(X,Y)$.
	\begin{compactenum}[(1)]
		\item The inclusion $\fld^\ast[X]	\hookrightarrow \fld^\ast[\mr{pyr}_{\mbf{a}} X]	$ induces, for all $m$,  an  isomorphism
		\[
		\AR^{m}(\Psi)  \longrightarrow \AR^{m}(\mr{pyr}_{\mbf{a}} \Psi).
		\]

		\item [(2)]  The multiplication 
		\[  \fld^\ast[\mr{pyr}_{\mbf{a}} X] \xrightarrow{\ \cdot  \x_{\mbf{a}}\ }\fld^\ast[\mr{pyr}_{\mbf{a}} X] \] 
		induces, for all $m$,
		an  isomorphism
		\[ \fld^m [\mr{pyr}_{\mbf{a}} \Psi]\ \cong \ \fld^{m+1}[\mr{pyr}_{\mbf{a}} X,X \cup \mr{pyr}_{\mbf{a}} Y] .\] 
		
	\end{compactenum}
	
	In particular, if $\Delta$ is a lattice ball, we have an isomorphism
	\[\AR^{m}(\Delta,\partial \Delta)\  \to\ \AR^{m}(\mr{pyr}_{\mbf{a}}\Delta,\mr{pyr}_{\mbf{a}}\partial \Delta)\  \to\ \AR^{m+1}(\mr{pyr}_{\mbf{a}} \Delta,\partial \mr{pyr}_{\mbf{a}} \Delta). \hspace*{\fill}\]
\end{lem}

\begin{proof}
Part $(1)$ is a consequence of
\[\fld^{\ast}[\Psi] \ \cong\ \bigslant{\fld^{\ast}[\mr{pyr}_{\mbf{a}} \Psi]}{\overline{h}\fld^{\ast}[\mr{pyr}_{\mbf{a}} \Psi]}.\]
Part $(2)$ is trivial since the multiplication with $\x_{\mbf{a}}$ induces a bijection of the monomials generating the modules in questions by Lemma~\ref{lem:pyropyro}:
Observe that 
\[\fld^{m+1}[\mr{pyr}_{\mbf{a}} X,X]\ \cong\  \x_{\mbf{a}}\fld^{m}[X],\]
and analogously for $Y$. The claim then follows from the short exact sequences 
\[0\ \longrightarrow\ \fld^{\ast}[\Psi] \ \longrightarrow\ \fld^{\ast}[X]\ \longrightarrow\ \fld^{\ast}[Y]\ \longrightarrow\ 0\]
and
\[0\ \longrightarrow\ \fld^{\ast}[\mr{pyr}_{\mbf{a}} X,X \cup \mr{pyr}_{\mbf{a}} Y] \ \longrightarrow\ \fld^{\ast}[\mr{pyr}_{\mbf{a}} X,X]\ \longrightarrow\ \fld^{\ast}[\mr{pyr}_{\mbf{a}} Y,Y]\ \longrightarrow\ 0,\]
combined using the snake lemma.
\end{proof}

\subsection{Reduction via pyramidal lifting}

Consider now the case when $X=\Delta$, where $\Delta$ is a lattice ball of dimension $d$. The crucial lemma is the following.

\begin{lem}[compare {\cite[Lemma 7.5]{AHL}}]\label{lem:midred}
Let $k< \frac{d+1}{2}$ and $\IR^\ast$ be a nonzero graded submodule of  $\AR^\ast(\Delta,\partial \Delta)$. 
We consider the induced graded submodule $ \x_{\mbf{a}}\IR^\ast$ of $\AR^{\ast}(\mr{pyr}_{\mbf{a}} \Delta, \partial \mr{pyr}_{\mbf{a}} \Delta)$.
We also  set $h=\x_{\mbf{a}}-\overline{h}$.
Then the following two are equivalent:
\begin{compactenum}
	\item [(1)] The Hall-Laman relations for the triple  $(k+1,  \x_{\mbf{a}},  \x_{\mbf{a}}\IR^\ast)$.
	
	\item  [(2)] The Hall-Laman relations for the triple $(k,  h,   \IR^\ast)$.
\end{compactenum}

This extends naturally to the Hall-Laman relations relative to a homogeneous polynomial $p$ in $\fld^\ast[\Delta]$.
We denote the degree of   $p$   by  $s$,  and assume  $k< \frac{d+1-s}{2}$. Then the following two are equivalent:
\begin{compactenum}
	\item [(3)]  The Hall-Laman relations for the quadruple  $(k+1,p,\x_{\mbf{a}},   \x_{\mbf{a}}\IR^\ast) $.								
	
	\item [(4)]  The Hall-Laman relations for the quadruple  $(k,p, h,   \IR^\ast)$.
\end{compactenum}
\end{lem}

Let us note that we actually do only need the case when $d-s-2k=1$. But the lemma deserves to be stated fully nevertheless.

\begin{proof}
In the present proof all vertical maps are coming from the pyramid Lemma~\ref{lem:pyramid}.	
Moreover, we remark that    $h=\x_{\mbf{a}}$ in  $\AR^{\ast}(\mr{pyr}_{\mbf{a}} \Delta)$. Of course the first part is a special case of the second part (by setting $p=1$), but we felt it to be didactically helpful to nevertheless go over everything.

\vspace{10pt}

CASE 1.  We assume that   $p=1$ and $\IR^\ast = \AR^{*}(\Delta,\partial \Delta )$. 
We consider the  commutative diagram
\[\begin{tikzcd}[column sep=5em]
	\AR^{k}(\Delta,\partial \Delta ) \arrow{r}{\ \ \ \ \cdot h^{d+1-2k}\ \ \ \  } \arrow{d}{\sim } & \AR^{d+1-k}(\Delta) \arrow{d}{\sim } \\
	\AR^{k+1}(\mr{pyr}_{\mbf{a}} \Delta, \partial \mr{pyr}_{\mbf{a}} \Delta)  \arrow{r}{\ \ \ \cdot \x_{\mbf{a}}^{d-2k}\ \ \ } & {\AR^{d+1-k}(\mr{pyr}_{\mbf{a}} \Delta)}
\end{tikzcd}
\]	
The horizontal map on the top being an isomorphism is equivalent to the horizontal map on the bottom being an isomorphism.

\vspace{10pt}
CASE 2.  We assume that   $p=1$ and $\IR^\ast$  is a nonzero graded submodule of   $ \AR^{*}(\Delta,\partial \Delta )$.
For $U$ a vector subspace of $\AR^{k}(\Delta,\partial \Delta )$ we define
\begin{equation}  \label{eqn!ofksmncd}
	\mr{ann}_1 (U) =  \{ x \in \AR^{d+1-k}(\Delta)\;  : \;   xz = 0   \; \text{ for all }  \;z \in U \}
\end{equation}
and for $W$ a vector subspace of $\AR^{k+1}(\mr{pyr}_{\mbf{a}} \Delta, \partial \mr{pyr}_{\mbf{a}} \Delta)$ we define
\begin{equation}  \label{eqn!tidkiaksd}
	\mr{ann}_1 (W) =  \{ x \in \AR^{d+1-k}(\mr{pyr}_{\mbf{a}} \Delta)\;  : \;  xz = 0   \; \text{ for all }  \;z \in W \}.
\end{equation}
We consider the  isomorphism  
\[ 
\tau_1 :  \AR^{d+1-k}(\Delta) \longrightarrow  {\AR^{d+1-k}(\mr{pyr}_{\mbf{a}} \Delta)}
\] 
of Part (1) of Lemma~\ref{lem:pyramid}. It has the property that for any element 
$u \in \fld^{d+1-k}[\Delta]$ it holds $\tau_1 [u] = [u]$.

We claim that  
\begin{equation}   \label{eqn!mikod}
	\tau_1 (\mr{ann}_1 ( \IR^{k})) =  \mr{ann}_1 (\x_{\mbf{a}}  \IR^{k}).
\end{equation}   
The inclusion $\tau_1 (\mr{ann}_1 ( \IR^{k}))  \subset  \mr{ann}_1 (\x_{\mbf{a}}  \IR^{k})$ is obvious.
We assume that $u \in \fld^{d+1-k}[\Delta]$ has the property that  $u$ is not an element 
of $ \mr{ann}_1 ( \IR^{k})$. By the perfect pairing of Equation~(\ref{eqn!map0033}) there exists
$w \in \AR^{k}(\Delta, \partial \Delta)$ such that $uw$ is a nonzero element of $\AR^{d+1}(\Delta, \partial \Delta)$,
Using the isomorphism of Part (3) of Lemma~\ref{lem:pyramid} we get that 
$\x_{\mbf{a}} uw$  is a nonzero element of $\AR^{d+1}(\mr{pyr}_{\mbf{a}} \Delta), \partial \mr{pyr}_{\mbf{a}} \Delta)$.
Hence, $\tau_1 (u)$ is not an element of  $\mr{ann}_1 (\x_{\mbf{a}}  \IR^{k})$.

Using   Equation~(\ref{eqn!mikod}), we have, similarly to 
Case 1,  a  commutative diagram
\[\begin{tikzcd}[column sep=5em]
	\IR^{k}    \arrow{r}{\ \ \ \ \cdot h^{d+1-2k}\ \ \ \  } \arrow{d}{\sim } & \AR^{d+1-k}(\Delta)/{\mr{ann}_1 ( \IR^{k})}      
	\arrow{d}{\sim } \\
	\x_{\mbf{a}}  \IR^{k}   \arrow{r}{\ \ \ \cdot \x_{\mbf{a}}^{d-2k}\ \ \ } & {\AR^{d+1-k}(\mr{pyr}_{\mbf{a}} \Delta)}/ {\mr{ann}_1 (\x_{\mbf{a}}  \IR^{k}) }
\end{tikzcd}
\]	
The horizontal map on the top being an isomorphism is equivalent to the horizontal map on the bottom being an isomorphism.

\vspace{10pt}

CASE 3.   We assume that   $p$ is a homogeneous  polynomial of degree $s$  and $\IR^\ast = \AR^{*}(\Delta,\partial \Delta )$. 
We consider the  commutative diagram
\[\begin{tikzcd}[column sep=5em]
	p \AR^{k}(\Delta,\partial \Delta ) \arrow{r}{\ \ \ \ \cdot h^{d+1-s-2k}\ \ \ \  } \arrow{d}{\sim } & p \AR^{d+1-k-s}(\Delta) \arrow{d}{\sim } \\
	p \AR^{k+1}(\mr{pyr}_{\mbf{a}} \Delta, \partial \mr{pyr}_{\mbf{a}} \Delta)  \arrow{r}{\ \ \ \cdot \x_{\mbf{a}}^{d-s-2k}\ \ \ } & p {\AR^{d+1-k-s}(\mr{pyr}_{\mbf{a}} \Delta)}
\end{tikzcd}
\]	
The horizontal map on the top being an isomorphism is equivalent to the horizontal map on the bottom being an isomorphism.

\vspace{10pt}

CASE 4.  We now just combine the reasonings of Case 2. and Case 3.
We assume that   $p$ is a homogeneous  polynomial of degree $s$ and $\IR^\ast$  is a nonzero 
graded submodule of   $ \AR^{*}(\Delta,\partial \Delta )$.  We set $t = d+1-k-s$.  
Arguing similarly as in Case 2, we have that under the natural map 
\[    
p \AR^{t}(\Delta) \longrightarrow p  {\AR^{t}(\mr{pyr}_{\mbf{a}} \Delta)}
\]  
of Part (1) of the   pyramid Lemma~\ref{lem:pyramid}
the submodule ${\mr{ann} ( p\IR^{k})} $ of $p \AR^{t}(\Delta)$ maps isomorphically onto the submodule  
${\mr{ann} (\x_{\mbf{a}} p \IR^{k}) }$ of $ p  {\AR^{t}(\mr{pyr}_{\mbf{a}} \Delta)}$,
where $\mr{ann}$ was defined in Equation~(\ref{dfn!aitidsd}). Hence, 
similarly to Case 3 we have a  commutative diagram
\[\begin{tikzcd}[column sep=5em]
	p \IR^{k}    \arrow{r}{\ \ \ \ \cdot h^{d+1-s-2k}\ \ \ \  } \arrow{d}{\sim } &  p \AR^{d+1-k-s}(\Delta)/{\mr{ann} ( p\IR^{k})}      
	\arrow{d}{\sim } \\
	\x_{\mbf{a}} p  \IR^{k}   \arrow{r}{\ \ \ \cdot \x_{\mbf{a}}^{d-s-2k}\ \ \ } &  p  {\AR^{d+1-k-s}(\mr{pyr}_{\mbf{a}} \Delta)}/ {\mr{ann} (\x_{\mbf{a}} p \IR^{k}) }
\end{tikzcd}
\]	
The horizontal map on the top being an isomorphism is equivalent to the horizontal map on the bottom being an isomorphism.
\end{proof}

\subsection{Consequences of anisotropy}
With this, we are ready to conclude the following consequence of Theorem~\ref{thm:ani}, that elevates anisotropy to absolute Hall-Laman relations:

\begin{thm}\label{thm:Halllef} We assume that the characteristic of $\fld$ is $2$ or $0$.
Assume  $\Delta$ is an IDP lattice ball or sphere of dimension $d$,  $m$ is a monomial of degree $s$ in $\fld^\ast[\Delta]$,
and  $k\le \nicefrac{(d+1-s)}{2}$.  Then, with respect to the generic linear system of parameters and the generic $\ell\in\AR^1(\Delta)$,
we have that $\AR^\ast(\Delta,\partial \Delta)$ satisfies the 
absolute Hall-Laman relations with respect to the triple $(k,m,\ell)$.	 
\end{thm}

\begin{proof}   The characteristic  $0$ case follows from
the characteristic  $2$ by arguing as in \cite[Section 3.2]{KLS}.  Hence, we can assume for the remainder of the proof that the field $\fld$ has
characteristic~$2$.
We set  $t = d+1-s-2k$ and do induction on $t \geq 0$.

STEP 1.   We assume $t =0$.   We assume  $\ell$ is a nonzero element of $\AR^1(\Delta)$.  
As noted above, the absolute Hall-Laman relations with respect to the triple $(k,m,\ell)$ are equivalent to proving that
all  $u \in \AR^{k}(\Delta,\partial \Delta)$ such that  $mu\in \AR^{k+s}(\Delta,\partial \Delta)$ is nonzero
it holds that  $mu^2\ell^{d+1-s-2k}$ is a nonzero element of $\AR^{d+1}(\Delta,\partial \Delta)$.  
Since $d+1-s-2 k = t = 0$,  it is enough to prove that $mu^2$ is a nonzero element of 
$\AR^{d+1}(\Delta,\partial \Delta)$, which is true by Theorem~\ref{thm:ani}.

\vspace{10pt}

STEP 2.   We assume $t =1$.  We set $Z= \mr{pyr}_{\mbf{a}} \Delta$.  By STEP 1, if 
$\ell$ is any element of $\AR^1(Z)$ we have that $\AR^\ast(Z,\partial Z)$ satisfies the 
absolute Hall-Laman relations with respect to the triple $(k+1,m,\ell)$.	 Hence, 
$\AR^\ast(Z,\partial Z)$ satisfies the 
absolute Hall-Laman relations with respect to the triple $(k+1, m ,\x_{\mbf{a}})$.
Using the equivalence of Part (3) and Part (4) of Lemma~\ref{lem:midred}
it follows that $\AR^\ast(\Delta ,\partial \Delta)$ satisfies the 
absolute Hall-Laman relations with respect to the triple $(k, m , h )$.

\vspace{10pt}

STEP 3.   We assume $t \geq 2$ and, by inductive hypothesis, that  
that the statement of the Theorem is true for the values $t-2$ and $t-1$.
We first prove that for generic $\ell \in \AR^1(\Delta)$ and for any $u$ in $\AR^k(\Delta,\partial \Delta)$ 
such that $mu\neq0 $ we have
\[mu\cdot\ \ell\neq 0.\]
For this, notice that because $\Delta$ is IDP and due to Poincar\'e duality
(see Equation~(\ref{eqn!map0033}))  there exists a lattice point $v$ in $\Delta$ 
such that 
\[\x_vm u\neq 0.\]
It follows from the inductive  case  for $t-1$ that the absolute Hall-Laman relations 
are true relative to the triple $(k,\x_v m,\ell)$, and hence 
\[(\x_vm) u^2 \ell^{d+1-s-2k-1} \neq 0.\] 
In particular, 
\[mu \cdot \ell\ \neq\  0\] 
since $d+1-s-2k-1 = t -1 \geq 1$. 

We set $u'=u\ell$. 
By the inductive case for $t-2$ we have that $\AR^\ast(\Delta ,\partial \Delta)$ satisfies the 
absolute Hall-Laman relations with respect to the triple $(k-1, m , \ell)$, where $\ell$
is the generic element of $\AR^1(\Delta)$.  Hence,   $mu'  \neq\  0$ implies that
\[
m(u')^2\ell^{d-s-2k-1}\ \neq \ 0.
\] 
Since  
\[
mu^2\ell^{d+1-s-2k} \ =\ m(u')^2\ell^{d-s-2k-1}
\]
it follows that 
\[
mu^2\ell^{d+1-s-2k} \ \neq \ 0. \qedhere
\]
\end{proof}

\begin{cor*}
If $\Delta$ is an IDP lattice ball or sphere of dimension $d$, and the characteristic of $\fld$ is $2$ or $0$, then 
the some Artinian reduction $\AR^\ast(\Delta)$ of  ${\fld}^\ast[\Delta]$ has the 
(relative) Lefschetz property, i.e., there exists a linear element $\ell\in\AR^1(\Delta)$ such that for all $k\leq \nicefrac{d+1}{2}$,
\[\AR^k(\Delta,\partial \Delta)\ \xrightarrow{\ \cdot \ell^{d+1-2k}\ }\ \AR^{d+1-k}(\Delta)\]
is an isomorphism. 	
\end{cor*}

\begin{proof} 
Using  Theorem~\ref{thm:Halllef} for the special case  $m=1$ 
there exists a linear element $\ell\in\AR^1(\Delta)$ such that
that $\AR^\ast(\Delta,\partial \Delta)$ satisfies the 
absolute Hall-Laman relations with respect to the triple $(k,1,\ell)$.
The result follows by applying Lemma~\ref{lem:triv}  for $p=1$ and
$\mathcal{I}^\ast = \AR^\ast(\Delta,\partial \Delta)$.
\end{proof} 

Specializing to $\Delta=P$ gives Theorem~\ref{thm:rellef}, and of course the corollary itself is just Theorem~\ref{thm:lefsphere}. We can conclude that we are left with the task of proving Theorem~\ref{thm:ani}, and move on.

\section{Kustin-Miller normalization of the volume map}\label{sec:deg}
In the present section, unless otherwise mentioned, we work
over a field $\fld$ of characteristic $2$.

In the setting  of lattice polytopes we know the canonical module thanks 
to Danilov and Stanley (it is simply the ideal generated by the interior 
lattice points of the cone \cite[Theorem 6.3.5]{Bruns-Herzog}).
We discuss  below a specific vector space isomorphism
of the top homogeneous  degree of the canonical module 
with the  field  $\widetilde{\fld}$. This will be an essential ingredient 
in our proof of Theorem~\ref{thm:ani}.  In the situation of classical algebraic geometry, 
there is a canonical such identification, which leads to a classical combinatorial formula in toric geometry \cite{Brion}. 
In our case, no such canonical identification seems to have been explored.

In other words, for lattice polytopes of dimension $d$  we will define a vector space isomorphism
$\vol: \AR^{d+1}(P,\partial P)\to\widetilde{\fld}$.  One main case of interest will be  
to understand, when $d$ is odd,  $\vol(u^2)$ as a function of $u \in  \AR^{(d+1)/2}(P,\partial P)  $ and 
the linear system of parameters $\theta_{i,j}$. The map $\vol$ allows us to give a more direct description of such objects as rational functions in the variables $\theta_{i,j}.$

The map $\vol$ is usually called the degree map, but to avoid confusion with 
the degree of a polynomial, we shall instead call this identification the \emph{volume map}.
(alluding to language used in algebraic and convex geometry \cite{KFT}).


\subsection{Normalizing the volume map}

In \cite{APP2}, we developed  a theory that, for  
Gorenstein standard graded algebras over a field of arbitrary characteristic, 
provides  a  volume normalization, that is, a vector space isomorphism from the
socle degree of its generic Artinian reduction to the field 
of rational functions  $\widetilde{\fld}= {\fld}(\theta_{i,j})$,
which is canonical in the group quotient $\widetilde{\fld}^\times/\fld^\times$. 
The last statement means that the volume normalization isomorphism 
is unique up to multiplication by a nonzero element 
of the field ${\fld}$.   In the present setting of  characteristic $2$
graded algebras associated to lattice polytopes  
we describe below a {\it unique} normalization of the volume map, 
which we call the \Defn{Kustin-Miller normalization}.

Since the  $\widetilde{\fld}$-vector space $\AR^{d+1}(P,\partial P)$ is
of dimension $1$, to determine the volume map uniquely
it suffices to give one nontrivial affine condition. In other words, 
it is enough to exhibit a nonzero element $u \in   \AR^{d+1}(P,\partial P)$ 
and state that the volume map is the unique  linear map 
\[
\vol :  \AR^{d+1}(P,\partial P) \longrightarrow \widetilde{\fld}
\]
satisfying $\vol (u) = 1$.

A good definition, of course, should  come with desirable properties, and we will discuss 
it here.  For further algebraic  justifications we refer the reader to \cite{APP2}.

\emph{The Kustin-Miller normalization:} We start with the simplest case: Recall that a lattice simplex is called unimodular 
if the associated semigroup algebra is isomorphic to a polynomial ring. In other words, affine integral combinations of the
vertices of the simplex generate the lattice.

Now, if $P$  is a unimodular lattice simplex, then we naturally have in connection to toric varieties (see \cite{Brion}) the following:

\textbf{Prototype.}  Assume $P$ is a unimodular  lattice simplex.  We set  
\[
u =    \det(\Theta_{|P})  \x_P   \in    \AR^{d+1}(P,\partial P),
\]
where  $\Theta=(\theta_{i,j})$ denotes the matrix of coefficients of the linear system of parameters
and $\Theta_{|P}$ denotes the submatrix of $\Theta$ obtained by keeping the columns indexed by 
an ordering of the elements of  $P$.  We remark that due to working in characteristic $2$, the value of 
$\det(\Theta_{|P})$ is independent of the choice of the ordering.  Then,  there exists a unique
$\widetilde{\fld}$- linear isomorphism $\vol: \AR^{d+1}(P,\partial P)\to\widetilde{\fld}$   such that $\vol(u) = 1$.
Consequently, 
\begin{equation}\label{eq:normie}
1\ =\ \vol(\x_P) \det(\Theta_{|P}).
\end{equation}

This ends up being a good definition, even in a specialization of the indeterminates $\theta_{i,j}$. For instance, if $\Sigma$ is a lattice sphere whose \Defn{facets} (the inclusion-maximal faces) 
are simply unimodular lattice simplices, then the ring obtained is simply the face ring (or Hochster-Reisner-Stanley ring), 
and the volume map defined in this way coincides with the canonical map from toric geometry.

\textbf{Simple case.} If  the boundary $\partial P$ of $P$  has a  facet $\tau$ which 
is a unimodular lattice simplex, then we obtain the desired volume 
normalization by matching the face ring picture. We set 
\[
u =  \sum_{p\in P\cap\mathbb{Z}^d}   \det(\Theta_{|\tau,p}) (\x_p\x_\tau)  \in    \AR^{d+1}(P,\partial P)
\]
where   $\Theta_{|\tau,p}$ denotes the submatrix of the matrix $\Theta$ obtained by keeping the columns indexed
by the sequence which is the concatenation of an ordering of  $\tau$ with  $p$.
Then,  there exists a unique $\widetilde{\fld}$-linear isomorphism 
$\vol: \AR^{d+1}(P,\partial P)\to\widetilde{\fld}$   such that $\vol(u) = 1$.
Consequently, 
\begin{equation}\label{eq:normie}
1\ =\ \sum_{p\in P\cap\mathbb{Z}^d} \vol(\x_p\x_\tau) \det(\Theta_{|\tau,p}).
\end{equation}

\textbf{General case.} In general, consider a flag 
\[ 
(\tau_i) =  (\tau_0, \tau_1, \dots , \tau_d) 
\]
of faces of $P$   such that $\tau_d=P$ 
and, for all $0 \leq i \leq d-1$, it holds that 
$\tau_i$ is a facet of $\tau_{i+1}$ (we call this a \Defn{full flag}).  
We call a subset  $\sigma$ of $P\cap\mathbb{Z}^d$
\Defn{coherent} with the flag  $(\tau_i)$ if,
for all $0 \leq i \leq d $, the intersection  $\sigma \cap \tau_i$
has cardinality $i+1$.  In particular, this implies that $\sigma$ has cardinality $d+1$.

We set 
\begin{equation} \label{eq:norm86712}
u = \sum_{\sigma \text{ coherent with } (\tau_i)}   \det(\Theta_{|\sigma})   \x_\sigma   \in    \AR^{d+1}(P,\partial P).
\end{equation}
By   Proposition~\ref{propos!nonvanishingofvolumenormalization}  below
$u$ is a  nonzero element of $\AR^{d+1}(P,\partial P)$. 
Hence,  there exists unique  $\widetilde{\fld}$-linear isomorphism 
$\vol: \AR^{d+1}(P,\partial P)\to\widetilde{\fld}$   such that $\vol(u) = 1$.
Consequently,
\begin{equation}\label{eq:norm2}
1\ =\ \sum_{\sigma \text{ coherent with } (\tau_i)} \vol(\x_\sigma) \det(\Theta_{|\sigma}).
\end{equation}
In Theorem~\ref{thm:transfer} we will prove that the element $u$ is independent of the choice of 
the flag $(\tau_i)$. 

\subsection{The normalization is well-defined}    \label{subs!welldefinedofolume}

\begin{thm} \label{thm:transfer} (We recall that 
the characteristic of the field $\fld$ is $2$.)  Assume $P$ is a lattice
polytope.  Then the element $u$
of  $\AR^{d+1}(P,\partial P)$ defined in 
Equation~(\ref{eq:norm86712})  is independent of the flag $(\tau_i)$ chosen, and 
it agrees with the Kustin-Miller normalization of \cite{APP2}.
\end{thm}


We shall  prove the first part
of Theorem~\ref{thm:transfer} in Subsection~\ref{subs!pfofconsist},
to convince the reader we are justified in choosing the volume normalization this way. 
More details, as well as a proof of  the second part, is provided in \cite{APP2}.

\subsection{The porcupine} Let us begin with a simple observation:

\begin{lem}
If $P$ is a $d$-dimensional lattice polytope, then the boundary of the pyramid 
over $P$ is a lattice sphere that contains $P$ as a facet. \qedhere
\end{lem}

We now give a \emph{construction}: Given the lattice polytope $P$ of dimension $d$, we want 
to obtain the \Defn{$d$-th generation porcupine} of $P$. Let us illustrate the definition 
by first defining the first and second generation.

We start by considering the pyramid over $P$. Let us call the apex point $\alpha_0$. This 
is also the \Defn{$1$st generation porcupine}.

The boundary of the porcupine has several facets, that is, maximal faces that are 
of the form $\mr{pyr}_{\alpha_{0,P}} F$, where $F$ is any facet of $\partial P$. 

Consider the pyramids over those facets, each with its own apex $\alpha_{1,F}$. We obtain 
pyramids $\mr{pyr}_{\alpha_{1,F}}\mr{pyr}_{\alpha_{0,P}} F$, each of which naturally attach to $\mr{pyr}_{\alpha_{0,P}} P$ along 
their base, the common face $\mr{pyr}_{\alpha_{0,P}} F$, resulting in a lattice ball consisting of $\mr{pyr}_{\alpha_{0,P}} P$, and the $\mr{pyr}_{\alpha_{1,F}}\mr{pyr}_{\alpha_{0,P}} F$.
This lattice ball is the \Defn{$2$-nd generation porcupine}.

\begin{dfn}
Consider a lattice polytope $P$ of dimension $d$. Let $P^{(i)}$ denote the collection of $i$-dimensional 
faces of $P$, and let $P^{(\ge i)}$ denote the faces of dimension at least $i$ (which includes $P$ itself if $i\le d$). 

The \Defn{$k$-th generation porcupine} over $P$, $1\le k\le d+1$, is the $(d+1)$-dimensional lattice ball 
\[\bigcup_{i \ge d-k+1}\  \bigcup_{(G_i \in P^{(i)} \subset G_{i+1} \in P^{(i+1)} \subset \cdots  \subset P \in P^{(d)})\ \text{flag}} 
\mr{pyr}_{\alpha_{d-i,G_i}} \mr{pyr}_{\alpha_{d-i+1,G_{i+1}}}\cdots  \mr{pyr}_{\alpha_{0,P}} G_i 
\]
where \[\mr{pyr}_{\alpha_{d-i,G_i}} \mr{pyr}_{\alpha_{d-i+1,G_{i+1}}}\cdots  \mr{pyr}_{\alpha_{0,P}} G_i\] 
is attached to \[\mr{pyr}_{\alpha_{d-i+1,G_{i+2}}} \mr{pyr}_{\alpha_{d-i+2,G_{i+2}}}\cdots  \mr{pyr}_{\alpha_{0,P}} G_{i+2} \]
along the common face \[\mr{pyr}_{\alpha_{d-i+1,G_{i+1}}}\cdots  \mr{pyr}_{\alpha_{0,P}} G_i.\]	 
We will denote this object by $\mr{porc}_k P$.
\end{dfn}	

\subsection{Balancing and locality}
We need a further small lemma, which describes what is sometimes called a \Defn{balancing identity}:

\begin{lem}\label{lem:balancing}
Consider any $1\leq i\leq d+1$, and $\x_I$ a monomial of degree $d$ of $\AR^\ast(X)$, where $X$ 
is some $d$-dimensional lattice complex. Then in $\AR^\ast(X)$ we have 
\[\sum_{p \text{ lattice point in } X} \theta_{i,p}\vol(\x_{I}\x_p)\ =\ 0 .\]
\end{lem}

\begin{proof}
The identity in the statement of the lemma  arises because
when constructing the Artinian reduction $\AR^\ast(X)$
we quotient by the linear elements 
\[   \theta_i=\sum_{p \in P \cap \mathbb{Z}^d} \theta_{i,p} \x_p.  \qedhere\]
\end{proof}

We now focus on Theorem~\ref{thm:transfer}. 
An intermediate word is useful: We are working to establish the 
volume map on $\AR^\ast(P,\partial P)$. The idea now is to think of P as 
a facet in a larger lattice sphere $M$	of the same dimension, and use properties 
of that larger space. Let us call this a locality principle: 
The volume is locally defined on $P$, but is consistent with the algebraic 
structure of the larger space. This is justified by the following Lemma.

\begin{lem}\label{lem:local}[Locality principle]
Consider $\Delta$ a lattice ball of dimension $d$, and $X \supset \Delta$ a lattice sphere or ball 
of the same dimension that contains $\Delta$ as subcomplex. Then the natural inclusion of semigroup algebras induces an 
injective map
\[
\AR^\ast(\Delta,\partial \Delta) \ \longhookrightarrow\  \AR^\ast(X,\partial X).
\]
\end{lem}

\begin{proof}
We begin by defining the lattice complex $X-\Delta$ as induced by those 
facets of $X$ not in $\Delta$. We have a short exact sequence
\[0\ \longrightarrow\ \widetilde{\fld}^\ast[\Delta,\partial \Delta]\ \longrightarrow\ 
\widetilde{\fld}^\ast[X,\partial X]\ \longrightarrow\ \widetilde{\fld}^\ast[(X-\Delta) 
\cup \partial X,\partial X] \ \longrightarrow\ 0.\]
All three modules in this sequence are Cohen-Macaulay, see \cite{BBR}, hence the sequence 
remains exact after Artinian reduction. The claim follows. 
\end{proof}

\begin{prp}   \label{propos!nonvanishingofvolumenormalization} 
Assume  $\Delta$ is a lattice ball of dimension $d$ and 
$(\tau_i)$ is a flag of faces of $\Delta$,
where the notion of flag is as defined above
hence $(\tau_i)$  is a maximal flag. Then the element $u$
of  $\AR^{d+1}(\Delta,\partial \Delta)$ defined in 
Equation~(\ref{eq:norm86712}) is nonzero.
\end{prp}

\begin{proof}  Combining the computations in  Subsection~\ref{subs!pfofconsist}
with Lemma~\ref{lem:local} 
the claim is reduced  to the case of a unimodular lattice simplex where it
is obvious.  
\end{proof}


\subsection{Proof of consistency}  \label{subs!pfofconsist}

We return to prove the well-definedness of the Kustin-Miller normalization. 
A direct   proof is rather uninformative, but 
we can use an indirect argument without getting our fingers dirty, and that informs what is really happening and how 
it connects to the naturality of the definition. However, for the reader preferring a down 
to earth and explicit proof, we refer them to the Appendix~\ref{sec:KMrevisit}.

\begin{proof}[\textbf{Proof of  the first part of \textbf{Theorem~\ref{thm:transfer}}}]
We recall that  we work over a field $\fld$   of   characteristic   $2$.
Moreover, by flag of faces of a lattice complex we will always mean a full flag in this proof.

Assume  $\Sigma$  is a $d$-dimensional lattice sphere or lattice ball.
We denote by $L(\Sigma)$ the set of lattice points of $\Sigma$. 
Given a  flag $T=(\tau_i)$ of faces of $\Sigma$
we set
\begin{equation}\label{eq:normolo}
	\mathcal{N}[T]\coloneqq\sum_{\sigma \text{ coherent with } (\tau_i)}  \det(\Theta_{|\sigma}) \x_\sigma  \in    \AR^{d+1}(\Sigma,\partial \Sigma).
\end{equation}

More generally,  given a finite sequence   $U=(U_0, \dots , U_d)$   such that 
$ \emptyset \not=  U_i \subset L(\Sigma)$  for all $i$ we set
\begin{equation}
	\mathcal{N}[U] \coloneqq  \sum_{r}   \det(\Theta_{|r}) \x_r   \in    \AR^{d+1}(\Sigma,\partial \Sigma),
\end{equation}
where the sum is over all  $r =(r_0, \dots , r_d)$ with  $r_i \in U_i$ for all $0 \leq i \leq d$. 

Since the field  $\fld$ has characteristic $2$, it is clear that 
if there exist  $i \not= j$ with $U_i = U_j$ then 
\begin{equation}    \label{eqn!eawreno1}
	\mathcal{N}[U] = 0
\end{equation}
due to pairwise cancellation of terms. 
Moreover, if  for some $i$ it holds that   $U_i$ is the disjoint union of two nonempty subsets 
$W_1, W_2$ then
\begin{equation}    \label{eqn!eawreno2}
	\mathcal{N}[U] = \mathcal{N}[Z_1] + \mathcal{N}[Z_2],
\end{equation}
where, for $1 \leq t \leq 2$,  $Z_t$ is obtained from $U$ by replacing
$U_i$ with $W_t$.
\\

Assume now $P$ is a lattice polytope. We denote by 
$M = \partial \mr{porc}_d P$  the boundary  of the $d$-th generation porcupine of $P$. 
It is a lattice sphere, and therefore the resulting algebra is Gorenstein \cite{BBR}.

Notice that it has several facets which are unimodular simplices, in particular those of the form
\[\mr{pyr}_{\alpha_{d-1,S}}\dots\mr{pyr}_{\alpha_{1,F}}\mr{pyr}_{\alpha_{0,P}} v,\]
where $v$ is a vertex of $P$ and $S\subset\dots \subset G\subset F$ is a  flag. We already 
know the natural normalization on them, and it states 
\[1\ =\  \vol( \mathcal{N}[(v, \mr{pyr}_{\alpha_{0,P}} v, \mr{pyr}_{\alpha_{1,F}}\mr{pyr}_{\alpha_{0,P}} v, \dots, 
\mr{pyr}_{\alpha_{d-1,S}}\dots\mr{pyr}_{\alpha_{1,F}}\mr{pyr}_{\alpha_{0,P}} v)]).\]
It suffices to prove that the normalization of Equation~\eqref{eq:norm2} is 
consistent with the normalization of the unimodular simplex.
Hence, we want to prove that for an arbitrary flag $(\tau_i)$ 
ending with $\tau_d=P$, the sum $\mathcal{N}[T]$ equals such a term.  \\

All equalities in the following  are in $ \AR^{d+1}(M,\partial M) =  \AR^{d+1}(M, \emptyset) =   \AR^{d+1}(M) $.

For simplicity of notation  we set  $\mr{p} = \mr{pyr}$ and,  for $i \geq 0$,  $w_i = a_{i, \tau_{d-i}}$.
We also  set  
\[
H_{(1,0)} =  \mathcal{N}[T]   = \mathcal{N}[(\tau_0,\tau_1,\dots,\tau_{d-1},\tau_d)],  \quad  
H_{(1,1)} = \mathcal{N}[(\tau_0,\tau_1,\dots,\tau_{d-1},  \mr{p}_{w_0} \tau_{d-1})],
\]
and for $2 \leq i \leq d-1$ we set 
\[
H_{(1,i)} =  \mathcal{N}[(\tau_0,\tau_1,\dots,       \tau_{d-i},   \mr{p}_{w_0} \tau_{d-i}, 
\mr{p}_{w_0} \tau_{d-i+1},   \dots ,     \mr{p}_{w_0} \tau_{d-1}  ) ].
\]

We have, by definition, $\mathcal{N}[T] = H_{(1,0)}$.  Using  Lemma~\ref{lem:balancing}  we get $H_{(1,0)} = H_{(1,1)}$. 
Combining Equations~(\ref{eqn!eawreno1}) and~(\ref{eqn!eawreno2}) we 
get that  $H_{(1,i)} = H_{(1,i+1)}$ for all $1 \leq i \leq  d-2$. Hence
\[
\mathcal{N}[T] = H_{(1,d-1)} =   \mathcal{N}[(\tau_0, \mr{p}_{w_0} \tau_{0},  \mr{p}_{w_0} \tau_{1},  
\mr{p}_{w_0} \tau_{2}   ,\dots,  
\mr{p}_{w_0} \tau_{d-1}  ) .
\]               
We set  
\[
H_{(2,0)} =  H_{(1,d-1)},  \quad  H_{(2,1)} = \mathcal{N}[(\tau_0, \mr{p}_{w_0} \tau_{0},  \mr{p}_{w_0} \tau_{1}   ,\dots,  
\mr{p}_{w_0} \tau_{d-2},   \mr{p}_{w_1}(  \mr{p}_{w_0} \tau_{d-2}) )]
\]
and for $2 \leq i \leq d-2$ we set  
\[
H_{(2,i)} =  \mathcal{N}[(\tau_0,    \mr{p}_{w_0} \tau_{0}  ,\dots,       ,  \mr{p}_{w_0} \tau_{d-i},   
\mr{p}_{w_1}(  \mr{p}_{w_0} \tau_{d-i}),  \mr{p}_{w_1}(  \mr{p}_{w_0} \tau_{d-i+1}),
\dots ,     \mr{p}_{w_1}(  \mr{p}_{w_0} \tau_{d-2})   ) ].
\]
Using  Lemma~\ref{lem:balancing}  we get $H_{(2,0)} = H_{(2,1)}$. 
Combining Equations~(\ref{eqn!eawreno1}) and~(\ref{eqn!eawreno2}) we get that  $H_{(2,i)} = H_{(2,i+1)}$ for all $1 \leq i \leq  d-3$.
Hence
\[
\mathcal{N}[T] = H_{(2,d-3)} =   \mathcal{N}[(\tau_0, \mr{p}_{u_0} \tau_{0}, q_2, q_3   ,\dots,   q_d  ) ],
\]
where   $q_{i+2} =   \mr{p}_{w_1} ( \mr{p}_{w_0} \tau_{i} )$   for all $0 \leq i \leq d-2$.

Continuing inductively,  we get 
\[
\mathcal{N}[T] = \ \mathcal{N}[(z_0, z_1, \dots  , z_d)], 
\]
where   $ z_0 =   \tau_0$  and   $z_{i+1} =   \mr{p}_{w_i}(z_i)$ for all $0 \leq i \leq d-1$.
This  finishes the proof  of  the first part of  Theorem~\ref{thm:transfer}.   \end{proof}


\section{Differential equations for the volume map}\label{sec:proofani}

Henceforth, the volume of a lattice polytope,  lattice ball or lattice sphere, is considered 
to be normalized with respect to the Kustin-Miller normalization map $\vol$ introduced in 
Subsection~\ref{sec:deg}. The anisotropy 
property we wish to prove, that is, Theorem~\ref{thm:ani}, is dependent of our choice 
of linear system of parameters $\Theta \x$, so we will consider $\vol$ as a rational function 
in the variables $\theta_{i,j}$. This allows us to formulate the following auxiliary lemma on the 
way to anisotropy.

For a finite sequence of lattice points $w=(w_1, \dots , w_s)$ of $P$, 
we denote by $|w|$ their sum 
which is an element of the lattice. In other words, $|w|=w_1 + \dots + w_s$.

\begin{lem}\label{lem:diffid} 
Assume  $k \geq 1, j \geq 0$  and that the field $\fld$ has characteristic $2$.
Assume $\Delta$ is a $(2k-1+j)$-dimensional lattice ball or sphere and $P$ is a 
facet of $\Delta$. Let $F=(w_1,\ldots,w_{2k+j})\subset (P\cap\mathbb{Z}^d)^{2k+j}$ be  
a sequence of $2k+j$ lattice points of $P$, and 
$\sigma=(\sigma_1,\ldots,\sigma_{j})\subset (P\cap\mathbb{Z}^d)^{j}$
a possibly empty, if $j=0$, sequence of $j$ lattice points of $F$. Assume that
\[
           |F|-|\sigma|=\sum_{i=1}^{2k+j}  w_i - \sum_{i=1}^{j} \sigma_i =2|G|
\] 
for some sequence $G = (g_1, \dots, g_{k+j})$ of $k+j$   lattice points  of $P$. Consider $u\in \fld^k[\Delta, \partial \Delta]$, or $u\in \fld^k[\Delta]$ if $\sigma$ is not contained in a face of $\partial \Delta$. Then we have 
\[
      \partial_F \vol_\sigma(u^2)\ =\ (\vol_\sigma(u \cdot\x_G ))^2.
\]
\end{lem}

\vspace{20pt}

Here,  $\partial_F$ is the $(2k+j)$-th order  differential operator obtained as the composition  of differentiation 
after the variables 
$\theta_{1,w_1}, \theta_{2,w_2}, \theta_{3,w_3},  \dots  ,  \theta_{2k+j,w_{2k+j}}$.
 In other words, 
\[
    \partial_F  =  \frac { \partial^{2k+j}}
          {\partial \theta_{1,w_1} \;  \partial \theta_{2,w_2} \; \partial \theta_{3,w_3}\; \dots 
            \; {\partial \theta_{ 2k+j,w_{2k+j}}}}.
\]
Moreover, by definition, $\; \x_\sigma = \prod_{i=1}^j \x_{\sigma_i}$,  \;
$\x_G = \prod_{i=1}^{k+j} \x_{g_i}$ and
\[
    \vol_{\sigma} (w) = \vol ( \x_{\sigma} \cdot w)
\]
for all $w \in \AR^{2k}(\Delta,\partial \Delta)$.

We postpone the proof of this lemma to Subsection~\ref{sub:diffid}, where we will derive 
it from the key identity of Parseval-Rayleigh type. 

From this differential identity, anisotropy follows at once.

\begin{proof}[\textbf{Proof of Theorem~\ref{thm:ani}}]
Consider $u\neq 0$ of degree $k\le \frac{d+1}{2}$ in $\AR^k(\Delta, \partial \Delta)$. 
The pairing 
\[
   \AR^k(\Delta, \partial \Delta)\ \times\ \AR^{d+1-k}(\Delta)    \ \longrightarrow\ \AR^{d+1}(\Delta, \partial \Delta)
\] 
is nondegenerate. Since $\Delta$ is IDP, hence each $P\in \Delta$ is IDP, there exist $P\in \Delta$, and sequences of 
lattice points $\sigma=(\sigma_1,\ldots,\sigma_{d+1-2k}) \subset (P\cap\mathbb{Z}^d)^{d+1-2k}$  
and $G= ( g_1,\ldots,g_k ) \subset (P\cap\mathbb{Z}^d)^k$ such that $\x_\sigma \cdot\x_G \cdot u$ is not zero.
Consequently, $\vol_{\sigma} ( u \x_G )$ is not zero.


We set  $F$ to be the concatenation of the family $G$, another copy of $G$, and finally $\sigma$.
Using Lemma~\ref{lem:diffid},  $\vol_\sigma(u^2)$ is not zero as
 \[
          \partial_{F}\vol_{\sigma}(u^2)\ =\ (\vol_\sigma(u \cdot\x_G ))^2.
\] 
Hence $u^2$ is also not zero.
\end{proof}

\section{Level properties and another Lefschetz/anisotropy theorem}\label{sec:anilevel}

In this section we consider Theorem~\ref{thm:level}. The algebraic statement is the following:

\begin{thm}\label{thm:levellef}
Assume that the field $\fld$ has characteristic $2$ or $0$.
Then the semigroup algebra of an IDP lattice polytope $P$ of dimension $d$ with $\cone^\circ(P)\cap (\mathbb{Z}^d\times \mathbb{Z})$ generated 
at height $\leq j$ satisfies an almost Lefschetz theorem: We have for some suitable Artinian reduction of $\fld^\ast[P]$ that
\[\AR^k(P)\ \xrightarrow{\ \cdot \ell^{d+1-j-2k}\ } \AR^{d+1-j-k}(P)\]
is an injection for some $\ell$ in $\AR^1(P)$ and every $k\leq \frac{d+1-j}{2}$.
\end{thm}

It is the consequence of two lemmata.

Before we go to the statements, we need the concept of an \Defn{interior simplex}, that is, a collection of lattice points that 
do not lie in a strict face of the polytope.  More precisely, assume  $P$ is an IDP lattice polytope of dimension $d$ and
$1 \leq j \leq d+1$. We denote by 
\[
     \mathbf{I}_{j} (P)  
\]
the set of $j$-tuples $v = (v_1,  \dots , v_j)$, with $v_i$ a lattice point of $P$ for all $1 \leq i \leq j$,
such that there is no facet of $F$ of $P$ with 
\[
     \{ v_1,  \dots , v_j \} \subset F.
\]
We remark that $v_a = v_b$ for $a \not= b$ is allowed. 
It is clear that if $\sigma \in  \mathbf{I}_{j} (P)$, then it holds that  $\x_\sigma \in  \fld^\ast[P, \partial P]$. 
Since $\AR^\ast(P, \partial P)$ is an $\AR^\ast(P)$-module, we set 
\[
   {\mr{ann}\, \x_\sigma}  =  \{  r \in  \AR^\ast(P) :  r [\x_{\sigma}] = 0_{\AR^\ast(P, \partial P)}  \},
\]
where $[\x_{\sigma}]$ denote the class of $\x_\sigma$ in  $\AR^\ast(P, \partial P)$. We have that
${\mr{ann}\, \x_\sigma}$ is a homogeneous ideal of $\AR^\ast(P)$.

\subsection{Two lemmata for interior simplices}

The first lemma is a version of partition of unity in the given situation (see also \cite{AY} for the case of face rings).
We will use the notation 
\[
      \AR^{\leq d+1-j}(P)  \  = \  \bigoplus_{i=0}^{d+1-j} \AR^i(P) \  \subset \AR^\ast(P).
\]

\begin{lem}  \label{lem!nfksgfd}  Assume the field $\fld$ has arbitrary characteristic  and
$P$ is an IDP lattice polytope  of dimension $d$ 
with $\fld[P,\partial P]$ generated at height $\leq j$.
Then the natural map
\[
       \AR^\ast(P)\ \longrightarrow \bigoplus_{\sigma \in \mathbf{I}_{j} (P)}\bigslant{\AR^\ast(P)}{\mr{ann}\, \x_\sigma}
\]
defined by 
\[
      u  \ \mapsto  \  ([u]_{\sigma} )_{\sigma \in \mathbf{I}_{j} (P)},
\] 
where $[u]_{\sigma}$  denotes the class of $u$ in ${\AR^\ast(P)}/{\mr{ann}\, \x_\sigma} $,
induces by restriction an injection
 \[
    \AR^{\leq d+1-j}(P) \ \longhookrightarrow\ \bigoplus_{\sigma \in \mathbf{I}_{j} (P)}\bigslant{\AR^\ast(P)}{\mr{ann}\, \x_\sigma}.
\]
\end{lem}

\begin{proof}
It is enough to prove that if  $t$ satisfies  $1 \leq t \leq  d+1-j$ and 
$0 \not= u \in  \AR^t(P)$
then there exists $\sigma \in   \mathbf{I}_{j} (P)$ such that the element 
$u \cdot \x_\sigma$ of $\AR^\ast(P, \partial P)$ is nonzero.
This follows as by  \cite[Theorem~6.31]{BG}, there exists  a homogeneous generating  set $q_1, \dots ,q_s$  
for the   $\AR^\ast(P)$-module $ \AR^\ast(P,\partial P)$ with  $\deg (q_i) \leq j$
for all $1 \leq i \leq s$. 
\end{proof}

The second is a Lefschetz type fact:
\begin{lem}[Lefschetz properties relative to an interior simplex]\label{lem:levellef}  
Assume that the field $\fld$ has characteristic $2$ or $0$,
 $P$ is an IDP lattice $d$-polytope and   $\sigma  \in   \mathbf{I}_{j} (P)$. We set 
 \[\mathcal{C}^\ast~:=~\bigslant{\AR^\ast(P)}{\mr{ann}\,  \x_\sigma}\]
 and assume we work in the generic Artinian reduction.
Then the following   Lefschetz type  property holds: We have that there exists
$\ell \in \mathcal{C}^1$ such that for all $k\leq \frac{d+1-j}{2}$ the 
map induced by multiplication with $\ell^{d+1-j-2k}$, 
\[ 
        \mathcal{C}^k \  \xrightarrow{\ \cdot \ell^{d+1-j-2k}\ }\   \mathcal{C}^{d+1-j-k},
\]
is a bijection. 
\end{lem}

\subsection{Interlude: The ideas behind Lemma~\ref{lem:levellef}}

We go quickly over the main ideas of the proof of this lemma, which is analogous to the ideas of Section~\ref{sec:ani}.
In fact, $\mathcal{C}$ is a Gorenstein ring, see Remark~\ref{rem:annGorenstein} and socle lives in degree $d+1-j$, and satisfies a Hall-Laman type property analogous to Theorem~\ref{thm:Halllef}, applied to the pairing
\[ \mathcal{C}^k \times  \mathcal{C}^{d+1-j-k}  \longrightarrow  
 \mathcal{C}^{d+1-j}.\]
Of course, we can understand this entirely within $\AR^\ast(P,\partial P)$ via the map 
\begin{equation}\label{eq:sl}
	\begin{array}{ccc}
		\AR^q(P)& \longrightarrow &\AR^{q+j}(P,\partial P) \\
		a	& {{\xmapsto{\ \ \ \ }}} &\x_\sigma a
	\end{array}
\end{equation}

The corresponding Hall-Laman statement, that immediately implies Lemma~\ref{lem:levellef}, is therefore:

\begin{lem}\label{lem:Hallleflev} We assume that the characteristic of $\fld$ is $2$ or $0$.
	Assume  $\Delta$ is an IDP lattice ball or sphere of dimension $d$,  $m$ is a monomial of degree $s$ in $\fld^\ast[\Delta,\partial \Delta]$,
	and  $k\le \nicefrac{(d+1-s)}{2}$.  Then, with respect to the generic Artinian reduction and the generic $\ell\in\AR^1(\Delta)$, we have that every $u \in \AR^k(\Delta)$ such that $mu\neq 0$ in $\AR^\ast(\Delta,\partial \Delta)$ also satisfies
	\[mu^2 \ell^{d+1-s-2k}\ \neq 0\ \text{in}\  \AR^\ast(\Delta,\partial \Delta).\]		 
\end{lem}

We skip going over the proof of this statement in detail, and simply note that Lemma~\ref{lem:Hallleflev} follows as in Section~\ref{sec:ani} from the pyramidal lifting reduction and an analogous anisotropy theorem:

\begin{thm}\label{thm:levelani}   
	If $X$ is an IDP lattice ball or sphere of dimension $d$, and the characteristic of $\fld$ is $2$ or $0$. Consider 
	an element $u\in\AR^{k}(X)$ with respect to the generic Artinian reduction.
	If $m \in \fld^\ast[X,\partial X]$ is a monomial of degree $\le d+1-2k$ such that $m u$ is nonzero in
	$\AR^{k+ \deg(m)}(X,\partial X)$, then 
	\[m u^2\ \neq\ 0, \]
	in other words   $m u^2$ is a nonzero element of $\AR^{2k+ \deg(m)}(X,\partial X)$.
\end{thm}

\begin{proof}    
This follows as Theorem~\ref{thm:ani} from Lemma~\ref{lem:diffid}.
\end{proof}

For completeness, the pyramidal lifting reduction in this case is stated as follows:

\begin{lem}\label{lem:midredlev}
	Consider a lattice ball $\Delta$ of dimension $d$. 
	
	Let $\IR^\ast$ be a nonzero graded ideal of  $\AR^\ast(\Delta)$. 
	
	We also consider the induced graded ideal $ \x_{\mbf{a}}\IR^\ast$ of $\AR^{\ast}(\mr{pyr}_{\mbf{a}} \Delta,\Delta)$.
	
	We also  set $h=\x_{\mbf{a}}-\overline{h}$, a homogeneous polynomial $p\in \fld^\ast[\Delta,\partial \Delta]$ of degree $s$ and a $k< \frac{d+1-s}{2}$.
	Then the following two are equivalent:
	\begin{compactenum}
		\item [(1)]  The pairing 
		\begin{equation*}
			\begin{array}{rccccc}
				& p\x_{\mbf{a}} \AR^k(\mr{pyr}_{\mbf{a}} \Delta)& \times & p\x_{\mbf{a}} \AR^k(\mr{pyr}_{\mbf{a}} \Delta)
				& \longrightarrow &\ \AR^{d+2}(\mr{pyr}_{\mbf{a}}\Delta,\partial \mr{pyr}_{\mbf{a}} \Delta) \\
				& p\x_{\mbf{a}}	a	&		& p\x_{\mbf{a}} b& {{\xmapsto{\ \ \ \ }}} &\ p (\x_{\mbf{a}}a) (\x_{\mbf{a}}b)\x_{\mbf{a}}^{d-s-2k}
			\end{array}
		\end{equation*}
		is nondegenerate on $p\x_{\mbf{a}} \IR^k\subset p\x_{\mbf{a}} \AR^k(\mr{pyr}_{\mbf{a}} \Delta) \subset \AR^{k+s+1}(\mr{pyr}_{\mbf{a}}\Delta,\partial \mr{pyr}_{\mbf{a}} \Delta)$.
		\item [(2)]  The pairing 
		\begin{equation*}
			\begin{array}{rccccc}
				& p \AR^k(\Delta)& \times & p \AR^k(\Delta)
				& \longrightarrow &\ \AR^{d+1}(\Delta,\partial \Delta) \\
				& p	a	&		& p b& {{\xmapsto{\ \ \ \ }}} &\ pab h^{d+1-s-2k}
			\end{array}
		\end{equation*}
		is nondegenerate on $p \IR^k \subset p \AR^k(\Delta) \subset \AR^{k+s}(\mr{pyr}_{\mbf{a}}\Delta,\partial \mr{pyr}_{\mbf{a}} \Delta)$.
	\end{compactenum}
\end{lem}

\begin{proof}
Consider the square
\[\begin{tikzcd}[column sep=5em]
	p \IR^{k}    \arrow{r}{\ \ \ \ \cdot h^{d+1-s-2k}\ \ \ \  } \arrow{d}{\sim } &  p \AR^{d+1-k-s}(\Delta)/{\mr{ann} ( p\IR^{k})}      
	\arrow{d}{\sim } \\
	\x_{\mbf{a}} p  \IR^{k}   \arrow{r}{\ \ \ \cdot \x_{\mbf{a}}^{d-s-2k}\ \ \ } &  p  {\AR^{d+1-k-s}(\mr{pyr}_{\mbf{a}} \Delta)}/ {\mr{ann} (\x_{\mbf{a}} p \IR^{k}) }
\end{tikzcd}
\]
where, as in the proof of Lemma~\ref{lem:midred}, the righthandside is considered as quotients of $\AR^{\ast}(\Delta)$ and $\AR^{\ast}(\mr{pyr}_{\mbf{a}} \Delta)$, respectively.
As in the proof of Lemma~\ref{lem:midred}, Case 4, the top horizonzal map is an isomorphism if and only if the bottom map is.	
\end{proof}

\begin{proof}[\textbf{Proof of Lemma~\ref{lem:Hallleflev}}]
	As in the proof of Theorem~\ref{thm:Halllef}, the case $d+1-s-2k=0$ is immediate from Theorem~\ref{thm:levelani}. The case $d+1-s-2k=1$ follows by combining Theorem~\ref{thm:levelani} and Lemma~\ref{lem:midredlev}.
	
	If $d+1-s-2k>1$, use the IDP property to see that there is a lattice point $v$ such that $\x_vm u \neq0$. Hence by induction, for generic $\ell$ we have $\ell u\neq 0$. Use induction on $d+1-s-2k$ to establish the claim.	
\end{proof}

This finishes the proof of Lemma~\ref{lem:levellef}.

\subsection{Returning to Theorem~\ref{thm:levellef}} 

So far, the arguments are analogous to those of Section~\ref{sec:ani}. Now we return to what is specific to this section.

\begin{proof}[\textbf{Proof of Theorem~\ref{thm:levellef}}]

We will again abuse the well-known fact that the existence  of a  degree $1$ homogeneous
element $ \ell$  such that multiplication by $\ell^{d+1-j-2k}$ is injective
is equivalent to that for the generic degree $1$ homogeneous
element $ \ell$  the multiplication by $\ell^{d+1-j-2k}$ induces an injective map.

For $\sigma  \in \mathbf{I}_{j} (P)$ we set \[\mathcal{C}^\ast_{\sigma}~=~\bigslant{\AR^\ast(P)}{\mr{ann}\,  \x_\sigma} .\]
Consider the commutative diagram 
\[\begin{tikzcd}[column sep=5em]
      \AR^{k}(P)  \arrow{r}{\ \ \ \ \cdot \ell^{d+1-j-2k}\ \ \ \  }   \arrow{d}{  }   &  \AR^{d+1-j-k}(P)   \arrow{d}{ }   \\
      \bigoplus_{{\sigma \in \mathbf{I}_{j} (P)}}  {\mathcal{C}^k_{\sigma}}  \arrow{r}{\ \ \ \ \cdot \ell^{d+1-j-2k}\ \ \ \  } &  
                       \bigoplus_{{\sigma \in \mathbf{I}_{j} (P)}}  {\mathcal{C}^{d+1-j-k}_{\sigma}} 
      \end{tikzcd}
\]	
where both horizontal maps are multiplication by $\ell^{d+1-j-2k}$ and the vertical
maps are as in the statement of   Lemma~\ref{lem!nfksgfd}, hence they are injective.
By Lemma~\ref{lem:levellef}  the lower horizontal map is injective for Zariski general $\ell$. 
Since both the lower  horizontal map and the left vertical map are injective for general $\ell$, we get 
by elementary set theory that  the upper  horizontal map is also injective for general $\ell$.  This finishes
the proof of  Theorem~\ref{thm:levellef}.   \end{proof}

\begin{rem}
	The results of this section extend (with similar proof) to lattice balls and spheres in which every face satisfies 
	the assumptions of Theorem~\ref{thm:levellef}.
	
	\begin{thm}\label{thm:levellefball}
		Assume that the field $\fld$ has characteristic $2$ or $0$. Consider a lattice ball $\Delta$ of dimension $d$ such that the ideal $\fld^\ast[\Delta,\partial \Delta]$ is generated 
		in degree $\leq j$. Then we have, in a suitable Artinian reduction, that
		\[\AR^k(\Delta)\ \xrightarrow{\ \cdot \ell^{d+1-j-2k}\ } \AR^{d+1-j-k}(\Delta)\]
		is an injection for some $\ell$ in $\AR^1(\Delta)$ and every $k\leq \frac{d+1-j}{2}$.
	\end{thm}
\end{rem}

\section{The Parseval-Rayleigh identities for complexes and differential identities}\label{sec:PC}

We now provide a proof of the Parseval-Rayleigh identities, and then conclude the differential identities from them.

\subsection{The Parseval-Rayleigh identity for lattice balls and spheres}\label{sec:PC2}

\begin{lem}[The Parseval-Rayleigh identity]\label{lem:parsevalcomplex}
	For a lattice $d$-ball or sphere $\Delta$, in $\AR^\ast(\Delta,\partial \Delta)$ over characteristic $2$, and $\sigma$ a family of lattice points in $\Delta$ and for $d+1+\#\sigma$ even, we have
	\begin{equation}
		\vol (\x_\sigma u^2)\ =\ \sum_{ \beta \in (\Delta\cap \mathbb{Z}^d)^{d+1}}  \vol (u\cdot\x_{\frac{\sigma+\beta}{2}})^2 \theta^\beta
	\end{equation}
	for all $u \in \AR^{\nicefrac{d+1-\#\sigma}{2}}(\Delta,\partial \Delta)$, or $u \in \AR^{\nicefrac{d+1-\#\sigma}{2}}(\Delta)$ if $\sigma$ is an interior simplex. 	
\end{lem}

Here, $\theta^\beta \coloneqq\prod_{1 \leq i \leq d+1} \theta_{i,\beta_i}$. We further remind ourselves: $\x_{\frac{\alpha+\beta}{2}}$ is naturally $0$ if $\frac{\alpha+\beta}{2}$ is not an element of $\cone(\Delta)\cap (\mathbb{Z}^d\times \{d+1\})$, which is in particular the case if $\alpha$ and $\beta$ do not lie in a common face of $\cone (\Delta)$. Moreover, $\Delta\cap \mathbb{Z}^d$ is identified with $\cone (\Delta)\cap (\mathbb{Z}^d\times \{1\})$, so that the index of the summands $\beta$ corresponds to $d+1$ lattice points of the latter.

It suffices to prove Lemma~\ref{lem:parsevalcomplex} for monomials, because, if $u= \sum_a\lambda_a \x_a$, then 
\begin{align*}
	\vol(\x_\sigma u^2) \ &= \ \vol\left(\x_\sigma\left(\sum_a\lambda_a \x_a\right)^2\right) \ = \  \sum_a  \lambda_a^2\ \vol(\x_\sigma\x_a^2) \\
	&= \  \sum_a  \lambda_a^2\ \sum_{ \beta \in (\Delta\cap \mathbb{Z}^d)^{d+1}\times \{1\}}   \vol(\x_a\cdot\x_{\frac{\sigma+\beta}{2}})^2 \theta^\beta \\
	&= \  \sum_{ \beta \in (\Delta\cap \mathbb{Z}^d)^{d+1}\times \{1\}} \vol\left(\sum_a  \lambda_a\  \x_a\cdot\x_{\frac{\sigma+\beta}{2}}\right)^2 \theta^\beta \\
	&= \ \sum_{ \beta \in (\Delta\cap \mathbb{Z}^d)^{d+1}\times \{1\}}  \vol (u\cdot\x_{\frac{\sigma+\beta}{2}})^2 \theta^\beta.
\end{align*}

Hence, we are left with proving Theorem~\ref{thm:parsevalcomplex}, which we restate here for convenience.

\begin{thm*}
	For a lattice $d$-ball or sphere $\Delta$, and $\alpha$ a lattice point of $\cone^\circ(\Delta)\cap (\mathbb{Z}^d\times \{d+1\})$, we have in $\AR^\ast(\Delta,\partial \Delta)$ over characteristic $2$
	\begin{equation}\label{eq:idc}
		\vol (\x_\alpha) \   =\    \sum_{ \beta \in (\Delta\cap \mathbb{Z}^d)^{d+1}}   \vol (\x_{\frac{\alpha+\beta}{2}})^2   \theta^\beta.
	\end{equation}
\end{thm*}

\begin{proof}
Notice that by locality (Lemma~\ref{lem:local}), we may assume that $\Delta$ is a sphere: If $\Delta$ is a ball, consider instead the union of $\mr{pyr}_{\mbf{a}} \partial \Delta$ and $\Delta$, which are identified along the common boundary $\partial \Delta$.

We may also assume that $\Delta$ has one facet that is a unimodular simplex, say $S$: There is at least one facet $F$ that does not contain $\alpha$. Remove $F$ from $\Delta$ and $\mr{porc}_d F$, and identify the remainders along the common boundary $\partial F$.

Now, we prove the result by establishing three facts:

\begin{compactenum}[(1)]
\item It is true for some $\overline{\alpha}$.
\item Find a system of linear equations that, together with the Kustin-Miller normalization, determine $\vol (\x_\bullet)$ uniquely.
\item The linear equations that are satisfied for $\vol (\x_\bullet)$ are satisfied for \[\sum_{ \beta \in (\Delta\cap \mathbb{Z}^d)^{d+1}}   \vol(\x_{\frac{\bullet+\beta}{2}})^2   \theta^\beta.\]
\end{compactenum}
For (1), consider the vertices of the unimodular simplex $S$, and consider $\overline{\alpha}$ to be their sum. It is clear that it satisfies Equation~\eqref{eq:idc}: We have  
\[\sum_{ \beta \in (\Delta\cap \mathbb{Z}^d)^{d+1}}   \vol (\x_{\frac{\overline{\alpha}+\beta}{2}})^2   \theta^\beta \ =\  \sum_{ \beta \in (S\cap \mathbb{Z}^d)^{d+1}}   \vol (\x_{\frac{\overline{\alpha}+\beta}{2}})^2   \theta^\beta.\]
It is easy to see that $\frac{\overline{\alpha}+\beta}{2}$ is a lattice point if and only if $\sum \beta_i =\overline{\alpha}$, that is, if $\beta_i$ ranges over the vertices of $S$. Hence, the right hand side equals
\[\vol(\x_S)^2 \det(\Theta_{|S})\ =\ \vol(\x_S)\]
where the equation follows from the Kustin-Miller normalization.

For (2), we consider the following system of equations arising from the linear system of parameters: For every $A$ in 
$\cone(\Delta)\cap (\mathbb{Z}^d\times \{d\})$, and every $i\in \{1,\dots,d+1\}$, we have 
the linear equation
\begin{equation}
0\ =\ \sum_{j \in \Delta\cap \mathbb{Z}^d} \theta_{i,j} \vol(\x_j \x_A).
\end{equation}
These are the linear equations determining $\vol$ because the semigroup algebra of $\Delta$ is Gorenstein. 

It remains to answer (3). Hence, we are left with verifying that 
\begin{equation}
0\ =\ \sum_{j \in\Delta \cap \mathbb{Z}^d} \theta_{i,j} \sum_{ \beta \in (\Delta \cap \mathbb{Z}^d)^{d+1}}   \vol (\x_{\frac{A+j+\beta}{2}})^2   \theta^\beta.
\end{equation}
For this, notice that the right side equals 
\begin{align*}
	&\sum_{j \in\Delta\cap \mathbb{Z}^d} \theta_{i,j} \sum_{\substack{ \beta \in (\Delta\cap \mathbb{Z}^d)^{d+1} \\  \beta_i= j } }   \vol (\x_{\frac{A+j+\beta}{2}})^2   \theta^\beta\\ \ =\ 
	&\sum_{j \in\Delta \cap \mathbb{Z}^d} \sum_{\substack{ \beta \in (\Delta\cap \mathbb{Z}^d)^{d+1} \\  \beta_i= j } }   \vol (\x_{\frac{A+j+\beta}{2}})^2 \theta_{i,j}^2 \prod_{\substack{1 \leq k \leq d+1,\\ k\neq i}} \theta_{k,\beta_k}.
\end{align*}
We see that the right side of this equation equals
\[\sum_{\substack{ B \in (\Delta\cap \mathbb{Z}^d)^{d} } } \prod_{\substack{1 \leq k \leq d+1,\\ k\neq i}} \theta_{\beta_k,k} \sum_{j \in\Delta \cap \mathbb{Z}^d}    \vol (\x_{\frac{A+2j+B}{2}})^2 \theta_{i,j}^2.\]
But 
\[\sum_{j \in\Delta \cap \mathbb{Z}^d}  \vol (\x_{\frac{A+2j+B}{2}})^2 \theta_{i,j}^2\ =\ \left(\sum_{j \in\Delta\cap \mathbb{Z}^d}    \vol (\x_{\frac{A+B}{2}} \x_j ) \theta_{i,j}\right)^2\ =\ 0\]
where the last equation is a squaring of the linear equations determining $\vol$.
\end{proof}

\subsection{The Parseval-Rayleigh identities imply the differential identity}\label{sub:diffid}
We now prove the differential equations.

\begin{proof}[\textbf{Proof of Lemma~\ref{lem:diffid} using Lemma~\ref{lem:parsevalcomplex}}]  
	We have \[\vol (\x_\sigma u^2)\ =\ \sum_{ \beta \in (\Delta\cap \mathbb{Z}^d)^{d+1}}  \vol (u\cdot\x_{\frac{\sigma+\beta}{2}})^2 \theta^\beta.\] Differentiating after $F$ gives the desired.
\end{proof}

\section{Beyond the integer decomposition property}\label{sec:beyond}

Before we discuss open questions, let us point out again that the integer decomposition property was used once, and once only: when concluding the anisotropy property from the differential relations/Parseval-Rayleigh identities. Let us return therefore to Theorem~\ref{thm:ani} as an example. As a reminder: it states that if $X$ is an IDP lattice ball or sphere of dimension $d$, and the characteristic of $\fld$ is $2$ or $0$, then the generic Artinian reduction $\AR^\ast(X, \partial X)$ of  $\fld^\ast[X,\partial X]$ has the anisotropy property for $\AR$ over $\fld(\theta_{i,j})$: for every nontrivial $u\in\AR^{k}(X,\partial X)$ of degree $k\le \nicefrac{d+1}{2}$, we have \[u^2\ \neq\ 0.\]

Now, if we want to remove the integer decomposition property, something interesting happens. Not only does anisotropy fail: it turns into its opposite. We have the following consequence of the Parseval-Rayleigh identities, which makes the difference most clear:

\begin{prp}[(An)Isotropy in lattice polytopes]\label{thm:aninonIDP}
If $X$ is a lattice ball or sphere of dimension $d=2k-1$, and the characteristic of $\fld$ is $2$, then the generic Artinian reduction $\AR^\ast(X, \partial X)$ of  $\fld^\ast[X,\partial X]$ has the following dichotomy: Consider $u\in\AR^{k}(X,\partial X)$ of degree $k$ that pairs with some $\x_{\nicefrac{F}{2}}$ whose square is generated in degree $1$, that is, there is an $F\in (X\cap \mathbb{Z}^d)^{d+1}$ such that $u\x_{\nicefrac{F}{2}}$ is not $0$. Then
\[u^2\ \neq\ 0.\]
Otherwise, that is, if $u$ pairs with no such element $\x_{\nicefrac{F}{2}}$ as above,
\[u^2\ =\ 0.\]
The latter in particular also applies under any specialization of the $\theta_{i,j}$: the element $u$ is isotropic.
\end{prp}

Let us summarize where the IDP comes in, and which theorems are general.

\begin{compactitem}[$\circ$]
\item The Kustin-Miller normalization, and its well-definedness, are independent of the integer decomposition property. Similarly, so are the Parseval-Rayleigh identities and the differential identities.
\item The crucial junction appears when we conclude anisotropy from the Parseval-Rayleigh identities. At this point, we needed to use the fact that every element $u$ in a semigroup algebra $\AR^\ast(X, \partial X)$ pairs with \emph{some} element of degree one.
\end{compactitem}

\section{Outlook and open questions}\label{sec:discussion}

The non-lattice cases of Stanley's conjecture remain. Even for IDP lattice polytopes or stronger yet, for lattice polytopes with a regular unimodular triangulation, we are left with a gap in the inequalities restricting $h^{\ast}$ if the interior of the cone is generated in higher degree. We conjecture that the unimodality of (the coefficients of) $h^{\ast}$ fails in general.

The intuition here is that lattice polytopes behave like triangulated disks, which can have non-unimodal $h$-vectors. The idea here could rely on constructing appropriate connected sums: as we saw above, Gorenstein polytopes have $h^\ast$-polynomials that peak at half of their socle degree (which is $d+1-s$, $s$ being the minimal dilation constant so that the polytope has an interior vertex). By considering the union of two polytopes with different socle degree, one could hope to turn a dromedary into a camel (though a mythical beast with more humps is not beyond our imagination, alas such a creature has to be high-dimensional). 

A word of caution, however, lies in an inequality for the $h^\ast$-polynomial arising from work of Eisenbud and Harris \cite{Stanleydomain}: we have that for any nonnegative $k$, and $s$ the degree of the $h^\ast$-polynomial, we have
\[h^\ast_0\ +\ \dots\ +\ h^\ast_k\ \le\ h^\ast_s\ +\ \dots\ +\ h^\ast_{s-k}.\]
This inequality is special to domains, and prevents us from introducing a hump below half the socle degree easily; it remains to understand the impact of this inequality in general. The most promising approach is then to look among polytopes whose interior is generated in high degree, and look for non-unimodality between half the degree of the $h^\ast$-polynomial and half the dimension of the polytope.
Another direction that we shall investigate in \cite{APP4} is the impact of restricted systems of parameters. Here, we proved the Lefschetz property for linear systems of parameters that correspond to orbifold Chow rings; this leads to unimodality results for the local $h^\ast$-vector.

\textbf{Acknowledgements.} 
We thank Giuliamaria Menara for proofreading, Benjamin Nill, Mykola Pochekai, Paco Santos and Alan Stapledon for helpful comments, and Eric Katz and Tim Römer for deep and insightful discussions.
K. A. and V. P. were supported by Horizon Europe ERC Grant number: 101045750 / Project acronym: HodgeGeoComb. We benefitted from experimentation with Macaulay2 \cite{M2}.

{\small
\bibliographystyle{myamsalpha}
\bibliography{ref}}

\appendix

\section{Characteristic $p$ (and the relation to Stanley-Reisner rings)}

Let us note that the similar identities apply in characteristic $p>0$. For instance, we have:
\begin{thm}\label{lem:parsevaldcomplexp}
	For a lattice $d$-ball or sphere $\Delta$, and $\alpha$ is a lattice point of $\cone(\Delta)\cap (\mathbb{Z}^d\times \{d+1\})$ we have in $\AR^\ast(\Delta,\partial \Delta)$ over characteristic $p$
	\begin{equation}\label{eq:id}
		\vol (\x_\alpha) \   =\    \sum_{ \substack{\beta \in \mathbb{Z}_{\ge 0}^{[d+1]\times (\Delta\cap \mathbb{Z}^d)}\\ \beta \cdot \mathbf{1}^{\Delta\cap \mathbb{Z}^d} =(p-1)\mathbf{1}^{[d+1]} }}   \vol (\x_{\frac{\alpha+\beta}{p}})^p   \frac{\theta^\beta}{\beta!}.
	\end{equation}
\end{thm}
The proof is analogous to the proof of the Parseval identities in characteristic two, see Section~\ref{sec:PC2}. Moreover, $\beta \cdot \mathbf{1}^{\Delta\cap \mathbb{Z}^d} =(p-1)\mathbf{1}^{[d+1]}$ means that as a matrix, the entries of each row of $\beta$ sum up to $p-1$.

Here, $\beta$ are matrices over the lattice points of $\Delta$ with integer entries and $d+1$ rows (that is, matrices just fitting the size of the linear system of parameters) and whose rows sum to $p-1$, and we have $\theta^\beta: =\prod \theta_{i,j}^{\beta_{i,j}}$ and $\beta! :=\prod \beta_{i,j}!$.

In particular, we have 

\begin{thm}[The Parseval-Rayleigh identity]\label{lem:parsevalcomplexp}
	For a lattice $d$-ball or sphere $\Delta$, in $\AR^\ast(\Delta,\partial \Delta)$ over characteristic $p$, and $\sigma$ an interior simplex and for $d+1-\#\sigma$ divisible by $p$, we have
	\begin{equation}
		\vol (\x_\sigma u^p)\ =\ \sum_{ \substack{\beta \in \mathbb{Z}_{\ge 0}^{[d+1]\times (\Delta\cap \mathbb{Z}^d)}\\ \beta \cdot \mathbf{1}^{\Delta\cap \mathbb{Z}^d} =\mathbf{p-1}^{[d+1]} }}  \vol(u\cdot\x_{\frac{\sigma+\beta}{p}})^p  \frac{\theta^\beta}{\beta!}.
	\end{equation}
	for all $u \in \AR^{\nicefrac{d+1-\#\sigma}{p}}(\Delta,\partial \Delta)$. 	
\end{thm}

In particular, we obtain generalizations of the Lefschetz properties for simplicial cycles (of \cite{APPHAL, AHKKS}) in characteristic $p$ to lattice complexes. In the special case of simplicial spheres, these identities were also obtained in \cite{KLS}. 
That said, in the case of simplicial complexes, the formula is comparatively straightforward, as we have explicit forms of the volume map; moreover, as we shall see in the next section, the Parseval-Rayleigh identities are equivalent to the differential identities that were already known in any positive characteristic.

In lattice polytopes, the Parseval-Rayleigh identity is the best thing we have to describe the volume map, and came before the differential identities.

But it is a proper generalization: Lattice complexes contain the case of simplicial complexes, and face rings (also known as Stanley-Reisner rings): By considering a simplicial complex as a lattice complex where each simplex is a unimodular lattice simplex, this specializes to previous $p$-anisotropy theorems in the aforementioned. We will explore these identities in an upcoming paper.

\section{Differential equations and the Euler Formula}

We have seen that the Parseval-Rayleigh identities imply the differential equations that we needed to prove our anisotropy theorem. Of course, as both differential equations and Parseval-Rayleigh identities are nontrivial nonhomogeneous relations that uniquely determine the fundamental class, they are equivalent. However, their connection is more direct: We have seen that the Parseval-Rayleigh identities imply the differential equations. We now also provide the other direction, even if it is unnecessary for our purposes. It is quite simple, yet beautiful.

Assume that $d\geq 2$, $m\geq 3$, $p$ is a prime number and $k$ is a field of characteristic $p$. Consider the polynomial ring 
\[
R=\fld^\ast[\theta_{(i,j)} \; :  \;   1 \leq i \leq d,  \; 1 \leq j \leq m].
\]
Assume $f,g\in R\setminus \{0\}$ such that, for all $i$ with $1\leq i \leq d$, 
the polynomials  $f, g$ are homogeneous with respect to the set of variables 
\[\theta_{(i,1)}, \theta_{(i,2)},\dots ,\theta_{(i,m)} \]
of degrees $\deg_i f$,  $\deg_i g$, respectively and
\[\deg_i f - \deg_i g= -1.\]
We denote by $\mathcal{A}\subseteq R$ the following set of monomials:
\[\mathcal{A}= \{z=\prod_{1\leq i\leq d, 1\leq j\leq m} \theta_{(i,j)}^{e_{i,j}}\colon 
e_{i,j}\geq 0 \,  \text{ and\, for\, all i, } \sum_{1\leq j\leq m} e_{i,j}= p-1 \}.\]

\begin{rem}
It is clear that each $z\in \mathcal{A}$ is a homogeneous element of $R$ of degree $d(p-1)$. Moreover, 
for each $i$ with $1\leq i\leq d$, $z$ is homogeneous with respect to the variables
\[\theta_{(i,1)}, \theta_{(i,2)},\dots ,\theta_{(i,d)}\]
of degree $p-1$.  Each $z\in \mathcal{A}$ defines the differential operator
\[\partial_{z}=\frac{\partial^{d(p-1)}}{\prod_{1\leq i\leq d, 1\leq j\leq m}(\partial \theta_{(i,j)}^{e_{i,j}})}.\]
\end{rem}

\begin{ex}
Assume that $d=2, m=6$, the characteristic of the field is equal to $3$
and  $z=\theta_{(1,1)} \theta_{(1,5)} \theta_{(2,1)}^2\in \mathcal{A}$.  Then,
\[\partial_{z}= \frac{\partial^{4}}{\partial \theta_{(1,1)} \partial \theta_{(1,5)} \partial \theta_{(2,1)} \partial \theta_{(2,1)} }.\]
\end{ex}

\begin{prp}\label{Prop!deg}
We have
\[\frac{f}{g}= (-1)^d \sum_{z\in \mathcal{A}} z \partial_{z}\left(\frac{f}{g}\right).\]
\end{prp}

\begin{proof}
For fixed $i$, with $1\leq i\leq m$, the rational function ${f}/{g}$ is homogeneous with respect to the variables
\[\theta_{(i,1)}, \theta_{(i,2)},\dots ,\theta_{(i,m)}\]
of degree equal to $-1$. Then, from the Euler formula for homogeneous rational functions (\cite[Theorem 1]{H}) 
\[[(p-1)!]^d\frac{f}{g}= \sum_{z\in \mathcal{A}} z \partial_{z}\left(\frac{f}{g}\right).\]
By Wilson's theorem $p$ divides $(p-1)! +1$.   The proposition follows.
\end{proof}

\begin{rem}
We refer the reader to the previous Section~\ref{sec:PC} for the setting of Parseval-Rayleigh Identities. Assume $D$ 
is a simplicial sphere of dimension $d-1$  or $P$ is an IDP lattice polytope of 
dimension  $d-1$ and $m$ is a monomial in the $\x_i$ of  degree $d$.   
We have that  $\vol(m)$ is an element of the field of fractions  $Q(R)$ of $R$. 
It is clear that $\vol(m)$ is, for fixed $m$, homogeneous with respect to the set of variables
\[\theta_{(i,1)}, \theta_{(i,2)},\dots ,\theta_{(i,m)}\]
of degree $-1$.  Hence,  Proposition~\ref{Prop!deg} implies that
\[
\vol(m)= (-1)^d \sum_{z\in \mathcal{A}} z \partial_{z}(\vol(m)).
\]
Hence the differential identities imply the Parseval-Rayleigh Identities.
\end{rem}

\section{The Kustin-Miller normalization, revisited}\label{sec:KMrevisit}

In this section we  give a  second proof that the volume  normalization of a lattice polytope $P$ induced by a 
boundary flag  is independent, up to sign, of the choice of the flag. As a bonus, we work in arbitrary characteristic, paying attention to signs. 

\vspace{1.5mm}

In the following  $P \subset \mathbb{R}^N$ denotes a $d$-dimensional lattice polytope with lattice 
point set  $\{1,\dots,m\}$. We denote by  $L(P)$  the set of lattice points of $\tau$. 
For a nonempty subset $Z$ of  $\mathbb{R}^N$  we denote by  $\operatorname{aff}(Z) $ the
smallest affine subset of  $\mathbb{R}^N$  containing~$Z$.

\begin{prp}\label{prp!Claim_1}
Assume   
\[ 
\tau = (\tau_0, \tau_1,  \dots ,  \tau_d = P) 
\] 
is a boundary flag of  $P$,  in the sense that for all  $1 \leq i \leq d$ we
have that  $\tau_{i-1}$  is a facet of the boundary of  $\tau_i$.
Suppose      $a_0 \in L(\tau_0)$    and for all  $1 \leq i \leq d$,
$ a_i \in L(\tau_i) \setminus L(\tau_{i-1})$.
Then,   for all  $0 \leq i \leq d$,   it holds that 
\[\operatorname{dim}   (\operatorname{aff} (a_0,\dots,a_i))  = i\]
and
\[ \operatorname{aff} (a_0,\dots,a_i)  =  \operatorname{aff} (\tau_i)  .\]
\end{prp}

\begin{proof}
It is well-known that 
$\operatorname{dim} (\operatorname{aff}(\tau_i)) = i$.
We use induction on $i$.  For $i = 0$,  we have  $ \tau_0 = \{a_0\}$
and the two claims are true.   Assume that  $1 \leq i \leq d-1$  and
the two claims are true for $i$.   Hence
$\operatorname{dim}   (\operatorname{aff} (a_0,\dots , a_i))  = i$
and
$( \operatorname{aff} (a_0,\dots ,a_i))  = (\operatorname{aff} (\tau_i))$.
Since   $\tau_i$  is a facet of  $\tau_{i+1}$, it holds that 
$ \tau_i  =   \tau_{i+1} \cap  H$ 
where  $H = \operatorname{aff} (\tau_i)$.   Since  
$a_i \in L(\tau_i) \setminus L(\tau_{i-1})$ 
we get that   $a_i \notin H$.   Consequently, 
$H$ is a proper subset of the $\operatorname{aff} ( H \cup  \{a_{i+1}\})$.
This implies that  
$\operatorname{dim}   (\operatorname{aff} (H \cup \{a_{i+1}\}))  = i+1$.
Since by the inductive hypothesis  
$ (\operatorname{aff} (a_0,\dots,a_i))  = (\operatorname{aff}(\tau_i))$ 
we get that 
$\operatorname{dim} (\operatorname{aff} (a_0,\dots,a_{i+1})) = i + 1$.
Taking into account   that    
$(\operatorname{aff}(a_0,\dots,a_{i+1}))  \subset  (\operatorname{aff}(\tau_{i+1}))$ 
and that 
$\operatorname{dim} (\operatorname{aff}(\tau_{i+1})) = i + 1$,
we get that 
$(\operatorname{aff}(a_0,\dots,a_{i+1})) =  (\operatorname{aff} (\tau_{i+1}))$,
which finishes the proof.
\end{proof}

We consider the $(d+1) \times m$   matrix  $M_{\theta} = [\theta_{i,j}]$, 
where   $\theta_i = \sum_{j=1}^m  \theta_{i,j} \x_j$ are the linear polynomials we use for the Artinian reduction.
We use the notations
\[ (a_0,\dots, a_d) =   [a_0,\dots, a_d] \x_{a_0}\x_{a_1} \dots \x_{a_d} \]  
and
\[  \mbf{R}_{a_0,\dots, a_{d-1}} =  \sum_{i=1}^{m} (a_0, \dots a_{d-1},i)   =
\sum_{i=1}^{m} [a_0,\dots a_{d-1},i] \x_{a_0}\x_{a_1}\dots \x_{a_{d-1}} \x_i,   \]
where $[a_0,  \dots , a_{d}] $ denotes the determinant of the 
$(d+1) \times (d+1)$ submatrix of $M_{\theta}$
specified by the  columns $a_0, \dots,  a_{d}$. 

\begin{prp} \label{prop!Claim_2rep}
\begin{itemize} 
\item[(i)] We have that 
\[
\mbf{R}_{a_0,  \dots, a_{d-1} } =  \operatorname{sign} (\sigma)   \mbf{R}_{b_0, \dots   ,b_{d-1} } 
\]
if $(b_0, \dots   ,b_{d-1} )$ is obtained from $(a_0, \dots   ,a_{d-1} )$  by a permutation
$\sigma$ and  $\operatorname{sign} (\sigma)$ denotes the 
sign of the permutation.   Moreover, 
\[
\mbf{R}_{a_0,  \dots, a_{d-1}} = 0 
\]
if there exists  $i\neq j$ with  $a_i=a_j$.

\item[(ii)]  Given  $a_0,\dots, a_{d-2}$ , we have that    
\[
\sum_{i=1}^{m} \mbf{R}_{a_0,a_1,\dots,a_{d-2},i} = 0.
\]
\end{itemize}
\end{prp}

\begin{proof}
Since  $[a_0,\dots, a_{d-1},i] $ is a determinant, (i) follows immediately.

We now prove (ii).  We have 
\[ \sum_{i=1}^{m} \mbf{R}_{a_0,a_1,\dots, a_{d-2}, i} =  \sum_{i=1}^{m}  \sum_{j=1}^{m} (a_0,a_1,\dots, a_{d-2}, i,j).\]
Since   for all $i$
\[(a_0,a_1,\dots, a_{d-2}, i,i)  = 0 \]   
and   when  $i \not= j$ we have 
\[(a_0,a_1,\dots, a_{d-2}, i,j)+ (a_0,a_1,\dots, a_{d-2}, j,i)=0 \]
the result follows.
\end{proof}

We keep assuming that  $P$ is a $d$-dimensional lattice polytope with lattice 
point set  $\{1,\dots, m\}$.   Assume 
$ T = (\tau_0, \tau_1,\dots, \tau_d = P)$  
is a boundary flag of $P$ in the above sense. 
We define the expression  $H_{T}$ as follows
\[H_{T}   =    \sum_{a_0,\dots,a_d} (a_0,\dots,a_d)\]
with the sum for all  $a_i \in L(\tau_i)$.

\begin{prp}\label{prp!Claim3}
We have that 
\[
H_{T}   =    \sum_ {a_0,\dots,a_d} (a_0,\dots,a_d)
\]
with the sum for $a_0\in L(\tau_0)$ and   $a_i\in L(\tau_i)\setminus L(\tau_{i-1})$  for $i>0$.
\end{prp}

\begin{proof}    
Assume  $a_i \in L(\tau_i)$, for $i=0\dots d$.
If  $a_1\in L(\tau_0)$  we get that  $a_1=a_0$. Hence,  $(a_0,\dots,a_d)=0$. Consequently,
\[H_{T}   =    \sum _{a_0,\dots,a_d} (a_0,\dots,a_d)\]
with the sum for  $a_0\in L(\tau_0)$,  $a_1\in L(\tau_1)\setminus L(\tau_0)$ and  $a_i\in L(\tau_i)$  for  $i\geq 2$.
Assume now that  $a_2\in L(\tau_1)$. 
If $a_2=a_0$ then $(a_0,\dots,a_d)=0$.
Similarly if $a_2=a_1$. Otherwise, both terms $(a_0, a_1, a_2,\dots, a_d)$ and $(a_0, a_2, a_1,\dots, a_d)$ 
appear in the sum defining  $H_{T}$ and they cancel each other.
Consequently,
\[H_{T}   =    \sum_{a_0,\dots,a_d} (a_0,\dots,a_d)\]
with the sum for  $a_0\in L(\tau_0)$,  $a_1\in L(\tau_1)\setminus L(\tau_0)$, 
$a_2\in L(\tau_2)\setminus L(\tau_1)$  and $a_i\in L(\tau_i)$ for $i\geq 3$.
Continuing on the same way  the result follows.
\end{proof}

\begin{prp}\label{prp!Claim4}
We have that 
\[
H_{T}   =    \sum_{a_0,\dots, a_{d-1}} \mbf{R}_{a_0,\dots, a_{d-1}} 
\]    
with the sum  for  $a_i \in L(\tau_i)$.
\end{prp}

\begin{proof}
It is clear from the definitions.       
\end{proof}

\begin{prp}\label{prp!Claim5}
We have that 
\[   H_{T}   =    \sum_{a_0,\dots, a_{d-1}} \mbf{R}_{a_0,\dots, a_{d-1}} \]
with the sum  for  $a_0\in L(\tau_0)$ and $a_i\in L(\tau_i)\setminus L(\tau_{i-1})$  for  $i>0$.
\end{prp}

\begin{proof}
It follows immediately from  Proposition~\ref{prp!Claim3}.
\end{proof}

Assume  $t < d $ and  we have two boundary flags  of $P$,
\[                T_1 = (\sigma_0, \sigma_1,  \dots ,  \sigma_d = P)\quad \text{and}\quad  T_2 = (\rho_0, \rho_1,  \dots ,  \rho_d = P),  \]
with the property that they only differ on the $t$-position, in the sense that      $\sigma_i =  \rho_i$  when  $i$ is different from $t$   and       $\sigma_t  \not=  \rho_t.$

We set  for $j > 0$ and  $j \notin \{t, t+1\}$
\[                  S_j = L(\sigma_j) \setminus  L(\sigma_{j-1}),\]
\[                 E_{sp} =   L(\sigma_{t+1}) \setminus (L(\sigma_t) \cup  L(\rho_t) ),\]
\[                  U_1 =  S_1 \times S_2 \times \dots \times S_{t-1},\ \ \text{and}\]
\[                 U_2 =  S_{t+2} \times S_{t+3} \times  \dots \times S_{d}.\]
To set $U$ depending on the value of $t$ and $d$:
\[ \text{ if   }  t=0  \, \,  \textit{and} \,  \,  d=1, \text{ then } \,  \, U= E_{sp}.\]
\[\text{ if   }t=0  \,  \,  \textit{and}  \,  \,  d\geq 2 ,\text{ then }  \,  \,   U= E_{sp}\times U_2.\]
\[\text{ if   }t=1 \,  \,  \textit{and} \,  \,    d< 3,\text{ then }  \,  \,   U= L(\sigma_0)\times E_{sp}.\]
\[\text{ if   } t=1 \,   \, \textit{and}  \,  \, d\geq 3,\text{ then }  \,   \,   U= L(\sigma_0)\times E_{sp}\times U_2.\]
\[\text{ if   } t\geq 2  \, \, \textit{and}  \,   \, d\geq t+2,\text{ then }  \,  \,  U = L(\sigma_0) \times U_1  \times E_{sp} \times U_2.\]
\[\text{ if   } t\geq 2  \, \,  and  \, \,  d< t+2,\text{ then }  \,   \,    U= L(\sigma_0)\times U_1\times E_{sp}.\]
where  $\times$ denotes the Cartesian product of sets.
For $z = (z_1,z_2, \dots, z_d)  \in U$,  we set 
\[
\mbf{R}_z =  \sum_{i=1}^m (z_1, z_2, \dots , z_d, i).
\]
Hence,  
\[\mbf{R}_z =  \sum_{i=1}^m [z_1,z_2, \dots, z_d, i]  \x_{z_1} \dots   \x_{z_d}   \x_i.\]

\begin{thm}\label{thm!theor101}
\begin{itemize}
\item[(i)]  We have the following equality of polynomials
\[                       H_{T_1} + H_{T_2}  =  (-1)^{t+d} \sum_{z \in U} \mbf{R}_z.\]

\item[(ii)]     Assume  $z = (z_1,z_2, \dots, z_d)  \in U$.  Then   
the set  $\{z_1, z_2, \dots, z_d\}$ is not contained in a facet of $P$.     
\end{itemize}       
\end{thm}

\begin{proof}
We prove (i) only for the case $t\geq 2$ and $d\geq t+2$. All the other cases can be proven in a similar way.

We assume first that $t \geq 2$ and $d> t+2$  and that we have two boundary flags,   

\[ T_1 = (\sigma_0, \dots, \sigma_{t-1}, \sigma_t, \sigma_{t+1}, \dots , \sigma_d = P) \quad \text{and} \quad  T_2 = (\sigma_0, \dots, \sigma_{t-1}, \rho_t, \sigma_{t+1}, \dots , \sigma_d = P),\]  
differing only on the $t$-th position.

We set  $L(\sigma_{0})=\{a\}$ and 
\[M_1 =  \{a\}   \times   (L(\sigma_1)\setminus L(\{a\}))  \times \dots \times (L(\sigma_{t-1})\setminus L(\sigma_{t-2}) ),\]  
\[M_2 = (L(\sigma_{t+2})\setminus L(\sigma_{t+1}) \dots  \times (L(\sigma_{d-1})\setminus L(\sigma_{d-2}). \]

We have 

\[U=  M_1 \times  (L(\sigma_{t+1})\setminus (L(\sigma_{t})\cup L(\rho_{t})) )  \times   M_2  \times (L(\sigma_d)\setminus (L(\sigma_{d-1}) ).\]

By Proposition~\ref{prp!Claim5},
\[ H_{T_1}= \sum_{(a_0,a_1,\dots,  a_{d-1}) \in W_1}  \mbf{R}_{a_0, a_1,  \dots,   a_{d-1}} \quad \text{and}\quad H_{T_2}= \sum_{(b_0,b_1,\dots  b_{d-1}) \in W_2 }   \mbf{R}_{b_0, b_1 \dots,   b_{d-1}},\]
where
\[W_1= M_1 \times  (L(\sigma_t)\setminus L(\sigma_{t-1}) )  \times   (L(\sigma_{t+1})\setminus L(\sigma_{t}))\times  M_2\] 
and
\[W_2=  M_1  \times  (L(\rho_t)\setminus L(\sigma_{t-1}) )  \times    (L(\sigma_{t+1})\setminus L(\rho_{t})) \times M_2.\]
We set 
\[V_1=  M_1 \times (L(\sigma_t)\setminus L(\sigma_{t-1}) )  \times   (L(\rho_t)\setminus L(\sigma_{t-1}) )  \times  M_2,\]
\[V_2=  M_1 \times (L(\sigma_t)\setminus L(\sigma_{t-1}) ) \times (L(\sigma_{t+1})\setminus (L(\sigma_{t})\cup L(\rho_{t}))\times M_2.\]

Since  $L(\sigma_{t+1})\setminus L(\sigma_{t})$   is the disjoint union of  
\[      L(\rho_t)\setminus L(\sigma_{t-1})  \,  \,   \textit{and} \,   \,  L(\sigma_{t+1})\setminus (L(\sigma_{t})\cup L(\rho_{t})) \]
it follows that  $W_1$  is the disjoint union of  $V_1$ and $V_2$.

We also set 
\[V_3=  M_1 \times (L(\rho_t)\setminus L(\sigma_{t-1}) )  \times   (L(\sigma_t)\setminus L(\sigma_{t-1}) )  \times  M_2,\]
\[V_4=  M_1 \times (L(\rho_t)\setminus L(\sigma_{t-1}) ) \times (L(\sigma_{t+1})\setminus (L(\sigma_{t})\cup L(\rho_{t}))  \times M_2.\]

Since    $L(\sigma_{t+1})\setminus L(\rho_{t})$   is the disjoint union of     
\[   L(\sigma_t)\setminus L(\sigma_{t-1})   \,   \,  \textit{and} \,  \,    L(\sigma_{t+1})\setminus (L(\sigma_{t})\cup L(\rho_{t})),  \]
it follows  that  $W_2$  is the disjoint union of $V_3$ and  $V_4$.

Therefore,  $H_{T_1}$ becomes a sum of two expressions, one for $V_1$ and one for $V_2$.   Similarly  
$H_{T_2}$   becomes a sum of two expressions, one for $V_3$ and one for $V_4$.   Moreover, 
the exchange of the $t$ and $t+1$ positions give a bijection between $V_1$ and $V_3$, and using 
the first part of Proposition~\ref{prop!Claim_2rep}, the corresponding
terms  in the sum    $H_{T_1} + H_{T_2}$   add to zero.  As a consequence,

\begin{align*}H_{T_1} + H_{T_2} \ &=\   \sum_{ u \in V_1} \mbf{R}_u  +   \sum_ {u \in V_2} \mbf{R}_u + \sum_ {  u \in V_3} \mbf{R}_u + \sum_ {  u \in V_4} \mbf{R}_u \\
  &=\   \sum_ { u \in V_2} \mbf{R}_u +  \sum_ { u \in V_4} \mbf{R}_u.
\end{align*}

Hence, 
\begin{align*}
 &H_{T_1} + H_{T_2} - (-1)^{(t+d)} \sum_{z \in U} \mbf{R}_z\\
=\ & \sum_{(a_1, \dots, a_{t-1}, a_t, a_{t+1}, a_{t+2}, \dots, a_{d-1})\in V_2}  \mbf{R}_{a, a_1, a_2, \dots, a_{d-1}}\\
+\ & \sum_{(b_1, \dots, b_{t-1},b_t, b_{t+1}, b_{t+2}, \dots, b_{d-1})\in V_4} \mbf{R}_{a, b_1,\dots, b_2,  \dots, b_{d-1}} - (-1)^{(t+d)} \sum_{z \in U} \mbf{R}_z
\end{align*}

This in turn equals
\begin{align*}
\ &  (-1)^{d-1-t} \sum_{(a_1, \dots, a_{t-1}, a_t, a_{t+1}, a_{t+2}, \dots, a_{d-1})\in V_2}  
\mbf{R}_{a, a_1, \dots, a_{t-1}, a_{t+1}, a_{t+2}, \dots , a_{d-1}, a_t}\\
+\ &	(-1)^{d-1-t} \sum_{(b_1, \dots, b_{t-1},b_t, b_{t+1}, b_{t+2}, \dots, b_{d-1})\in V_4} 
\mbf{R}_{a, b_1, \dots, b_{t-1}, b_{t+1}, b_{t+2}, \dots ,  b_{d-1}, b_{t}}\\
-\ & (-1)^{(t+d)} \sum_{z \in U} \mbf{R}_z\\
=\ & (-1) (-1)^{t+d} \sum_{(z_1, \dots ,z_{d-1}\in U} \sum_{i=1}^{m}  \mbf{R}_{z_1, z_2, \dots, z_{d-1}, i}\\
=\ &0
\end{align*}
where for the final equality  we used  that
\[     (-1)^{d-1-t} =  (-1) (-1)^{d}  (-1)^{-t} = (-1) (-1)^{d}  (-1)^{t} = (-1) (-1)^{t+d}\]
and  Proposition ~\ref{prop!Claim_2rep}.

We now assume that $t \geq 2$ and  $d= t+2$  and that we have two boundary flags   

\[T_1 = (\sigma_0, \dots, \sigma_{t-1}, \sigma_t, \sigma_{t+1}, \sigma_{t+2} = P),\]  
\[T_2 = (\sigma_0, \dots, \sigma_{t-1}, \rho_t, \sigma_{t+1}, \sigma_{t+2} = P).\]  

We set  $L(\sigma_{0})=\{a\}$ and

\[M_1 =  \{a\}   \times   (L(\sigma_1)\setminus \{a\})  \times \dots \times (L(\sigma_{t-1})\setminus L(\sigma_{t-2}) ). \] 

We observe that in this case $M_2$ defined above does not exist.

We have 

\[U=  M_1 \times  (L(\sigma_{t+1})\setminus (L(\sigma_{t})\cup L(\rho_{t})) )  \times   (L(\sigma_{t+2})\setminus L(\sigma_{t+1})).\]

By  Proposition~\ref{prp!Claim5},  

\[H_{T_1}= \sum_{(a_0,a_1,\dots,  a_{t+1}) \in W_1}  \mbf{R}_{a_0, a_1,  \dots,   a_{t+1}}, \text{and}\]

\[H_{T_2}= \sum_{(b_0,b_1,\dots  b_{t+1})\in W_2 }   \mbf{R}_{b_0, b_1 \dots,   b_{t+1}}\]

where,
\[W_1= M_1 \times  (L(\sigma_t)\setminus L(\sigma_{t-1}) )  \times   (L(\sigma_{t+1})\setminus L(\sigma_{t}))\]
and
\[W_2=  M_1  \times  (L(\rho_t)\setminus L(\sigma_{t-1}) )  \times    (L(\sigma_{t+1})\setminus L(\rho_{t})).\] 
We set 
\[V_1=  M_1 \times (L(\sigma_t)\setminus L(\sigma_{t-1}) )  \times   (L(\rho_t)\setminus L(\sigma_{t-1}) ),\]  
\[V_2=  M_1 \times (L(\sigma_t)\setminus L(\sigma_{t-1}) ) \times (L(\sigma_{t+1})\setminus (L(\sigma_{t})\cup L(\rho_{t})) ). \]
Since  $L(\sigma_{t+1})\setminus L(\sigma_{t})$   is the disjoint union of  
\[L(\rho_t)\setminus L(\sigma_{t-1}) \, \,    \textit{and} \,  \,  L(\sigma_{t+1})\setminus (L(\sigma_{t})\cup L(\rho_{t}))\] 
it follows that  $W_1$  is the disjoint union of  $V_1$ and $V_2$.
We also set 
\[V_3=  M_1 \times (L(\rho_t)\setminus L(\sigma_{t-1}) )  \times   (L(\sigma_t)\setminus L(\sigma_{t-1}) ), \] 
\[V_4=  M_1 \times (L(\rho_t)\setminus L(\sigma_{t-1}) ) \times (L(\sigma_{t+1})\setminus (L(\sigma_{t})\cup L(\rho_{t})) ).\]

Since    $L(\sigma_{t+1})\setminus L(\rho_{t})$   is the disjoint union of     
\[ L(\sigma_t)\setminus L(\sigma_{t-1})   \, \,   and \,   \,   L(\sigma_{t+1})\setminus (L(\sigma_{t})\cup L(\rho_{t})),\]  
it follows  that  $W_2$  is the disjoint union of  $V_3$ and $V_4$.

Therefore,  $H_{T_1}$ becomes a sum of two expressions, one for $V_1$ and one for $V_2$.   Similarly  
$H_{T_2}$   becomes a sum of two expressions, one for $V_3$ and one for $V_4$.   Moreover, 
the exchange of the $t$ and $t+1$ positions give a bijection between $V_1$ and $V_3$, and using the first part of Proposition~\ref{prop!Claim_2rep}, the corresponding
terms  in the sum    $H_{T_1} + H_{T_2}$   add to zero.  As a consequence,

\begin{align*}
H_{T_1} + H_{T_2} \ &=\  \sum_{ u \in V_1} \mbf{R}_u  +   \sum_ {u \in V_2} \mbf{R}_u + \sum_ {  u \in V_3} \mbf{R}_u + \sum_ {  u \in V_4} \mbf{R}_u \\
&=\ \sum_ { u \in V_2} \mbf{R}_u +  \sum_ { u \in V_4} \mbf{R}_u .
\end{align*}

Hence, 
\begin{align*}
	&\ H_{T_1} + H_{T_2} - (-1)^{(2t+2)} \sum_{z \in U} \mbf{R}_z\\
	=&\ \sum_{(a_1, \dots, a_{t-1}, a_t, a_{t+1})\in V_2}  \mbf{R}_{a, a_1, a_2, \dots, a_{t}, a_{t+1}}\\
	+&\ \sum_{(b_1, \dots, b_{t-1},b_t, b_{t+1})\in V_4} \mbf{R}_{a, b_1,\dots, b_2,  \dots, b_{t}, b_{t+1}}  - (-1)^{(2t+2)} \sum_{z \in U} \mbf{R}_z\\
	=&\ (-1) \sum_{(a_1, \dots, a_{t-1}, a_t, a_{t+1})\in V_2}  
	\mbf{R}_{a, a_1, \dots, a_{t-1}, a_{t+1}, a_t}\\
	+&\   (-1) \sum_{(b_1, \dots, b_{t-1},b_t, b_{t+1}, b_{t+2}, \dots, b_{d-1})\in V_4} 
	\mbf{R}_{a, b_1, \dots, b_{t-1}, b_{t+1}, b_{t}}\\
	-&\ (-1)^{(2t+2)} \sum_{z \in U} \mbf{R}_z\\
	=&\ (-1) \sum_{(z_1, \dots , z_{t+1})\in U} \sum_{i=1}^{m}  \mbf{R}_{z_1, z_2, \dots, z_{t+1}, i}\\
	=&\ 0
\end{align*}
where for the final equality  we used  the second part of  Proposition~\ref{prop!Claim_2rep}.

We now prove (ii).

We recall that every facet of a $d$-dimensional lattice polytope $P$ has dimension $d-1$.
Assume that $\{z_1, z_2, \dots, z_d\}$n is contained in a facet of $P$. Then  $\{z_1, z_2, \dots, z_d\}$ should have dimension less or equal to $d-1$.
By Proposition ~\ref{prp!Claim_1}, we get a contradiction since the affine span of  $\{z_1, z_2,\dots , z_d\}$ has dimension $d$.   
\end{proof}
Denote by  $\mathcal{A}(P, \partial P)$ the Artinian reduction of  $k[P, \partial P]$ parametrized by $\theta_{i,j}$. That is, 
\[\mathcal{A} (P, \partial P) = \mathcal{I}_{\partial P} / (\mathcal{I}_{\partial P} \mathcal{I}_{lins} + \mathcal{I}_P),\]
where we denote  by $\mathcal{I}_P$ and  $\mathcal{I}_{\partial P}$ the ideals of $P$ and of the boundary of $P$ respectively 
and by $\mathcal{I}_{lins}= (\theta_1, \dots , \theta_{d+1})$ the ideal of the linear relations, where the $\theta_{i,j}$ parametrize the coefficients as usual.

\begin{prp}
We have
\[  H_{T_1} + H_{T_2}=0 \]

in  $\mathcal{A}(P, \partial P)$.
\end{prp}

\begin{proof}

Recall that by the first part of  Theorem~\ref{thm!theor101},

\[H_{T_1} + H_{T_2}  =  (-1)^{t+d} \sum_{z \in U} \mbf{R}_z\]

where 

\[\mbf{R}_z =  \sum_{i=1}^m [z_1,z_2, \dots, z_d, i]  \x_{z_1} \x_{z_2} \dots  \x_{z_d}  \x_i\]

for $z = (z_1,z_2, \dots, z_d)  \in U$.

It is enough to prove that $\sum_{z \in U} \mbf{R}_z$ belongs to the ideal  $\mathcal{I}_{\partial P} \mathcal{I}_{lins} + \mathcal{I}_P$.

By the second part of Theorem~\ref{thm!theor101},   
\[ \x_{z_1}  \x_{z_2}  \dots   \x_{z_d}  \]
is an element of the ideal   $\mathcal{I}_{\partial P}$  of the boundary of the $P$. By Lemma~\ref{lem:balancing}.   The result follows.
\end{proof}

\section{Locality in lattice sheaves and almost pullbacks}\label{sec:locality}

The next three sections provide a second proof of the Parseval-Rayleigh identities, and provide some helpful results along the way that are of independent interest.
The locality Lemma~\ref{lem:balancing} and Lemma~\ref{thm:transfer} imply that for a facet $P$ of a lattice 
sphere or ball $X$, the volume map on $(P,\partial P)$ coincides with the volume map on $(X,\partial X),$ restricted to $P$.

It is natural to ask more questions. For instance, consider the following case:

Consider for instance a lattice polytope $P$ of dimension $d$, and a lattice polytope $Q$ of 
the same dimension inside it. 
Consider a monomial $m$ of 
$\fld^{d+1}[Q, \partial Q]$. What can be said of the relation of $\vol_Q(m)$ 
and $\vol_P(m)$ in $\AR^\ast(Q,\partial Q)$ resp.\ $\AR^\ast(P,\partial P)$? Unlike in the previous case, they do not coincide.

\begin{ex}  Consider the $1$-dimensional lattice polytopes    $Q=[1,2]  \subset P=[1,3]$ and the monomial
	$m = \x_1 \x_2 \in \fld^{d+1}[Q, \partial Q]$ and assume that 
	the field $\fld$ has characteristic $2$.
	An easy computation gives that 
	\[  
	\vol_Q (m) = \frac{1}{[1,2]}, \quad  \quad  \vol_P (m) = \frac{[2,3]}{[1,3]^2 + [1,2][2,3]},  
	\]
	where  $[a,b] = \det  \Theta_{|(a,b)}$  and $\Theta_{|(a,b)}$ denotes the 
	submatrix of $\Theta = [\theta_{i,j}]$ obtained by keeping the columns indexed by $a$ and $b$.
\end{ex}

We will ask such questions here, and prove two lemmata that we feel can be of independent interest. Let us for this purpose introduce another parameter $t$. Given 
a set of lattice points $V$ in a polytope $P$, we wish to study the following variation of the generic 
linear system of parameters: Instead of using the linear system of parameters $\theta_{i,j}$, we use a modified system of parameters
\[\theta^V_{i,j}[t]=\left\{\begin{array}{cl} t\theta_{i,j}& \text{if } j\in V \\
	\theta_{i,j}& \text{otherwise}.
\end{array}
\right.\]
We denote the corresponding Artinian reduction by $\AR^\ast(P,\partial P)[\theta^V_{i,j}[t]]$

\subsection{Finer lattices}

Consider the following situation: $P$ is a $d$-dimensional lattice polytope with lattice $\mathbb{Z}^d$. And $\Lambda$ is some finer lattice: it contains $\mathbb{Z}^d$ as a strict subset. Of course, we could consider 
$\AR^\ast(P,\partial P)[\mathbb{Z}^d]$, that is, $P$ and the semigroup algebra with respect to the lattice $\mathbb{Z}^d$. But we could equally consider $\AR^\ast(P,\partial P)[\Lambda]$. 

A prototype of such a situation is to consider $P$, and a positive dilate $nP$. 

In either case, how do they relate? The answer is actually easy:

\begin{lem}\label{lem:refined}
	If $V$ consists of those lattice points of $\Lambda\setminus \mathbb{Z}^d$ that lie in $P$, we  have 
	\[\AR^{d+1}(P,\partial P)[\Lambda][\theta^V_{i,j}[0]]\ =\ \AR^{d+1}(P,\partial P)[\mathbb{Z}^d].\]
	Moreover, if $m$ is a monomial of degree $d+1$ in the $\fld^\ast[P,\partial P]$, then, marking down the obvious dependencies, we have 
	\[\vol_{P,\Lambda, [\theta^V_{i,j}[t]]}(m)- \vol_{P,\mathbb{Z}^d, [\theta_{i,j}]}(m)\]
	is a rational function that vanishes at $t=0$.
\end{lem}	

Both of these facts are obvious. We come to a slightly more intricate case.

\subsection{Bigger polytopes}

We now go back to the original situation: 

\begin{lem}\label{lem:bigger}
	Consider a lattice polytope $P$ of dimension $d$, and a lattice polytope $Q$ of the same dimension inside it. Consider a monomial $m$ of 
	$\fld^{d+1}[Q, \partial Q]$, and assume additionally that $Q$ is obtained from $P$ by cutting the latter with a halfspace delimited by hyperplane $H$. Let $V$ denote the lattice points of $P$ not in $Q$. Then 
	\[\vol_{P,\Lambda, [\theta^V_{i,j}[t]]}(m)- \vol_{Q,\mathbb{Z}^d, [\theta_{i,j}]}(m)\]
	is a rational function that vanishes at $t=0$.
\end{lem}	

\begin{rmk}
	In fact, the additional assumption is not necessary, but that requires some further elaboration that we will only be able to discuss in the next section. We content ourselves with this version, which is enough for our purposes.
\end{rmk}

\begin{proof}
	$P$ and $Q$ have a vertex in common that is not in $H$. In particular, there is a full flag of faces of $P$ that restricts to a full flag of faces of $Q$. We consider the Kustin-Miller normalization of 
	\[\AR^\ast(P,\partial P)[\mathbb{Z}^d][\theta^V_{i,j}[t]]\]
	and
	\[\AR^\ast(Q,\partial Q)[\mathbb{Z}^d][\theta^V_{i,j}]\]
	with respect to this flag. 
	
	Notice that in the former, the Kustin-Miller normalization and linear relations do not uniquely determine $\AR^{d+1}(P,\partial P)$ if $t=0$; in general $\theta^V_{i,j}[0]$ is not a linear system of parameters for the semigroup algebra $\fld^{d+1}[P, \partial P]$, so 
	$\vol_{P,\Lambda, [\theta^V_{i,j}[t]]}(m)$ is not well-defined. However, if we restrict to the image of 
	$\fld^{d+1}[Q, \partial Q]$ in 
	\[\bigslant{\fld^{d+1}[P, \partial P]}{\theta^V_{i,j}[t]},\]
	it is of dimension $1$ independent of whether $t\neq 0$ or $t=0$. Hence, $\vol_{P,\mathbb{Z}^d, [\theta^V_{i,j}[t]]}(m)$, defined using the Kustin-Miller normalization on this subspace, is well-defined and has no pole at $t=0$. In particular, at this point, it coincides with $\vol_{Q,\mathbb{Z}^d, [\theta_{i,j}]}(m)$.
	
	To see this in terms of linear algebra, let us examine briefly what the linear equations determine $\vol_{P,\mathbb{Z}^d, [\theta^V_{i,j}[t]]}(m)$. If $\mbf{m}_P$ is the vector of monomials in $\fld^{d+1}[P, \partial P]$, with the elements of $\fld^{d+1}[Q, \partial Q]$ coming last. The linear relations between them are the relations of the form 
	\[\sum_{j} \theta^V_{i,j}[t] x_j \x_I\]
	for $\x_I$ monomials in $\fld^{d}[P, \partial P]$, with the relations corresponding to $\x_{I}$ in $\fld^{d}[Q, \partial Q]$ coming last purely by convention. Then, finally, we have the affine equation coming from the Kustin-Miller normalization.
	
	In other words, by passing to a basis, we can write this in matrix form as
	\[\left(\begin{array}{cc} 
		A & C \\
		B & \mr{M}^Q (t)
	\end{array}\right) \mbf{m}_P= e_P
	\]
	where $e_P=(0,\cdots,0,1)$ is a vector of appropriate length, and $\mr{M}^Q (t)$ is a matrix with entries depending on $t$ such that \[\mr{M}^Q (0) \mbf{m}_Q = e_Q \]
	is the analogous affine system of equations that uniquely defines $\vol_{Q,\mathbb{Z}^d, [\theta_{i,j}]}(m)$. 
	
	By performing column operations on those elements of  $\fld^{d+1}[P, \partial P]$ not in $\fld^{d+1}[Q, \partial Q]$, we obtain an equivalent equation
	\[\left(\begin{array}{cc} 
		A' & C \\
		0 & \mr{M}^Q
	\end{array}\right) \mbf{m}_P(t)= e_P.
	\]
	where $\mbf{m}_P(t)$ has entries that are rational functions in $t$ such that $\mbf{m}_P(0)$, restricted to the entries in $\fld^{d+1}[Q, \partial Q]$, coincides with $\mbf{m}_Q$ since $\mr{M}^Q(0)$ is invertible. Hence we obtain the desired.
\end{proof}

\section{Parseval-Rayleigh identities}\label{sec:parseval}

While the differential equation of Lemma~\ref{lem:diffid} is superficially similar to identities proven in the case of simplicial cycles in \cite{APP, APPHAL, PP}, where they follow immediately from the known formulas for the volume map in toric varieties, the case of lattice polytopes is much harder: we understand the volume map only indirectly, using a nonhomogeneous equation that takes the form of an identity of the Parseval-Rayleigh type. We consider lattice polytopes of dimension $d$ in $\mathbb{Z}^d$.

Assume $v = (v_1,  \dots , v_{d+1} )   \in  (P\cap\mathbb{Z}^d)^{d+1}$. We set as usual
\[\x_v\ =\ \prod_{1 \leq i \leq d+1}  \x_{v_i}.\]
Note that since we are working in the semigroup algebra, this only depends on \[|v|=v_1+v_2+\dots+v_{d+1},\] the sum over the entries of $v$ within the semigroup $\cone(P)\cap(\mathbb{Z}^d\times \mathbb{Z})$.

\begin{lem}\label{lem:parsevald}
	For a lattice $d$-polytope $P$, and $\alpha$ a lattice point $\cone^\circ(P)\cap (\mathbb{Z}^d\times \{d+1\})$, we have in $\AR^\ast(P,\partial P)$ over characteristic $2$
	\begin{equation}\label{eq:id}
		\vol (\x_\alpha) \   =\    \sum_{ \beta \in (P\cap \mathbb{Z}^d)^{d+1}\times \{1\}}   \vol (\x_{\frac{\alpha+\beta}{2}})^2   \theta^\beta.
	\end{equation}
\end{lem}

Here, we follow the convention $\theta^\beta \coloneqq\prod_{1 \leq i \leq d+1} \theta_{i,\beta_i}$. Moreover, $\vol (\x_{\frac{\alpha+\beta}{2}})$ is defined to be $\vol (\x_\gamma)$ if there is an  $\x_\gamma \in \fld^\ast[P]$ such that $\x_\alpha \x_\beta   = \x_\gamma^2$, and 0 otherwise.

This specializes to the following identity for $\alpha=\sigma + 2 \alpha'$, which explains the naming of this identity:

\begin{lem}\label{lem:parsevalpre}
	For a lattice $d$-polytope $P$, in $\AR^\ast(P,\partial P)$ over characteristic $2$, and $\sigma$ a family of lattice points of $P$ and for $d+1+\#\sigma$ even, we have
	\begin{equation}
		\vol (\x_\sigma \x_{\alpha'}^2)\ =\ \sum_{ \beta \in (P\cap \mathbb{Z}^d)^{d+1}\times \{1\}}  \vol (\x_{\alpha'}\cdot\x_{\frac{\sigma+\beta}{2}})^2 \theta^\beta
	\end{equation}
	where $\sigma+\beta$ denotes the concatenation of the families $\sigma$ and $\beta$. 
\end{lem}

From here we conclude identities deserving their name:

\begin{lem}[The Parseval-Rayleigh identity]\label{lem:parseval}
	For a lattice $d$-polytope $P$, in $\AR^\ast(P,\partial P)$ over characteristic $2$, and $\sigma$ a family of lattice points and for $d+1+\#\sigma$ even, we have
	\begin{equation}
		\vol (\x_\sigma u^2)\ =\ \sum_{ \beta \in (P\cap \mathbb{Z}^d)^{d+1}\times \{1\}}  \vol (u\cdot\x_{\frac{\sigma+\beta}{2}})^2 \theta^\beta
	\end{equation} 
	for all $u$ in $\fld^\ast[P,\partial P]$, and for $u\in \fld^\ast[P]$ if $\sigma$ is an interior simplex.
\end{lem}

\subsection{Beginning the proof: The case of the simplex} To simplify the proof of the Lemma~\ref{lem:parsevald}, we first prove a variant:

\begin{lem}\label{lem:parsevalSimplex}
	Consider a lattice $d$-simplex $S$ that is a dilation of a unimodular simplex, and let $\alpha$ denote the sum $\sum v$, where $v$ ranges over the vertices of $S\times\{1\}$. Then
	\begin{equation}\label{eq:iddescent}
		\vol (\x_\alpha) \   =\    \sum_{ \beta \in (S\cap \mathbb{Z}^d)^{d+1}\times \{1\}}   \vol (\x_{\frac{\alpha+\beta}{2}})^2   \theta^\beta.
	\end{equation}
\end{lem}

We move the proof to the next section, and first conclude the proof of the Parseval-Rayleigh identities for general polytopes.

\subsection{The general case}

We now obtain the proof of the general Parseval-Rayleigh identities. For this, we only need to use Theorems~\ref{lem:bigger} and Lemma~\ref{lem:refined}, as well as Lemma~\ref{lem:parsevalSimplex} of course.

Let us start with a pair of a rational polytope $Q$, and an interior point $\alpha$. Consider a polytope $P$ containing $Q$.

We say that the triple $(P,Q;\alpha)$ is \Defn{tight} if:
\begin{compactitem}[$\circ$]
	\item $Q=P \cap H$, where $H$ is some halfspace, and both $P$ and $Q$ are of the same dimension, and
	\item For every point $v$ in $P$, the point $\nicefrac{v+\alpha}{2}$ lies in the interior of $Q$.
\end{compactitem}

The two crucial lemmata now are the following:

\begin{lem}\label{lem:steps}
	Given any rational polytope $Q$ and a rational interior point $\alpha$, there exists a finite sequence $(P_i)$ $i=0,\cdots, n$ of rational polytopes such that	
	\begin{compactitem}[$\circ$]
		\item $P_n$ is a dilation of a unimodular simplex, and $\alpha$ is its barycenter.
		\item $P_0=Q$.
		\item $(P_{i+1},P_i;\alpha)$ is tight.
	\end{compactitem}
\end{lem}

\begin{proof}
	This is clear by gradually moving out the hyperplanes defining $Q$ until we are left with a simplex. Moving them out further ensures we can make $\alpha$ the barycenter. Doing this in discrete, small enough steps gives the finite sequence.
\end{proof}

We now return to lattice polytopes:

\begin{lem}\label{lem:hereditary}
	Consider $Q$ a lattice polytope of dimension $d$, and $\alpha$ a point in $\cone^\circ(Q)\cap (\mathbb{Z}^{d}\times \{d+1\})$, and $P$ a lattice polytope so that the triple $(P,Q;\nicefrac{\alpha}{d+1})$ is tight. Assume that the Parseval-Rayleigh identity of Lemma~\ref{lem:parsevald} holds for $P$ and $\alpha$. Then it holds for $Q$ and $\alpha$.
\end{lem}	

\begin{proof}
	This is an immediate consequence of Lemma~\ref{lem:bigger}.
\end{proof}	

Let us similarly note a curious consequence of Lemma~\ref{lem:refined}.

\begin{lem}\label{lem:hereditary2}
	Consider $Q$ a lattice polytope of dimension $d$ with respect to a lattice $\Lambda$, let $\Lambda' \supset \Lambda$ denote some finer lattice, and $\alpha$ a point in $\cone^\circ(Q)\cap (\Lambda'\times \{d+1\})$. Assume that the Parseval-Rayleigh identity of Lemma~\ref{lem:parsevald} holds for $Q$ and $\alpha$ with respect to $\Lambda'$.
	Then
	\begin{compactenum}[(1)]
		\item Consider the set $V$ of lattice points of $Q \cap \Lambda'$ not in $\Lambda$. Then, if $\alpha \notin \cone^\circ(Q)\cap (\Lambda\times \{d+1\})$, we have
		\[\vol_{Q,\Lambda', [\theta^V_{i,j}[0]]}(\x_\alpha)\ =\ 0.\]	
		\item The Parseval-Rayleigh identities hold for $Q$ and $\alpha$ with respect to $\Lambda$.
	\end{compactenum}
\end{lem}

\begin{rem}  \label{lem:dslkaadf}
	It is in fact not hard to see that the first conclusion of Lemma~\ref{lem:hereditary2} can be proven directly, and holds independently of the characteristic of the underlying field. Consider the semigroup algebra over a lattice polytope $Q$ with respect to a field of arbitrary characteristic. 
	
	Notice that $\theta^V_{i,j}[0]$ is a linear system of parameters for $\fld^\ast[Q]$ over $\Lambda'$, and that therefore the volume at a point $\alpha \notin \cone^\circ(Q)\cap (\Lambda\times \{d+1\})$ is well-defined. However, it is also a linear system over $\Lambda$. In particular, for $\alpha \notin \cone^\circ(Q)\cap (\Lambda\times \{d+1\})$, the normalization is independent of the normalization of volume map. Hence 
	\[\vol_{Q,\Lambda', [\theta^V_{i,j}[0]]}(\x_\alpha)\ =\ 0\]	 
	for all $\alpha \notin \cone^\circ(Q)\cap (\Lambda\times \{d+1\})$.
\end{rem}

\begin{ex} We give an example to demonstrate Lemma~\ref{lem:hereditary2}. Assume $\fld$ is a field of characteristic $2$,
	$\Lambda'= \mathbb{Z}^2$,  $\Lambda = 2 \Lambda'$ and $P \subset \mathbb{R}^2$ is the convex hull of 
	the set of points $\{(0,0),(0,2),(2,0)\}$. The lattice point set   of $P$ with respect to the lattice 
	$\Lambda'$ is the set
	\[
	\{  q_1 = (0,0),  \;  q_2 = (0,1),  \;  q_3 = (0,2), \;
	q_4 = (1,0),  \;  q_5 = (1,1),  \;  q_6 = (2,0) \},
	\]
	while the lattice point set   of $P$ with respect to the lattice 
	$\Lambda$ is the set  $\{  q_1, q_3, q_6 \}$.  For $1 \leq i \leq 3$, we denote by
	$\theta_i$ a general linear combination of the variables $x_1,x_3,x_6$.  Then  $\Theta\x = (\theta_1, \dots , \theta_3)$
	is a linear system of parameters for both $\fld^\ast_{\Lambda}[P]$ and $\fld^\ast_{\Lambda'}[P]$ and the 
	following holds:  Assume $m= \prod_{i=1}^6 x_i^{a_i}$ is  a monomial in $\fld^\ast[x_1, \dots , x_6]$ of degree $3$
	such that $m \in I_{\partial P}$.  If  $\; \sum_{i=1}^6 a_iq_i \notin \Lambda$,  then the class of $m$ in
	$\AR^3_{\Lambda', \Theta}(P,\partial P)$ is zero. Finally, we mention 
	that by Remark~\ref{lem:dslkaadf}  the same results hold
	in any characteristic.
\end{ex}

\begin{proof}[\textbf{Proof of Lemma~\ref{lem:parsevald}}]
	Consider a given lattice polytope $Q$ and $\alpha$ a point in $\cone^\circ(Q)\cap (\mathbb{Z}^{d}\times \{d+1\})$. We may assume that $\alpha=(0,\cdots,0,d+1)$.
	
	We want to prove Lemma~\ref{lem:parsevald} for $Q$ with respect to $\alpha$. Consider a sequence $P_i$ as given by Lemma~\ref{lem:steps}. Since the polytopes involved are rational, we can find a sufficiently large dilation $NP_i$ such that all involved polytopes are lattice polytopes. 
	
	We conclude from Lemma~\ref{lem:hereditary} and Lemma~\ref{lem:parsevalSimplex} that the Parseval-Rayleigh identities hold for $NQ$ with respect to $\alpha$. It follows that the Parseval-Rayleigh identities hold for $Q$ by Lemma \ref{lem:hereditary2}. 	
\end{proof}

\begin{rem}
	Note that this iterative reduction we used here immediately gives a strengthening of Lemma~\ref{lem:bigger}:	
\end{rem}

\begin{prp}\label{prp:bigger}
	Consider a lattice polytope $P$ of dimension $d$, and a lattice polytope $Q$ of the same dimension inside it. Consider a monomial $m$ of 
	$\fld^{d+1}[Q, \partial Q]$. Let $V$ denote the lattice points of $P$ not in $Q$. Then 
	\[\vol_{P,\Lambda, [\theta^V_{i,j}[t]]}(m)- \vol_{Q,\mathbb{Z}^d, [\theta_{i,j}]}(m)\]
	is a rational function that vanishes at $t=0$.
\end{prp}

\section{Parseval-Rayleigh identities on the simplex and the Liouville property}

In order to prove the Parseval-Rayleigh identity for simplices, we recall a basic identity for the volume:

\begin{lem}\label{lem:basic}
	Consider any two elements $I,J \in (P \cap \mathbb{Z}^d)^{d}$, where at least one point of $I$ lies in the interior of $P$, and an index $\mu \in [d+1]$. Then we have
	\[\mbf{R}[J,\mu,|I|]\ \coloneqq\ \sum_{p \in P \cap \mathbb{Z}^d} \theta^{J\ast_\mu p} \vol(\x_{I}\x_p)\ =\ 0.\]
\end{lem}

Here, $J\ast_\mu p$ is the vector family $J$ to which $p$ is inserted at place $\mu$. This is an immediate consequence of Lemma~\ref{lem:balancing}. To understand the principle of the proof, we first observe a curious property of the Liouville type:

\subsection{A Liouville type theorem}

Assume, somewhat uncleanly, that $P=\mathbb{R}^d$. That is, not only is $P$ unbounded, it is the entire space. In this case, it still makes sense to investigate what happens if we consider volume functions subject to the identities of Lemma~\ref{lem:basic}.

\begin{prp}\label{prp:liouville}
	Assume $P=\mathbb{R}^d$, then
	\begin{equation}\label{eq:iddescent}
		0 \   =\    \sum_{ \beta \in \mathbb{Z}^d)^{d+1}}   \vol (\x_{\frac{\alpha+\beta}{2}})^2   \theta^\beta.
	\end{equation}	
\end{prp}	

We obtain at once.

\begin{cor}
	Under the above conditions, $\vol \equiv 0$.
\end{cor}

\begin{proof}[{\textbf{Proof of Proposition~\ref{prp:liouville}}}]
	Fix a $\mu$ in $[d+1]=\{1,\cdots,d+1\}$. Consider the sum
	\[\sum_{\beta \in (\mathbb{Z}^d)^{d+1}\times \{1\}} \vol(x_{\frac{\alpha+\beta}{2}}) \mbf{R}\left[\beta_{|\bar{\mu}},\mu,\frac{\alpha+\beta}{2}-\beta_{|\mu}\right]\]
	where $\beta_{|\mu}$ is the restriction of $\beta$ to the entry $\mu$ and $\beta_{|\bar{\mu}}$ denotes the restriction to all remaining entries.	
	It is easy to see that this equals 
	\[\sum_{ \beta \in (\mathbb{Z}^d)^{d+1}\times \{1\}}   \vol (\x_{\frac{\alpha+\beta}{2}})^2   \theta^\beta\]
	which therefore equals $0.$ The key geometric insight here is that terms of the form $\vol(m)\vol(n)$ appear exactly twice unless $m=n$. 
\end{proof}

We move on

\subsection{The simplex: a (slightly) special linear system}

Now, let us consider the case of the simplex: We have a lattice $d$-simplex $S$ that is a dilation of a unimodular simplex, and let $\alpha$ denote the sum $\sum v$, where $v$ ranges over the vertices of $S\times\{1\}$. For simplicity, we call this the affine barycenter. 

$S$ is an $n$-fold lattice dilation of a unimodular simplex (around the barycenter). Consider instead the $(n+d+1)$st dilation $T$, which contains $S$ as a lattice simplex comfortably in its interior (and the same barycenter). 

Instead of a fully generic linear system, we have the following special linear system: The linear system is fully generic on the lattice points of $S \subset T$, and the vertices of $T$ (that is, generated by algebraically independent variables). We call this the slightly special linear system.

\begin{lem}
	If the Parseval-Rayleigh identities hold for $T$ with respect to the slightly special linear system, then they hold for $S$. 
\end{lem}

\begin{proof}
	This follows at once from Lemma~\ref{lem:bigger} (or more directly from Proposition~\ref{prp:bigger})
\end{proof}

\subsection{The Parseval-Rayleigh identities of the simplex: a subdivision into regions}

Now, let us consider the case of the simplex: We have the lattice $d$-simplex $T$ that is a dilation of a unimodular simplex, and we remind ourselves that we consider this $T$ with the slightly special linear system of parameters whenever we think about the (Artinian reduction) associated semigroup algebra. We will highlight where this peculiarity of the linear system comes in later.

Let $\alpha$ denote the sum $\sum v$, where $v$ ranges over the vertices of $T\times\{1\}$. We also order the vertices $v$ of $T$ from $1$ to $d+1$, so we have an indexing $v_i$.

Before we note how to prove the Parseval-Rayleigh identity of the simplex, let us note what not to do: What is wrong with considering the sum:
\[\sum_{\beta \in (T\cap \mathbb{Z}^d)^{d+1}\times \{1\}} \vol(x_{\frac{\alpha+\beta}{2}})\mbf{R}\left[\beta_{|\bar{\mu}},\mu,\frac{\alpha+\beta}{2}-\beta_{|\mu}\right]?\]
The issue lies in the choice of $\frac{\alpha+\beta}{2}-\beta_{|\mu}$: it may not lie in $\cone^\circ(T)$. And not even in its closure.

Consider a subset $M$ of the vertices $v_i$ of $T\times\{1\}$. Let $T_{|\overline{M}}\times \{1\}$ denote the simplex formed by those vertices of $T\times\{1\}$ not in $M$. Let $\alpha_M$ sum of the vertices of $T_{|\overline{M}}$. The cone over $T_{|\overline{M}}\times \{1\}$ is subdivided into $d+1-\# M$ simplicial cones of dimension $d+1-\# M$: It is exactly consisting of the simplicial cones over $\alpha_M$ over size $d-\# M$ subsets of $(v_i)$ without $M$. This extends to a subdivision of $\R^d\times \R$ into $d+1-\# M$ polyhedra that have the linear span over $\alpha_M$ and the points of $M$ in common. They are indexed by unique vertex of $T_{|\overline{M}}\times \{1\}$ they do not contain.

For any $t \in (\mathbb{Z}^d)\times \{k\}$, and $M$ as above, let $\mbf{L}_{M} (t)$ denote the index of the polyhedron containing $t$ (remember that they are indexed by $S_{|\overline{M}}\times \{1\}$). If there are several such polyhedra, take the least of the vertices (recalling that they are ordered).

Now, consider a $\beta \in (\mathbb{Z}^d)^{d+1} \times  \{1\}$. We associate to it an ordering $\mbf{L}(\beta)$ of the vertices $(v_i)$.

We start with $\mbf{L}_1(\beta):=\mbf{L}_{\emptyset} (|\beta|)$. We set $\mbf{L}_{i+1}(\beta):=\mbf{L}_{\mbf{L}_{i}(\beta)} (|\beta_{|\overline{\mbf{L}_{i}(\beta)}}|)$, where $\beta_{|\overline{\mbf{L}_{i}(\beta)}}$ is $\beta$ with the entries at indices ${\mbf{L}_{i}(\beta)}$ deleted.

\begin{ex}
	Consider the projectivized $T$ to be the convex cone of the vertices $v_1=(0,0,1),$ $v_2=(0,1,1),$ and $v_3=(1,0,1)$ in that order. Consider $\beta \in (\mathbb{Z}^d)^{d+1} \times  \{1\}$ to be the vertices $(1,0,1)$, $(0,1,1)$, and $(0,0,1)$ in that order. Then the associated ordering $\mbf{L}(\beta)$ is $(v_1,v_3,v_2).$
\end{ex}

\subsection{(Co)-admissibility}

We are now considering the sum \begin{multline}\label{eq:summysum}
	\sum_{J \in (T\cap \mathbb{Z}^d)^{d}}\ \sum_{\substack{(z_J,\mu) \in  (\mathbb{Z}^d\times\{{d+1}\})\times [d+1] \\ (z_J,\mu) \text{ admissible } }} \vol(x_{\frac{\alpha+z_J}{2}})\mbf{R}\left[J,\mu,\frac{\alpha-z_J}{2}+|J|\right]\\
	+\  \sum_{J \in (T\cap \mathbb{Z}^d)^{d}}\ \sum_{\substack{(z_J,\mu) \in  (\mathbb{Z}^d\times\{{d+1}\})\times [d+1] \\ (z_J,\mu) \text{ co-admissible } }} \vol(x_{\frac{\alpha+z_J}{2}})\mbf{R}\left[J,\mu,\frac{\alpha-z_J}{2}+|J|\right].
\end{multline}

Here, $(z_J,\mu)$ is \Defn{admissible at $J$} if $z_J-|J|\in \cone(T)$, and if we consider $\hat{J}=J\ast_\mu z_J-|J|$, then the associated $\mbf{L}(\hat{J})$ has $J$ as its first $d$ entries, and $\mu$ is the final index.

Moreover, $(z_J,\mu)$ is \Defn{co-admissible} if $z_J-|J|\notin \cone(T)$, and if we consider $\hat{J}=J\ast_\mu z_J-|J|$, then the associated $\mbf{L}(\hat{J})$ has $J$ as its first $d$ entries, and $\mu$ is the final index, \emph{and} there is a $w$ in $T\times \{1\}$, and $\mu'$ an index in $[d+1]$ such that with respect to $J'=(J\ast_\mu w)_{\overline{\mu'}}$, the pair $(2|J|+2w-z_J, \mu')$ is admissible at $J'$. We call $J'$ and  $(2|J|+2w-z_J, \mu')$ the mirror of $(z_J,\mu)$.

We define these sums whether $\mbf{R}[\ast,\ast,\ast]$ is well-defined or not. This is what we address now.

\subsection{Analyzing (co)-admissibility: The polyhedron $T_J$}

We notice two things: In Equation~\eqref{eq:summysum}, we have $\nicefrac{(\alpha-z_J)}{2}+|J|\in \cone(T)$. This is rather simple: Given a $(z_J,\mu)$ that is admissible with respect to a $J$, verifying that $\nicefrac{(\alpha-z_J)}{2}+|J|\in \cone(T)$ is a simple calculation. If the triple $(J,z_J, \mu)$ is co-admissible, we have to work harder, but we can easily describe a polyhedron containing the value of $z_J-|J|$.
Since $\nicefrac{(\alpha-z_J)}{2}+|J|= \nicefrac{(\alpha+z_J)}{2}-(z_J-|J|)$, this gives the desired. 

Project $J$ to $\pi J\subset \mathbb{R}^d \times \{0\}$ by deleting the last coordinate in each member of the family. Consider the Minkowski sum $T_J:=\pi J + T + (-T)$, where $\pi J$ is thought of as convex hull of its elements. We make two observations:

\begin{lem}
	$z_J-|J|\in T_J\times\{1\}$ if  $(z_J,\mu)$ is co-admissible at $J$.
\end{lem}	
\begin{proof}
	It suffices to consider the projection to a line, that is, the case $d=1$, where the proof is easy.
\end{proof}

The second observation is similarly immediate:

\begin{lem}\label{lem:coad}
	For all $t\in T_J\times\{1\}$, and all $z=|J|+t$, we have $\nicefrac{(\alpha-z)}{2}+|J|\in \cone(T)$. 
	Moreover, $\nicefrac{(\alpha-z)}{2}+|J|\in ´\partial\cone(T)$ if and only if for the mirror $J'$, we have $\nicefrac{(\alpha-2J+2w-z_j)}{2}+|J|'\in \partial\cone(T)$.
\end{lem}

\subsection{The role of the special linear system}
This finishes the argument almost, but leaves some cases where $\nicefrac{(\alpha-z_J)}{2}+|J|\in \partial \cone(S)$, that is, cases in which a summand of \eqref{eq:summysum} is still not well-defined. Let us first consider the admissible triples $(J,z_J, \mu)$: Notice that $\nicefrac{(\alpha-z_J)}{2}+|J|\in \partial \cone(S)$ if and only if $J\ast_\mu (z_J-|J|)$ enumerates the vertices of $T$, because we chose only the defining vertices of $T$ to have nontrivial linear system of parameters. It is not hard to see that these terms exactly correspond to the Kustin-Miller normalization.

As for coadmissible triples, this follows analogously by the second part of Lemma~\ref{lem:coad}. All remaining terms have $\nicefrac{(\alpha-z_J)}{2}+|J|\in \cone^\circ(S)$.\hspace*{\fill} \qed

\section{The Parseval-Rayleigh identities for complexes, 
	revisited}

The following is a slightly different formulation of the arguments contained
in Section~\ref{sec:PC}. It is purely for those who found the derivation of Section~\ref{sec:PC} too quick.

Assume $\fld$ is a field of characteristic $2$ and  $X$ is a lattice ball or sphere
of dimension $d$.
We have that $X$  contains/consists of
a finite number of facets $\{F_1, \dots ,F_s \}$ for some $s \geq 1$,
with each $F_i$ being a lattice polytope.
Without loss of generality we can assume that 
each $F_i$ is a subset  of $\mathbb{R}^d$ and the lattice of $F_i$ is
$\mathbb{Z}^d$.

For each $F_i$, we consider the cone 
\[  
\cone(F_i)  \subset  \mathbb {R}^{d} \times \mathbb {R}.
\] 
Moreover, given $j \geq 0$ we define 
\[
\cone(X)  \cap  (\mathbb {Z}^{d} \times \{j\})
\]
as the "union" (in the obvious sense),  for $1 \leq i \leq s$,  
of the sets    $\; \cone(F_i) \cap  ( \mathbb {Z}^{d} \times \{j\})$.
We will use the notation    
\[
\mathcal{U}^j (X)  =  \cone(X)  \cap  (\mathbb {Z}^{d} \times \{j\}).
\]
For $z \in  \mathcal{U}^j (X)$, we denote by  $ \x_{z}$ the corresponding 
homogeneous  of degree $j$ element of the $\fld$-algebra  $\fld^\ast[X]$.

\begin{rem}   \label{rem!sjakfuerta} 
	Assume $j \geq 0$.   By the definition of the ring  $\fld^\ast[X]$ 
	we have that  the set 
	\[
	\{  \x_{z} :  z \in  \mathcal{U}^j (X) \}
	\]   
	is a basis  of $\fld^j[X]$ as $\fld$-vector space, where 
	$\fld^j[X]$ denotes the $j$-th graded part   of the algebra  $\fld^\ast[X]$.
	Moreover, if $z_1, z_2 \in  \mathcal{U}^j (X)$ satisfy $z_1 \not= z_2$ then
	$\x_{z_1} \not= \x_{z_2}$. 
	
\end{rem}

For simplicity, in the following we assume that $X$ is a $d$-dimensional
lattice sphere with facet set $\{F_1, \dots ,F_s \}$. 

We set  $\widetilde{\fld} = \fld(\theta_{i,j})$.  We  also set 
$V_1 = \widetilde{\fld}^{d+1}[X]$, $V_2 = \AR^{d+1}(X)$ and denote by $\pi : V_1 \to V_2$
the natural projection.
By Remark~\ref{rem!sjakfuerta},  the set 
\[
\{  \x_{z} :  z \in  \mathcal{U}^{d+1} (X) \}
\]   
is a basis  of $V_1$ as $\fld$-vector space.  

Assume $a \in  \mathcal{U}^{d+1} (X) $ and  $b =(v_1, \dots , v_{d+1})$, with each 
$v_i$ a lattice point of $X$.  We set 
\[
\x_b =  \prod_{1 \leq i \leq  d+1 }  \x_{(v_i,1)}  \in  \widetilde{\fld}^{d+1}[X].
\]

We will now define two $\widetilde{\fld}$-linear functions  $f_1, f_2 :  V_1 \to \widetilde{\fld}$. 
Using Remark~\ref{rem!sjakfuerta}, it is enough to specify the value
$f_i (\x_{a}) \in  \widetilde{\fld}$ for all $ a \in  \mathcal{U}^{d+1} (X)$.

We define 
\[
f_1 (\x_{a}) =  (\vol  \circ \pi) (\x_{a}) 
\]
and
\[
f_2 (\x_{a}) =   \sum_{\beta \in \mathcal{U}^{d+1} (X)}  
((\vol  \circ \pi) (\x_{\frac{\alpha+\beta}{2}}))^2   \theta^\beta,   
\]  
where for $\beta = (\beta_1, \dots , \beta_{d+1})$ we set
$\theta^\beta = \prod_{1 \leq i \leq d+1} \theta_{i,\beta_i}$.

\begin{rem}   \label{rem!oldfopsog} 
	By linear algebra, to prove that $f_1 = f_2$, it is enough to prove
	that    
	\[
	\operatorname{Ker} (f_1) \subset \operatorname{Ker} (f_2)
	\]  
	and that there exists  $a \in V_1 \setminus \operatorname{Ker} (f_1) \;$ such that 
	\[ 
	f_2(a) = f_1(a).
	\]
\end{rem}

NOTATION:   Assume $1 \leq i \leq d+1$  and $A  \in   \mathcal{U}^d (X)$. We set 
\[
\theta_i = \sum_{p \in  \mathcal{U}^1 (X)} \theta_{i,p} \x_p  \in  \widetilde{\fld}^{1}[X]
\]
and
\[
q(i,A) =    \x_A  \theta_i  \in V_1.
\]

\begin {prp}    \label{prop!kclsjfjk}
The vector subspace $ \operatorname{Ker} (f_1)$ of $V_1$ is 
generated by the set    
\[
\{ q(i,A) \; :  \;   1 \leq i \leq d+1,  \  A  \in   \mathcal{U}^d (X) \}.
\]
\end{prp}

\begin{proof}  Since $\vol :  \AR^{d+1}(X)  \to  \widetilde{\fld}$ is an isomorphism 
of $\widetilde{\fld}$-vector spaces, we have that 
\[
\operatorname{Ker} (f_1) =  \operatorname{Ker} (\pi).
\]
Since 
\[
\AR^\ast(X) = \widetilde{\fld}^\ast[X] / (\theta_1, \dots , \theta_{d+1}) 
\]
the set in the statement of the proposition generates $\operatorname{Ker} (\pi)$ and
the result follows.  \end{proof}

\begin {prp}   \label{prop!udttgsy}  We have 
\[
\operatorname{Ker} (f_1) \subset \operatorname{Ker} (f_2).
\]  
\end{prp}

\begin{proof}   
Using Proposition~\ref{prop!kclsjfjk}  it is enough to prove  that $f_2(q(i,A)) = 0$ 
for all  $1 \leq i \leq d+1$  and $A \in   \mathcal{U}^d (X)$. 

For simplicity of notation  we set $S =  \mathcal{U}^1 (X)$.

Assume $1 \leq i \leq d+1$  and $A \in   \mathcal{U}^d (X)$.  Since
\[
f_2 (q(i,A)) \  =\  f_2 ( \theta_i \x_A )  \    =\ 
\sum_{p \in S}  \theta_{i,p} f_2 ( \x_p \x_A  )
= \     \ \sum_{p \in S} \theta_{i,p} \sum_{ \beta \in S^{d+1}}  
(\vol (\x_{\frac{A+p+\beta}{2}}))^2   \theta^\beta
\]
we are left with verifying that 
\[  \ \sum_{p \in S} \theta_{i,p} \sum_{ \beta \in S^{d+1}}  
(\vol (\x_{\frac{A+p+\beta}{2}}))^2   \theta^\beta \ =\ 0.   \]

If $\beta_i \not= p $ in the above left hand side we get two equal entries
for the sum, the first for  $p=p$ and $\beta_i = \beta_i$ and the second for 
$p=\beta_i$ and $\beta_i = p$, which they cancel since we work in characteristic $2$.
Consequently,  the left hand side equals 
\[\sum_{p \in S} \theta_{i,p} \sum_{\substack{ \beta \in S^{d+1} \\  \beta_i= p } } 
(\vol (\x_{\frac{A+p+\beta}{2}}))^2   \theta^\beta\ =\ \sum_{p \in S}
\sum_{\substack{ \beta \in S^{d+1} \\  \beta_i= p } }  
(\vol (\x_{\frac{A+p+\beta}{2}}))^2 \theta_{i,p}^2 \prod_{\substack{1 \leq k \leq d+1,\\ k\neq i}} \theta_{k,\beta_k}.\]
We denote by  $B = (B_1, \dots , B_d)$ the sequence obtained from $\beta$ by removing the $i$-th entry.
Then  the right hand side of the last equation equals
\[\sum_{\substack{ B \in S^{d} } }   \prod_{1 \leq k \leq  i-1}  \theta_{k, B_k} 
\prod_{i \leq k \leq  d}  \theta_{k+1, B_k}  \sum_{p \in S}    (\vol (\x_{\frac{A+2p+B}{2}}))^2 \theta_{i,p}^2.\]
But 
\[\sum_{p \in S}  (\vol (\x_{\frac{A+2p+B}{2}}))^2 \theta_{i,p}^2\ =\
\left(\sum_{p \in S}    \vol (\x_{\frac{A+B}{2}} \x_p ) \theta_{i,p}\right)^2\ =\ 0. \]
The reason for the last equality is that 
\[
\sum_{p \in S}    \x_{\frac{A+B}{2}} \x_p  \theta_{i,p} =    \x_{\frac{A+B}{2}}  \theta_{i},  
\]
hence its  class in  $\AR^\ast(X)$ is zero. 
\end{proof}

\begin{lem}\label{lem:parsevaldcomplex}
Assume $X$ is  a lattice $d$-sphere  and $\alpha \in \mathcal{U}^{d+1}(X)$. 
Then  we have in $\AR^\ast(X)$ that
\begin{equation}\label{eq:id}
\vol (\x_\alpha) \   =\    \sum_{ 
	\beta \in (\mathcal{U}^{1}(X))^{d+1}}   \vol (\x_{\frac{\alpha+\beta}{2}})^2   \theta^\beta.
	\end{equation}
\end{lem}

\begin{proof}

Combining locality and  the porcupine construction 
we may also assume that $X$ has one facet that is a unimodular simplex, say $S$: There is 
at least one facet $F$ 
that does not contain $\alpha$. Remove $F$ from $X$ and $\mr{porc}_d F$, and identify the remainders along 
the common boundary $\partial F$.

Combining  Lemma~\ref{prop!udttgsy}  and Remark~\ref{rem!oldfopsog} 
it is enough to prove that the statement is true for some $\alpha$.
For that,  consider the vertices  $v_1, \dots , v_{d+1}$ of the unimodular simplex $S$, 
and set  $\alpha = (v_1, \dots , v_{d+1})$. 
It is clear that the equation  in the statement of the present Lemma is satisfied by $\alpha$. 
\end{proof}

\begin{lem}[The Parseval-Rayleigh identity]\label{lem:parsevalcomplexDOUBLED} 
(Recall we work in characteristic~2)
For a lattice $d$-sphere $X$ 
and $\sigma$ a finite sequence of lattice points in $X$ and for $d+1+\#\sigma$ even, we have
\begin{equation}
\vol (\x_\sigma u^2)\ =\ \sum_{ \beta \in (\mathcal{U}^{1}(X))^{d+1}}
\vol (u\cdot\x_{\frac{\sigma+\beta}{2}})^2 \theta^\beta
\end{equation}
for all $u \in \AR^{\nicefrac{d+1-\#\sigma}{2}}(X)$.  	
\end{lem}

\begin{proof}   We set  $S = (\mathcal{U}^{1}(X))^{d+1}$.
We use Lemma~\ref{lem:parsevaldcomplex}  and argue as following
\begin{align*} 
\vol(\x_\sigma u^2) \ &= \ \vol\left(\x_\sigma\left(\sum_a\lambda_a \x_a\right)^2\right) 
\ = \  \sum_a  \lambda_a^2\ \vol(\x_\sigma\x_a^2) \\
&= \  \sum_a  \lambda_a^2\ \sum_{ \beta \in S}   
\vol(\x_a\cdot\x_{\frac{\sigma+\beta}{2}})^2 \theta^\beta \\
&= \  \sum_{ \beta \in S} \vol\left(\sum_a 
\lambda_a\  \x_a\cdot\x_{\frac{\sigma+\beta}{2}}\right)^2 \theta^\beta \\
&= \ \sum_{ \beta \in S}  \vol (u\cdot\x_{\frac{\sigma+\beta}{2}})^2 \theta^\beta. \qedhere
\end{align*}
\end{proof}

\end{document}